\documentclass[12pt]{amsart}
\usepackage{longtable}
\usepackage{amsfonts,amssymb,amscd,amsmath,epsfig}
\usepackage{latexsym}
\usepackage{mathrsfs}
\usepackage{tocvsec2}
\usepackage[%dvipdfm,
hyperindex=true,pagebackref=true,bookmarks=true,
colorlinks=true,linkcolor=blue,citecolor=red]{hyperref}
\begin{document}

%MACOS FOR LECTURES ON DAHA
\renewcommand{\tilde}{\widetilde}
\renewcommand{\hat}{\widehat}

\newcommand{\BR}{{\mathbb R}}
\newcommand{\BQ}{{\mathbb Q}}
\newcommand{\BC}{{\mathbb C}}
\newcommand{\BP}{{\mathbb P}}
\newcommand{\BZ}{{\mathbb Z}}
\newcommand{\BN}{{\mathbb N}}
\newcommand{\BS}{{\mathbb S}}

\newcommand{\cH}{{\mathcal H}}
\newcommand{\cA}{{\mathcal A}}
\newcommand{\cB}{{\mathcal B}}
\newcommand{\ccF}{{\mathfrak F}}
\newcommand{\cD}{{\mathcal D}}
\newcommand{\cL}{{\mathcal L}}
\newcommand{\cF}{{\mathcal F}}
\newcommand{\cP}{{\mathcal P}}
\newcommand{\cX}{{\mathcal X}}
\newcommand{\cY}{{\mathcal Y}}
\newcommand{\cS}{{\mathcal S}}
\newcommand{\cSol}{\hbox{$\mathcal Sol$}}
\newcommand{\cT}{\hbox{$\mathcal T$}}

\newcommand{\Z}{{\mathbb Z}}
\newcommand{\Q}{{\mathbb Q}}
\newcommand{\N}{{\mathbb N}}
\newcommand{\C}{{\mathbb C}}
\newcommand{\R}{{\mathbb R}}
\newcommand{\X}{{\mathbb X}}
\newcommand{\Y}{{\mathbb Y}}

\newcommand{\CH}{{\mathcal H}}
\newcommand{\CA}{{\mathcal A}}

\def\HH{\mbox{${\mathcal H}$\kern-5.2pt${\mathcal H}$}}

\newcommand{\binomial}[2]{\genfrac{(}{)}{0pt}{}{ #1 }{ #2 }}
\newcommand{\qbinomial}[2]{\genfrac{[}{]}{0pt}{}{ #1 }{ #2 }_q }
\newcommand{\qbinom}[3]{\genfrac{[}{]}{0pt}{}{ #1 }{ #2 }_{ #3 } }

%%SPECIAL SEC 1.0

\def\der{\partial}
\def\tensor{\otimes}
\def\gam{\gamma} \def\Gam{\Gamma}
\def\del{\delta} \def\Del{\Delta}
\def\kap{\kappa}
\def\lam{\lambda} \def\Lam{\Lambda}
\def\Comp{{\mathbb C}}
\def\sM{{\mathcal M}}

\newtheorem{theorem}{Theorem}[section]
\newtheorem{maintheorem}[theorem]{Main Theorem}
\newtheorem{proposition}[theorem]{Proposition}
\newtheorem{definition}[theorem]{Definition}
\newtheorem{lemma}[theorem]{Lemma}
\newtheorem{corollary}[theorem]{Corollary}
\newtheorem{notation}[theorem]{Notation}
\newtheorem{remark}[theorem]{Remark}
\newtheorem{example}[theorem]{Example}

\newtheorem{theorem }{Theorem}[section]
\newtheorem{maintheorem }[theorem]{Main Theorem}
\newtheorem{proposition }[theorem]{Proposition}
\newtheorem{definition }[theorem]{Definition}
\newtheorem{lemma }[theorem]{Lemma}
\newtheorem{corollary }[theorem]{Corollary}
\newtheorem{notation }[theorem]{Notation}
\newtheorem{remark }[theorem]{Remark}
\newtheorem{example }[theorem]{Example}

\newtheorem{ maintheorem }[theorem]{Main Theorem}
\newtheorem{ theorem}{Theorem}[section]
\newtheorem{ proposition}[theorem]{Proposition}
\newtheorem{ definition}[theorem]{Definition}
\newtheorem{ lemma}[theorem]{Lemma}
\newtheorem{ corollary}[theorem]{Corollary}
\newtheorem{ notation}[theorem]{Notation}
\newtheorem{ remark}[theorem]{Remark}
\newtheorem{ example}[theorem]{Example}

\newtheorem{thm}{Theorem}[section]
\newtheorem{prop}[thm]{Proposition}
\newtheorem{lem}[thm]{Lemma}
\newtheorem{cor}[thm]{Corollary}
\newtheorem{conj}[thm]{Conjecture}
\newtheorem{con}[thm]{Conjecture}
\newtheorem{dfn}[thm]{Definition}
\newtheorem{df}[thm]{Definition}
 \newcommand{\rem}{{\bf Comment.\ }}
 \newcommand{\rmk}{{\bf Comment.\ }}
 \newcommand{\exmp}{{\bf Example.\ }}
 \newcommand{\ex}{{\bf Example.\ }}
 \newcommand{\prob}{{\bf Problem.\ }}

\newtheorem{note}{Note} 
\renewcommand{\thenote}{}
\newtheorem*{acka}{Acknowledgments}
\newtheorem{ack}{Acknowledgments}
\renewcommand{\theack}{}
\renewcommand{\appendixname}{\bf Appendix}
\renewcommand{\proof}{{\em Proof.\ }}

\hyphenation{
ap-pen-dix as-ymp-tot-ic at-trib-uted at-trib-ut-able
Bry-li-n-sky com-mu-ta-tion de-ge-ne-rate
de-riv-a-tive dis-trib-ute equi-vari-ant ex-tra-or-di-nary  
geo-met-ric griev-ance griev-ous grad-ed ho-lo-no-my ho-mo-thetic
in-fin-ite-ly in-fin-i-tes-i-mal Ha-rish Cha-n-dra mul-ti-plic-able 
non-euclid-ean non-iso-mor-phic non-smooth par-a-digm 
par-a-bol-ic pa-rab-o-loid pa-ram-e-trize phe-nom-e-non 
post-script pseu-do-dif-fer-en-tial pseu-do-fi-nite 
qua-drat-ics quad-ra-ture Han-kel rec-tan-gle semi-def-i-nite 
set-up wide-spread Euler-ian Feb-ru-ary Gauss-ian Grothen-dieck 
Hamil-ton-ian Her-mi-t-ian her-mi-t-ian Jan-u-ary 
Japan-ese Ka-shi-wa-ra Kor-te-weg Le-gendre No-vem-ber Rie-mann-ian 
Sep-tem-ber Za-mo-lo-d-chi-kov Kni-zh-nik quan-tum Op-dam
Mac-do-nald Ca-lo-ge-ro Su-ther-land Mo-ser 
Ol-sha-net-sky  Pe-re-lo-mov in-de-pen-dent ope-ra-tors 
cy-clo-to-mic ra-tio-nal de-gen-er-a-tion 
in-ter-est-ing de-for-ma-tions de-for-ma-tion pro-ce-dure 
fol-lows ope-ra-tors  pre-serve suf-fices ap-proach 
for-mu-las con-sider its com-ple-tion cor-re-spond-ing 
au-to-mor-phism be-cause pro-por-tional fi-nal-ly let-ting 
equi-v-a-lence ge-n-er-al-ized Mac-do-nald iden-ti-ties 
cor-re-s-pond sub-dia-grams par-ti-tion na-t-u-ral-ly 
or-dered stan-dard de-for-ma-tion ar-gu-ment com-bined 
sphe-r-i-cal rep-re-sen-ta-tions tri-go-no-me-t-ric
ge-n-er-al-ly speak-ing pri-m-it-ive ir-re-du-cible 
sum-ma-tion  rep-re-sen-ta-tives pro-por-ti-o-na-li-ty
ultra-sphe-ri-cal Ro-gers}

\def\ffor{\quad\hbox{ for }\quad}
\def\wwhen{\quad\hbox{ when }\quad}
\def\wwhere{\quad\hbox{ where }\quad}
\def\aand{\quad\hbox{ and }\quad}
\def\for{\  \hbox{ for } \ }
\def\iif{ \ \hbox{ if } \ }
\def\when{ \ \hbox{ when } \ }
\def\where{\  \hbox{ where } \ }
\def\and{\  \hbox{ and } \ }
\def\and{\  \hbox{ and } \ }
\def\oor{\  \hbox{ or } \ }
\def\proof{{\em Proof. \  }}

\def\equal{\stackrel{\,\mathbf{def}}{= \kern-3pt =}}

\def\la{\lambda}
\def\La{\Lambda}
\def\om{\omega}
\def\Om{\Omega}
\def\Th{\Theta}
\def\th{\theta}
\def\al{\alpha}
\def\be{\beta}
\def\ga{\gamma}
\def\ep{\epsilon}
\def\up{\upsilon}
\def\Up{\Upsilon}
\def\de{\delta}
\def\De{\Delta}
\def\ka{\kappa}
\def\kapp{\hbox{\bf \ae}}
\def\si{\sigma}
\def\Si{\Sigma}
\def\Ga{\Gamma}
\def\ze{\zeta}
\def\io{\iota}
\def\bio{b^\iota}
\def\aio{a^\iota}
\def\twio{\tilde{w}^\iota}
\def\hwio{\hat{w}^\iota}
\def\gio{\g^\iota}
\def\Bio{B^\iota}

\def\del{\delta}
\def\pa{\partial}
\def\vp{\varphi}
\def\ve{\varepsilon}
\def\inf{\infty}

\def\vph{\varphi}
\def\vps{\varpsi}
\def\vPh{\varPhi}
\def\vep{\varepsilon}
\def\vpi{{\varpi}}
\def\vth{{\vartheta}}
\def\vsi{{\varsigma}}
\def\vrh{{\varrho}}

\def\bph{\bar{\phi}}
\def\bsi{\bar{\si}}
\def\bvp{\bar{\varphi}}

\newcommand{\bS}{{\mathbf S}}
\newcommand{\bH}{{\mathbf H}}
\newcommand{\bF}{{\mathbf F}}
\newcommand{\bE}{{\mathbf E}}

\def\tal{\tilde{\alpha}}
\def\tbe{\tilde{\beta}}
\def\tde{\tilde{\delta}}
\def\tpi{\tilde{\pi}}
\def\txi{\tilde{\xi}}
\def\tPi{\tilde{\Pi}}
\def\tPhi{\tilde{\Phi}}
\def\tV{\tilde{V}}
\def\tJ{\tilde{J}}
\def\tla{\tilde{\lambda}}
\def\tga{\tilde{\gamma}}
\def\tGa{\tilde{\Gamma}}
\def\tvs{\tilde{{\varsigma}}}
\def\tu{\tilde{u}}
\def\tU{\tilde{U}}
\def\tw{\widetilde w}
\def\tW{\widetilde W}
\def\tB{\tilde B}
\def\tv{\tilde v}
\def\tV{\tilde V}
\def\tz{\tilde z}
\def\tb{\tilde b}
\def\ta{\tilde a}
\def\tih{\tilde h}
\def\trh{\tilde {\rho}}
\def\tx{\tilde x}
\def\tf{\tilde f}
\def\tg{\tilde g}
\def\tG{\tilde G}
\def\tk{\tilde k}
\def\tl{\tilde l}
\def\tL{\tilde L}
\def\tD{\tilde D}
\def\tR{\tilde R}
\def\tP{\tilde P}
\def\tH{\tilde H}
\def\tp{\tilde p}

\def\hH{\hat{H}}
\def\hh{\hat{h}}
\def\hR{\hat{R}}
\def\hY{\hat{Y}}
\def\hX{\hat{X}}
\def\hP{\hat{P}}
\def\hT{\hat{T}}
\def\hV{\hat{V}}
\def\hG{\hat{G}}
\def\hF{\hat{F}}
\def\hw{\widehat{w}}
\def\hW{\widehat{W}}
\def\hu{\hat{u}}
\def\hs{\hat{s}}
\def\hv{\hat{v}}
\def\hb{\hat{b}}
\def\hB{\widehat{B}}
\def\hze{\hat{\zeta}}
\def\hsi{\hat{\sigma}}
\def\hrh{\hat{\rho}}
\def\hth{\hat{\theta}}
\def\hy{\hat{y}}
\def\hx{\hat{x}}
\def\hz{\hat{z}}
\def\hg{\hat{g}}
\def\he{\hat{e}}
\def\hE{\widehat{E}}

\def\B{\mathbf{B}}
\def\I{\mathbf{I}}
\def\P{\mathbf{P}}
\def\G{\mathbf{G}}
\def\S{\mathbf{S}}
\def\F{\mathbf{F}}
\def\one{\mathbf{1}}
\def\Sn{\mathbf{S}_n}
\def\0{\mathbf{0}}
\def\H{\mathbf{H}}
\def\V{\mathbf{V}}

\def\f{\mathcal{F}}
\def\çF{\mathcal{F}}
\def\o{\mathcal{O}}
\def\t{\mathcal{T}}
\def\r{\mathcal{R}}
\def\l{\mathcal{L}}
\def\m{\mathcal{M}}
\def\k{\mathcal{K}}
\def\n{\mathcal{N}}
\def\d{\mathcal{D}}
\def\p{\mathcal{P}}
\def\cP{\mathcal{P}}
\def\a{\mathcal{A}}
\def\h{\mathcal{H}}
\def\c{\mathcal{C}}
\def\y{\mathcal{Y}}
\def\e{\mathcal{E}}
\def\v{\mathcal{V}}
\def\z{\mathcal{Z}}
\def\x{\mathcal{X}}
\def\s{\mathcal{S}}
\def\g{\mathcal{G}}
\def\u{\mathcal{U}}
\def\w{\mathcal{W}}
\def\i{\mathcal{I}}
\def\j{\mathcal{J}}
\def\b{\mathcal{B}}

\def\lan{\langle}
\def\llb{(\!(}
\def\ran{\rangle}
\def\rrb{)\!)}
 \def\dim{{\hbox{\rm dim}}_{\mathbb C}\,}
\def\lng{\hbox{\rm{\tiny lng}}}
\def\sht{\hbox{\rm{\tiny sht}}}
\def\sph{\hbox{\rm{\tiny sph}}}
\def\inv{\hbox{\rm{\tiny inv}}}

\def\br#1{\langle #1 \rangle}

\def\rank{\hbox{rank}}
\def\gl{\mathfrak{gl}_N}
%\def\sgn{\hbox{sgn}}
%\font\germ=eufb10 %at 12pt 
%\def\mathfrak#1{\hbox{\germ #1}}

\newcommand{\Aut}{\operatorname{Aut}}
\newcommand{\Hom}{\operatorname{Hom}}
\newcommand{\End}{\operatorname{End}}
\newcommand{\Ind}{\operatorname{Ind}}
\newcommand{\ad}{\operatorname{ad}}
\newcommand{\pr}{\operatorname{pr}}
\newcommand{\aweyl}{\tilde{\mathbb S}_n}
\newcommand{\hec}{{\mathcal H}^t_n}
\newcommand{\Func}{{\mathcal F}({\mathbb C}^n,{\mathcal H}^t_n)}
\newcommand{\tr}{\operatorname{tr}}
\newcommand{\Out}{\operatorname{Out}}
\newcommand{\Rad}{\operatorname{Rad}}
\newcommand{\Spec}{\operatorname{Spec}}
\newcommand{\id}{\operatorname{id}}
\newcommand{\Int}{\operatorname{Int}}
\newcommand{\ct} {\operatorname{ct}}

\newcommand{\rat}{{\mathbb Q}}
\newcommand{\real}{{\mathbb R}}
\newcommand{\cplx}{{\mathbb C}}
\newcommand{\zint}{{\mathbb Z}}

\newcommand{\sq}{\phantom{1}\hfill$\qed$}
\newcommand{\Rea}{\Re}
\newcommand{\Ima}{\Im}

\newcommand{\st}{\bowtie}
\newcommand{\modd}{\mbox{\,mod\,}}
\newcommand{\lr}{\langle}
\newcommand{\rr}{\rangle}
\newcommand{\eps}{\varepsilon}
\newcommand{\phk}{\phi^{(k)}}
\newcommand{\psk}{\psi^{(k)}}
\newcommand{\Res}{\mbox{Res}\;}
\newcommand{\sgn}{\mbox{sgn}}
\newcommand{\mn} {\left\{ \begin{array}{c}m\\
n\end{array}\right\}}

\def\sX{\mathscr{X}}
\def\sH{\mathscr{H}}
\def\sY{\mathscr{Y}}
\def\TT{\mathfrak{T}}
\def\JJ{\mathfrak{J}}
\def\HH{\mathfrak{H}}
\def\FF{\mathfrak{F}}
\def\GG{\mathfrak{G}}
\def\CC{\mathfrak{C}}
\def\LL{\mathfrak{L}}

\def\BB{\mathfrak{B}}
\def\AA{\mathfrak{A}}
\def\ZZ{\mathfrak{Z}}
\def\HH{\hbox{${\mathcal H}$\kern-5.2pt${\mathcal H}$}}
\def\HHH{\hbox{${\mathbb H}$\kern-4.2pt${\mathbb H}$}}
\def\tHH{\widetilde{\HH\ }}

\font\smm=msbm10 at 12pt 
\def\symbol#1{\hbox{\smm #1}}
\def\lsmash{{\symbol n}}
\def\rsmash{{\symbol o}}
\def\#{\sharp}

\font\tenbf=cmbx10
\font\tenrm=cmr10
\font\tenit=cmti10
\font\ninebf=cmbx9
\font\ninerm=cmr9
\font\nineit=cmti9
\font\eightbf=cmbx8
\font\eightrm=cmr8
\font\eightit=cmti8
\font\sevenrm=cmr7
\font\sevenbf=cmbx7

%END MACROS

\title [Modules over plane curve singularities]
{Modules over plane curve singularities in any
ranks and DAHA}
\author[Ivan Cherednik]{Ivan Cherednik $^\dag$}
\author[Ian Philipp] {Ian Philipp}
%\date{February 2, 2014}

\begin{abstract}
We generalize the construction of geometric superpolynomials
for unibranch plane curve singularities from our prior
paper from rank one to any ranks; explicit formulas are
obtained for torus knots. The new feature is the
definition of counterparts of Jacobian factors (directly
related to compactified Jacobians) for higher ranks, which
is parallel to the classical passage from invertible
sheaves to vector bundles over algebraic curves. This is
an entirely local theory, connected with affine
Springer fibers for non-reduced (germs of) spectral curves.
We conjecture and justify numerically the connection of our
geometric polynomials in arbitrary ranks with the corresponding
DAHA superpolynomials for any algebraic knots colored by columns.
\end{abstract}

\thanks{$^\dag$ \today.
\ \ \ Partially supported by NSF grant
DMS--1363138}

\address[I. Cherednik]{Department of Mathematics, UNC
Chapel Hill, North Carolina 27599, USA\\
chered@email.unc.edu}

\address[I. Philipp] {Department of Mathematics, UNC
Chapel Hill, North Carolina 27599, USA\\
iphilipp@live.unc.edu}

 \def\sht{\raisebox{0.4ex}{\hbox{\rm{\tiny sht}}}}
 \def\bysame{{\bf --- }}
 \def\~{{\bf --}}

%%%rk>1 PAPER ONLY:
\def\rk{r\!k} % \def\rk{{\mathsf \varrho}}
\renewcommand{\dim}{\hbox{dim}_k}   %%%OTHERWISE dim_C
\renewcommand{\deg}{\hbox{deg}_{\o}}   %%%OTHERWISE dim_C

 \def\rr{{\mathsf r}}
 \def\cc{{\mathsf c}}
 \def\ss{{\mathsf s}}
 \def\mm{{\mathsf m}}
 \def\pp{{\mathsf p}}
 \def\ll{{\mathsf l}}
 \def\aa{{\mathsf a}}
 \def\bb{{\mathsf b}}
 \def\NS{\hbox{\tiny\sf ns}}
 \def\ssum{\hbox{\small$\sum$}}
\newcommand{\comment}[1]{}
\renewcommand{\tilde}{\widetilde}
\renewcommand{\hat}{\widehat}
\renewcommand{\V}{\mathbb{V}}
\renewcommand{\F}{\mathbb{F}}
\newcommand{\dagx}{\hbox{\tiny\mathversion{bold}$\dag$}}
\newcommand{\ddagx}{\hbox{\tiny\mathversion{bold}$\ddag$}}
\newtheorem{conjecture}[theorem]{Conjecture}
\newcommand*\toeq{
\raisebox{-0.15 em}{\,\ensuremath{
\xrightarrow{\raisebox{-0.3 em}{\ensuremath{\sim}}}}\,}
}
\newcommand{\unknot}{\hbox{\tiny\!\raisebox{0.2 em}{$\bigcirc$}}}
%\WarningFilter{latex}{Text page} %%%MAKE IT % FOR ARXIV
\newcommand*{\vect}[1]{\overrightarrow{\mkern0mu#1}}
\newcommand*{\medcap}{\mathbin{\scalebox{.75}
{\ensuremath{\bigcap}}}}
\newcommand*{\medcup}{\mathbin{\scalebox{.75}
{\ensuremath{\bigcup}}}}
%\WarningFilter{latex}{Text page} %%%MAKE IT % FOR ARXIV

\vskip -0.0cm
%\par
%{\centering
%\medskip
%\par}
%\vskip -0.0cm
\maketitle
\vskip -0.0cm
%\smallskip


{\tiny
\centerline{{\bf MSC} (2010): 14H50, 17B22, 17B45, 20C08, 20F36,
22E50, 22E57, 30F10, 33D52, 33D80, 57M25.}
}
\smallskip

\vskip -0.5cm
\renewcommand{\baselinestretch}{0.95}
{\small
%%%\pagenumbering{gobble}
\tableofcontents   %%% REMOVE % IN ARXIVE!!!!
}
\renewcommand{\baselinestretch}{1.0}

\renewcommand{\natural}{\wr}

\vfill\eject
\setcounter{section}{-1}
\setcounter{equation}{0}
\section{\sc Introduction}

We generalize the {\em geometric
superpolynomials\,}
for unibranch plane curve singularities from \cite{ChP}
from rank one to any ranks; explicit formulas are proved
for torus knots. The key is the passage from
the Jacobian factors (local
factors of compactified Jacobians) in rank one
to higher ranks, which
is parallel to the classical passage from the Jacobians
to vector bundles over smooth algebraic curves. 
This is a local theory, closely connected with
affine Springer fibers with {\em non-reduced\,} spectral
curves and related to the {\em Riemann Hypothesis\,} from
\cite{ChRH} when the Young diagrams are columns.
% (their germs, to be precise).

\subsection{\bf Flagged Jacobian factors}
{\em Flagged Jacobian factors\,} 
$\j_{\r,\rk}^{\ell}$ of rank $\rk$ are defined in 
(\ref{jrk-def}), are certain {\em 
increasing\,} full $\ell$\~flags of {\em standard\,} modules 
$\m\subset\o^{\rk}$, where $\r$ is the ring of singularity, $\o$ its 
normalization ring; standard modules are those generating 
$\o^{\rk}$ over $\o$. When $\ell\!=\!0$ the dimension of  
$\dim \j_{\r,\rk}^{\ell=0}=\rk^2\de$ for the arithmetic genus 
$\de$ of the singularity, which agrees with the classical formula  
$\rk^2(g\!-\!1)\!+\!1$ for smooth projective curves
of genus $g\!>\!1$. To match the
DAHA superpolynomials we do not 
divide by the action of (constant) $GL_{\rk}$ in $\o^{\rk}$,
which enlarges the dimension by  $\rk^2\!-\!1$ versus the
classical formula. 

The geometric superpolynomials  conjecturally
coincide with DAHA  superpolynomials of the corresponding
algebraic knots colored by $\rk$\~columns (columns of hight $\rk$)
and are given directly in terms
of  $\j_{\r,\rk}^{\ell}$. The latter are
quasi-projective varieties for $\rk\!=\!1$ and 
constructible sets for $\rk\!>\!1$. Their projectivization is
natural, but the boundary must be omitted
to ensure the connection with the DAHA-superpolynomials.
Presumably the boundary can be interpreted via standard modules
of lower ranks and is also related to DAHA
(to be discussed somewhere).
\vskip 0.2cm

The role of vector bundles and $Bun_G$ is tremendous in
modern mathematics and physics. 
Our approach can be
potentially extended to any reductive $G$.
% but we stick
%to $G=GL_{\rk}$ in
%this paper and fix one closed fiber up to proportionality
%(do not ``divide" by the action of $G$) in our counterpart 
%of $Bun_G$. 
The spaces $Bun_G$ are generally 
stacks and the  {\em stability\,} is of great importance 
for their theory. Our $\j_{\r,\rk}^{\ell}$ are constructible sets
for $\rk\!>\!1$
(projective varieties minus some explicit quasi-projective 
subvarieties). The {\em liftability to standard modules\,} is a 
certain analog
of (semi)stability, but its exact algebro-geometric
meaning is not clarified at the moment. Due to the
local setting (local singularities instead of curves), our
$\j_{\r,\rk}^{\ell}$ are of a much more combinatorial nature than
the classical $M(\rk,d)$ classifying
(stable) vector bundles of rank $\rk$ and degree $d$
over smooth projective curves. 
It is even not impossible that they are always
paved by affine spaces for $G=GL$, but we do not
conjecture this. There is some parallelism with
the theory of Schubert varieties in the
classical and affine settings, though  $\j_{\r,\rk}^{\ell}$
(and the combinatorics involved)
are of course more subtle due to the $\r$\~invariance.
%\vskip 0.2cm

\subsection{\bf Main Conjecture}
Conjecture \ref{CONJ:MAIN}
states that the geometric superpolynomials
$\h^{sin\!g}_{\r,\rk}(q,t,a)$ and $\h^{mot}_{\r,\rk}(q,t,a)$,
%of the liftable part of $\j_{\r,\rk}^{\ell}$,
describing flags of {\em standard\,}
rank-$\rk$ modules, coincide
up to some normalization factors ($q,t$\~powers)
with each other and with the DAHA
superpolynomials of the corresponding algebraic knots
colored by $\rk$-column diagrams.
%\smallskip

\vfil The {\em singular\,} geometric superpolynomials 
$\h^{sin\!g}_{\r,\rk}(q,t,a)$ are defined in (\ref{conjcoh}) in 
terms of the Borel-Moore (relative) homology $H^{\!BM}_{2i}$, of 
{\em quasi-projective varieties\,} $\j_{\r, \rk}^{\ell}[dev]$ 
formed by modules $\m$ of fixed {\em deviation\,} $dev$, which 
is the lattice index $[\m:\r^{\rk}]$ equal to 
$\hbox{dim}(\o^{\rk}/\r^{\rk}) - \hbox{dim}(\o^{\rk}/\m)$;
$\r$ is the ring of a given plane curve singularity. Up to 
a power of $t$, they are the sums of the 
corresponding $\,q^{dev\,} t^{-i}a^\ell.$  
\vskip 0.2cm

The {\em motivic\,} geometric superpolynomial
 $\h^{mot}_{\r,\rk}(q,t,a)$
from (\ref{conjmot}) is
a variant of $\h^{sin\!g}_{\r,\rk}(q,t,a)$, where
instead of employing homology, we count points of
$\j_{\r,\rk}^{\ell}$ over finite fields $\mathbb F_{p^m}$
with the weights $q^{dev} a^\ell$. We then set
$t=1/p^m$, which coincides with the previous definition
if $\j_{\r,\rk}^{\ell}$ is paved by affine spaces.
\vskip 0.2cm

The third superpolynomial is the DAHA one
 $\h^{daha}_{\l}(\om_{\rk}; q,t,a)$,
where $\l$ is the link of the corresponding
singularity up to isotopy.
To be exact, the construction of $\h^{daha}_{\l}$
is entirely algebraic; it is in terms of certain pairs of
relatively prime numbers. When these pairs are interpreted
as Newton's pairs of plane curve singularities, the
theorems from \cite{CJJ,ChD1,ChD2} state that the corresponding
DAHA superpolynomials actually depend only on $\l$, which 
is a non-trivial statement even for torus knots (though
not too difficult to prove), i.e. they are
{\em topological invariants\,} of plane
curve singularities.
\vskip 0.1cm

By construction, our geometric
superpolynomials, ORS polynomials and orbital integral
in the theory of affine Springer fibers
depend on the {\em analytic\,}
type of singularities. 
%We see no {\em a priori\,}
%reasons for these constructions to depend only on
%the topological type of the corresponding link
%(in type $A$).
%This implication of our conjecture is of obvious importance.
The analytic classification of plane curve singularities
is much more involved than their topological classification.
The justification of
our Main Conjecture is within reach (but not too simple).
One needs to compare the
recurrence relations for the DAHA superpolynomials (following
\cite{CJ}) with the behavior
of the geometric superpolynomials under monoidal transformations.
This open a road to verifying the Riemann Hypothesis from
\cite{ChRH} in the case of columns.

%This would give the topological 
%invariance of geometric superpolynomials.

%\vskip -0.5cm

\subsection{\bf Connections and applications}
The singular superpolynomials and (conjecturally)
those for DAHA
%via $\h^{sing}_{\r,\rk}(q,t,a)=\h^{daha}_{\l}(\om_{\rk}; q,t,a)$,
provide a lot of geometric information about $\r$ and 
$\j_{\r,\rk}^{\ell}$. For instance, their Betti numbers can be
obtained as $a=0,q=1$ (when $\ell=0$), the Alexander polynomial
occurs when
$a=-1, t=q$, the HOMFLY-PT polynomial of the corresponding
link $\l$ is for
$a\mapsto -a, t=q$ and so on. Concerning the motivic
version, the specialization $a=0,q=1,t=1/p^m$ gives the
orbital integrals in the nil-elliptic case (for $G=GL_{\rk}$).     
The whole polynomial
$\h^{daha}_{\l}(\om_{1}; q,t,a)$ (for $\rk=1$)
is conjectured to coincide
with the stable {\em Khovanov-Rozansky polynomial\,} of $\l$,
uncolored and reduced, upon some recalculation of
parameters and normalization. See  \cite{Kh,KhR1, Ras,
Rou, Ma, ChD2,ChP}.
%\smallskip
%\vskip 0.2cm

We will not discuss these and other {\em connection
conjectures\,} in this paper, as well as the general
DAHA-Jones theory. Concerning the latter,
see \cite{CJ,GoN,CJJ,ChD1};
the most comprehensive work at this point is \cite{ChD2}, the case
of arbitrary colored iterated torus links. Note that we
consider only algebraic {\em knots\,} colored by columns.
Algebraic {\em links\,} are doable too; a challenge is to
extend this and our previous paper \cite{ChP} to any
rectangle diagrams for arbitrary algebraic links. The
{\em positivity\,} for such links
conjectured in \cite{ChD1,ChD2} matches very well
our geometric approach (and the $RH$ Conjecture in \cite{ChRH}).
%\vskip 0.2cm

See \cite{ORS} and \cite{GORS} for the {\em ORS polynomials\,};
we compare them at the end of \cite{ChP}
with our geometric (singular and motivic) superpolynomials,
though these constructions are different. See also \cite{ChRH}
about connections with Kapranov's zeta function and beyond.
Some origins of our geometric superpolynomials, including
\cite{LS,Pi,GM2}, are discussed in \cite{ChP}.
Papers \cite{ObS, Ma, MoS} are about the important particular
case $t=q$  (HOMFLY-PT polynomials). For physics approaches, see
\cite{DGR,AS,GS,DMMSS,FGS,GGS}; \cite{AS} triggered \cite{CJ},
the super-duality \cite{GS} generally results from $M_5$\~theory. 
%\smallskip

There are also connections with the Hilbert schemes of $\C^2$;
see e.g. \cite{GoN,GORS}. Rational DAHA are related
to torus knots \cite{Gor2, EGL}. There are some recent
developments with Soergel modules \cite{Mel}. Let us also
mention the theory of compactified Jacobians with direct
links to Fundamental Lemma; see  \cite{MY,MS}.
\vfil
\eject

\setcounter{equation}{0}
\section{\sc DAHA superpolynomials}
In this section we will briefly summarize the necessary
theory for the reader to understand the definition
of DAHA superpolynomials for $GL_n$. Note that
DAHA-Jones polynomials are well defined for all
simple root systems and weights (see \cite{ChD1}).
The theory of DAHA superpolynomials is so far for type $A$.
Though there are some developments (mostly conjectures)
for classical root systems and even a couple of
examples in type $E$;  see \cite{CJ, ChP}.

\subsection{\bf DAHA of type
\texorpdfstring{{\mathversion{bold}$GL_n$}}
{GL}}
Let us begin with the definition of DAHA itself. See \cite{C101}
for details (and justifications).

\begin{definition}\label{daha}
The double affine Hecke algebra $\HH$ (referred to as DAHA) of
$GL_n$ is an algebra over $\C_{q,t}\equal\C(q^{1/2},t^{1/2})$
generated by $T_i^{\pm 1},X_j^{\pm 1}$ and $\pi$,
where $1\leq i \leq n-1$, $1\leq j \leq n$, which are
subject to the following relations:

(i)\,\, $(T_i-t^{1/2})(T_i+t^{-1/2})=0$, \quad
$T_iT_{i+1}T_i=T_{i+1}T_iT_{i+1}$,

(ii)\, $T_iT_k=T_kT_i$ \hspace{1mm} \mbox{if}\quad $|i-k|>1$,

(iii) $X_jX_k=X_kX_j$ \hspace{1mm} $(1\leq j,k\leq n)$,
\quad $T_iX_iT_i=X_{i+1}$ \hspace{1mm} $(i<n)$,

(iv)\, $\pi X_i=X_{i+1}\pi$ $(1\!\leq\! i\!\leq\! n\!-\!1)$
\mbox{ and }
$\pi^n X_i=q^{-1}X_i \pi^n$ $(1\!\leq\! i\! \leq\! n)$,

(v)\,\, $\pi T_i=T_{i+1}\pi$ $(1\leq\! i\!\leq\! n-2)$
\mbox{ and }
$\pi^n T_i\,=T_i \pi^n$ $(1\leq\! i\! \leq n\!-\!1)$.

\vskip -0.5cm
\sq
\end{definition}
The subalgebra of $\HH$ generated by
$T_i$ for  $1\leq i\leq n-1$
is isomorphic to the (non-affine) Hecke algebra of
type $GL_n$. The following {\em self-duality\,} is the
key in DAHA theory. Let us define pairwise commutative
$Y_i$ $(1\leq i\leq n)$ inductively as follows:
\[
Y_1=\pi T_{n-1}\cdots T_1 \mbox{ and }
Y_{i+1}=T^{-1}_iY_iT^{-1}_i.
\]
The subalgebras $\h_X$ and $\h_Y$ generated by
 $T_i,X_j$ and $T_i,Y_j$, where
$1\leq i\leq n-1, 1\leq j\leq n$,
are isomorphic to the extended
affine Hecke algebra of type $GL_n$ (and to each other) and
the aforementioned duality is given by an
 anti-involution of $\HH$, sending
$$
\varphi:X_j\leftrightarrow Y^{-1}_j (1\!\le\! j\!\le\! n),\,\,
 T_i \mapsto T_i
(1\!\le\! i\!\ <\! n) \hbox{\,\, and fixing \,\, }
q^{1/2},t^{1/2}.
$$

{\sf Projective $PSL_2(\Z)$}. Along with $\varphi$,
the symmetries of $\HH$ include
an action of the braid group $B_3$ on three strands.
Its two generators act as follows:
\begin{align}\label{tau+tau-}
&\tau_+(X_i)=X_i, \quad \tau_+(Y_1\cdots Y_i)\,=
q^{-i/2}(X_1\cdots X_i)(Y_1\cdots Y_i),\\
&\tau_-(Y_i)=Y_i,\, \quad \tau_-(X_1\cdots X_i)=\,
q^{i/2}\,(Y_1\cdots Y_i)\,(X_1\cdots X_i),\notag
\end{align}
where $\tau_+,\tau_-$ fix $T_i$ and $q^{1/2},t^{1/2}$.
They satisfy
$\tau_+\tau_-^{-1}\tau_+$ $=\tau_{-}^{-1}\tau_+\tau_{-}^{-1}$,
the defining relations of $B_3$.
This group is naturally a central extension
of $PSL_2(\Z)$; the covering map is
\begin{align}\label{lift}
\tau_+\mapsto \begin{pmatrix}
1 & 1\\
0 & 1
\end{pmatrix},\quad \tau_-\mapsto \begin{pmatrix}
1 & 0\\
1 & 1
\end{pmatrix}.
\end{align}
We will call it {\em projective $PSL_2(\Z)$\,};
this is due to Steinberg. See
\cite{C101}, Section 3.7 for details.
We note that $\tau_+$ is the conjugation
of $Y_i$ by the {\em Gaussian\,} $\ga=q^{\sum_{i=1}^n
x_i^2/2}$, where we put formally $X_i=q^{x_i}$.

\subsubsection{\sf Polynomial representation}
We use here and below the standard basis
$\{\ep_j, 1\leq j\leq n\}$ in $\R^n$, supplied with
the inner product $(\ep_i,\ep_j)=\de_{ij}$.
Let $\al_i=\ep_i-\ep_{i-1}$ be
simple roots for $\mathfrak{gl}_n$,
$\{\al\!=\!\ep_i\!-\!\ep_j,\, i\!<\!j\}$ the set of
positive roots,\,  $\rho=(1/2)\sum_{\al>0}\al$,\,
$W\!=\!\mathbb{S}^n$ the symmetric group,
$s_i=(i,i+1)$ simple reflections. The fundamental
roots are
$\omega_i=\ep_1+\ep_2+...+\ep_i$. Note that
$(\rho,\ep_i)=(n-(2i-1))/2,\ $
$(\rho,\om_i)=i(n-i)/2,\ (\rho,\ep_i-\ep_j)=j-i.$
The lattices
$Q,P$ are correspondingly
formed by  integer linear combinations of $\{\al\},\{\om_i\}$.
Dominant weights are
$b=\sum b_i\omega_i$ $b_i\in \Z_+$; the notation is
$b\in P_+$.
\vskip 0.2cm

{\sf Affine Weyl group.}
Let $\tilde{R}=\{[\al,j], \al\in R,j\in \Z\}$ be
the affine root system of type $A$. Affine
positive roots are for $j>0$ or with $\al>0$ if $j=0$.
The {\em extended Weyl group\,} $ \hW$ is the semidirect
product $W\lsmash P$, where
the action in $\C^{n+1}\ni [z,\ze]$ is:
\begin{align}
&(wb)([z,\ze])\ =\ [w(z),\ze-(z,b)] \for w\in W, b\in P.
\label{ondthr}
\end{align}
It contains $\tW\equal \lan s_i, 0\le i\le n-1\ran$ where
$\al_0=[1,-\th]$, where $\th=\ep_1-\ep_n$ is
the maximal positive root.
$\hW$ is naturally isomorphic to
$\tW\rsmash \Pi$ for $\Pi\equal P/Q$.
Setting $\hw = \pi^m\tw \in \hW$ for $m\in \Z,\,
\tw\in \tW$,  the length
\,$\ell(\hw)$ of $\hw$ is the length of any reduced decomposition
of $\tw$ in terms of simple reflections
$s_i (0\!\le\! i\!\le\! n).$ Then $\Pi$ is the subgroup of
$\hW$ of the elements of length $0$.
\vskip 0.2cm

{\sf Coinvariant.}
For $\hw=\pi^m\tw\in \hW$ and a reduced decomposition
$\tw=s_{i_\ell}\cdots s_{i_1}$ of length $\ell=\ell(\hw)$, we set
$T_{\hw}=\pi^m T_{i_\ell}\cdots T _{i_1}$.
There is a unique presentation of any
$H\in \HH$ as a sum of $(\prod X_j^{u_j})T_w
(\prod Y_j^{v_j})$,
where $u_j,v_j\in \Z, w\in W$; this is a DAHA counterpart of the
PBW Theorem.
We use this  to define
the \emph{coinvariant}
\begin{align}\label{coinvdef}
\{\cdot\}_{ev}:\HH\to \C_{q,t},\,
\{X_j\}_{ev}\!=\! t^{-(\rho,\ep_j)},
\{Y_j\}_{ev}\!=\! t^{(\rho,\ep_j)},
\{T_i\}_{ev}\!=\! t^{1/2}.
\end{align}
By construction:
$\{\varphi(H)\}_{ev}=\{H\}_{ev}$ for $H\in \HH$.

The coinvariant factors through the \emph{polynomial
representation\,} $\mathscr{X}$ of $\HH$.
Let $\chi$ be the following
one-dimensional character of the affine Hecke algebra $\h_Y$:
$\chi(T_i)=\{T_i\}_{ev}=t^{1/2}$,
$\chi(Y_j)=\{Y_j\}_{ev}$. Then
$\mathscr{X}\equal\text{Ind}_{\h_Y}^{\HH}(\chi)$.
We will use the following projection of $\HH$
onto $\mathscr{X}$:
\begin{align}\label{downproj}
\HH\ni H\mapsto H\!\Downarrow\,\equal\,
\hat{H}(1)\in \mathscr{X},
\end{align}
where $\hat{H}$ is the image
of $H$ in $\text{End}(\mathscr{X})$. The {\em coinvariant\,}
is then $\{H\}_{ev}=\{H\!\!\Downarrow\}_{ev}$, where
the latter evaluation sends
$X_j\mapsto \{X_j\}_{ev}$.

\subsubsection{\sf Macdonald Polynomials}
The polynomial representation is given explicitly by
the {\em Demazure-Lusztig operators\,}:
\begin{align}
&T_i\  = \  t^{1/2} s_i\ +\
(t^{1/2}-t^{-1/2})(X_{\al_i}-1)^{-1}(s_i-1),
\ 0\le i\le n.
\label{Demazx}
\end{align}
The elements $X_b$ become the multiplication operators
in $\mathscr{X}$;  $\pi^m\in \Pi$ act via the
general formula
$\hw(X_b)=X_{\hw(b)}$ for $\hw\in \hW$,
where $X_{[b,j]}\equal q^j X_b$.
The starting point of
the DAHA theory was adding $T_0$ here for
$\al_0=[1,-\th]$, where $\th=\ep_1-\ep_n$ is
the maximal positive root.
\vskip 0.2cm

{\sf The $E$ and $P$\~polynomials.}
The \emph{nonsymmetric Macdonald polynomials} $E_b$
for $b\in P$
are eigenvectors of $\hat{Y_j}$ (the images of $Y_j$ in
$\text{End}(\mathscr{X})$). This definition fixes
them uniquely up to
proportionality {\em for generic $q,t$\,}; the normalization
is to make  the coefficient of
$X_b$ equal to $1$.
The {\em symmetric Macdonald polynomials\,}
$P_b$ are defined
for {\em dominant\,} $b=\sum b_i\omega_i$, $b_i\in\Z_+$.
In terms of $E_b$, one has:
$P_b={\bf P}E_b$, where
\[
{\bf{P}}\equal\bigl(\sum_{w\in W}t^{\ell(w)}\bigr)^{-1}
\sum_{w\in W}t^{\ell(w)/2}T_w
\hbox{\,\, is the $t$-symmetrizer}.
\]
See \cite{C101}, (3.3.14).
The coefficient of $X_b$ is $1$ in $P_b$ as well.

Actually, the fundamental weights $\,\om_i\,$ are the key
for us in this paper. Then $E_{\om_i}=X_{\om_i}$ and
$P_{\om_i}= \sum_{b\in W(\om_i)} X_{b}$ for
$1\le i\le n$. For instance, $P_{\om_1}=P_{\ep_1}=
\sum_{i=1}^n X_i.$

The following
evaluation formula (the Macdonald
Evaluation Conjecture) is important for DAHA superpolynomials:
\begin{align}
\label{macdeval}
&P_{b}(t^{-\rho})=
t^{-(\rho,b)}
\prod_{\al>0}\prod_{j=0}^{(\al,b)-1}
\Bigl(
\frac{
1- q^{j}t X_\al(t^{\rho})
 }{
1- q^{j}X_\al(t^{\rho})
}
\Bigr),
\end{align}
where $X_{[b,j]}(t^{c})\equal q^j t^{(b,c)}$.
The evaluation can be here $P_b(t^{\pm\rho})$, but
the corresponding evaluation formulas for $E_b$ must be
exactly at $t^{-\rho}$; see (3.3.16) from \cite{C101}.

\subsection{\bf DAHA-Jones theory}
\subsubsection{\sf Algebraic knots}
In this section, we briefly review the necessary aspects of
algebraic knots. See e.g. \cite{EN} for details.

Let $\c$ be the germ of a unibranch plane curve singularity
at $0\in \C^2$. The intersection of $\c$
with a small $\S^3$ centered at $0$ is the corresponding
\emph{algebraic knot} $\l=\l_\c$, called the link of $\c$.
 All algebraic knots are
\emph{iterated torus knots}, but generally
iterated torus knots are not algebraic.

The Newton-Puiseux presentation of $\c$
is as follows:
\begin{align}\label{yxcurve}
y = x^{\frac{\ss_1}{\rr_1}}
(c_1 + x^{\frac{\ss_2}{\rr_1\rr_2}}
\bigl(c_2 + \ldots +
x^{\frac{\ss_l}{\rr_1\rr_2\cdots\rr_l}}\bigr))
\hbox{\, centered at\, } 0.
\end{align}
The pairs
$\{\rr_i,\ss_i\}$ are known as \emph{Newton pairs}.
They are positive integers such that
$\text{gcd}(\rr,\ss)=1$; $c_j$ here are
assumed sufficiently general.

The \emph{torus knot} $T(\rr,\ss)$ is for one pair
$\{\rr,\ss\}$. Topologically, it lies on
the surface of a torus and wraps $\rr$ times around the torus
longitudinally and $\ss$ times around latitudinally. Due to the
isotopy $T(\rr,\ss)=T(\ss,\rr)$, we may assume $\rr>\ss$.
For any $l$, let
\begin{align}\label{iterrss}
&\vec{\rr}=\{\rr_1,\ldots \rr_{l}\},\,\
\vec{\ss}=\{\ss_1,\ldots \ss_{l}\},\\
&\aa_1=\ss_1,\,\aa_{i}=\aa_{i-1}\rr_{i-1}\rr_{i}+\ss_{i}\,\
(i=2,\ldots,l). \notag%\label{Newtonpair}
\end{align}
The pairs $\{\aa_j,\rr_j\}$ are topological invariants
of $\l_\c$, the Newton pairs are not.

The knot $\l_\c$ can be described topologically as follows:
\begin{align}\label{Knotsiter}
\t(\vec{\rr},\!\vec{\ss})\equal
C\!ab(\vec{\aa},\!\vec{\rr})(\unknot)=
\bigl(C\!ab(\aa_l,\rr_l)\cdots C\!ab(\aa_2,\rr_2)\bigr)
\bigl(T(\rr_1,\!\ss_1)\bigr).
\end{align}
Here $\unknot$ is unknot and
$C\!ab$ is the {\em cabling\,} defined as follows.
For any oriented knot $K$ embedded in the 3-sphere consider
a small torus around $K$ which does not intersect itself.
Then $C\!ab(\aa,\bb)(K)$  is the knot created by placing
$T(\aa,\bb)$ on this torus with a proper framing. We omit 
details, but let us
remark that we can assume now that $\rr_1>\ss_1$ only
for the first pair of indices; see e.g. \cite{EN,ChD1}.

\subsubsection{\sf Main construction}
Now we will review DAHA-Jones polynomials
 and DAHA superpolynomials.
The main source for this content is \cite{ChD1}.

Generally, for any polynomial $F$ in fractional powers of
$q$ and $t$, its
{\em tilde-normalization} $\tilde{F}$ is the result of $F$ divided 
by the lowest $q, t$-monomial, assuming it is well defined.
We will use the notation $q^\bullet t^\bullet$ for
a monomial factor (possibly with fractional exponents)
in terms of $q,t$. We will also use the tilde-normalization for
superpolynomials, where the lowest $q,t$-monomial is chosen
from the those terms without the third parameter $a$ (i.e.
at $a=0$).

The following theorem is from \cite{ChD1}; it is stated there
for arbitrary root systems (see also references
there).  Recall (\ref{coinvdef})
and (\ref{downproj}), the definitions of the coinvariant
and the projection onto $\mathscr{X}$. 

\begin{theorem}\label{JONITER}
(i) Given two strictly
positive sequences $\vec\rr,\vec\ss$ of length $l$ as
in (\ref{iterrss}), let $\ga_i\in PSL_2(\Z)$ be the element whose
representative is any  matrix with the first column
$(\rr_i,\ss_i)^{tr}$ ($tr$ indicates transpose).
Let $\hat{\ga_i}$ be a lift of $\gamma_i$ to the
projective $PSL_2(\Z)$ ($=B_3$) via
the homomorphism (\ref{lift}).
For a dominant weight $b$ of $GL_n$,
the {\sf DAHA-Jones polynomial} is
\begin{align}\label{jones-dit}
&J\!D_{\,\vec\rr,\,\vec\ss\,}^{GL_n}(b\,;\,q,t)\ =\
J\!D_{\,\vec\rr,\,\vec\ss\,}(b\,;\,q,t)\, \equal\\
\Bigl\{\hat{\ga_{1}}&\Bigr(
\cdots\Bigl(\hat{\ga}_{l-1}
\Bigl(\bigl(\hat{\ga}_l(P_b)/P_b(t^{-\rho})\bigr)\!\Downarrow
\Bigr)\!\Downarrow\Bigr) \cdots\Bigr)\Bigr\}_{ev}.\notag
\end{align}
It does not depend
on the particular choice of the lifts \,$\ga_i$
and $\hat{\ga}_i$.
The {\sf tilde-normalization}
$\tilde{J\!D}_{\vec\rr,\vec\ss}\,(b\,;\,q,t)$ is well defined
and is a polynomial in terms of $\,q,t$ with the
constant term $1$.

(ii) Now consider a given dominant $\,b=$
$\sum_{i=1}^n b_i \om_i$ from $(i)$ as a (dominant) weight
for any $\mathfrak{gl}_m$ provided
 $m\ge n-1$. We make the convention that
$\om_{n}=0$ upon the restriction to $\mathfrak{gl}_{n-1}$.

Given sequences $\vec\rr,\vec\ss$ as above, there exists a
{\sf DAHA\~superpolynomial\,}
$\h_{\,\vec\rr,\,\vec\ss}\,(b\,;\,q,t,a)$ belonging to
$\Z[q,t^{\pm 1},a]$ and satisfying the relations
\begin{align}\label{jones-sup}
\h_{\,\vec\rr,\,\vec\ss}\,(b;q,t,a\!=\!-t^{m+1})\!=\!
\tilde{J\!D}_{\,\vec\rr,\,\vec\ss\,}^{GL_{m}}(b;q,t)
\hbox{\, for any } m\!\ge\! n\!-\!1;
\end{align}
its $a$\~constant term %(the coefficient of $a^0$)
is automatically tilde-normalized.
Here one sufficiently large $m$ uniquely
determines $\h_{\, \vec\rr,\,\vec\ss}$. This polynomial
depends only on the isotopy class of
$\t(\vec{\rr},\!\vec{\ss})$ and the color $b$ (the corresponding
Young diagram). \sq
\end{theorem}

It is important (theoretically and practically)
that one can use here  non-symmetric $E_b$ in place of the
symmetric $P_b$; this results in the same super-polynomial.
The $E$\~polynomials are actually beyond Lie theory;
so using $E$\~polynomials is an advantage
for the theory of DAHA superpolynomials, and a clear challenge
for other approaches. Computationally,
the programs written in terms of $E$\~polynomials are much faster
(though not $n!$ times faster).
\vskip 0.2cm

\subsubsection{\sf Khovanov-Rozansky homology}
Concerning the topological
invariance of $\h_{\,\vec\rr,\,\vec\ss}\,(b;q,t,a)$,
it is expected from \cite{ChD1}
that $\h_{\vec{\rr},\vec{\ss}}$
coincides with the
\emph{reduced stable Khovanov-Rozansky} polynomial of the
corresponding knot; see
\cite{KhR1}, \cite{KhR2}, \cite{Kh}. To be precise, it is
conjectured that
\begin{align}\label{khrconj}
\h_{\,\vec\rr,\,\vec\ss}\,(\omega_1\,;\,q,t,a)_{st}=
\tilde{K\!h\!R}_{\hbox{\tiny stab}}(q_{st},t_{st},a_{st}).\
\end{align}
The subscript $st$ indicates that we have made the
substitution to the \emph{standard topological
parameters} given by
\begin{align}\label{qtareli}
&t=q^2_{st},\  q=(q_{st}t_{st})^2,\  a=a_{st}^2 t_{st},\notag\\
&q^2_{st}=t,\  t_{st}=\sqrt{q/t},\  a_{st}^2=a\sqrt{t/q}.
\end{align}
$\tilde{K\!h\!R}_{\hbox{\tiny stab}}$ is the \emph{reduced}
$K\!h\!R_{\hbox{\tiny stab}}$ polynomial
divided by the smallest power of $a_{st}$ and then by
$q_{st}^\bullet t_{st}^\bullet\,$ such that
$\,\tilde{K\!h\!R}_{\hbox{\tiny stab}}(a_{st}\!=\!0)
\in \Z_+[q_{st},t_{st}]$ and the constant term of this
polynomial equals $1$.

The stable KhR\~polynomials are sufficiently understood only
in the uncolored case, so the DAHA construction must be restricted
to $b=\om_1$ here. The minuscule weights may be
managed in the {\em categorification theory\,},
though explicit formulas (to compare with our
ones) are still not available. The reduced setting is also
a problem. Generally, formulas for
KhR\~polynomials are quite a challenge; for uncolored
{\em torus\,} knots, there is recent progress based on Soergel
modules \cite{Mel}.

\setcounter{equation}{0}
\section{\sc Main conjecture}
\subsection{\bf Jacobian factors}
The extension of the theory of compactified
Jacobians and its (entirely) local version, the theory
of Jacobian factors, to modules of arbitrary rank
is not very straightforward.
We mainly follow \cite{PS,GP} for the rank-one definitions,
but there are changes due to $\rk>1$. For instance the 
Jacobian factor now becomes  a (locally) constructible set,
not a projective variety as for $\rk=1$
(in the absence of flags). Its definition as a single space 
becomes significantly less direct than for 
$\rk=1$, but is still not too involved. The points are
$\r$\~modules of {\em deviation\,} $0$,
``non-liftable" to standard ones, which are the
modules of full rank over the 
normalization $\o$ of the ring $\r$ of a given singularity.
Some explicit quasi-projective subvarieties must be
removed from the corresponding {\em compactified
Jacobian factor\,}. 

This compactification is of obvious importance, but not
actually necessary to define our geometric
superpolynomial. We need only its image in the
Grothendieck ring of varieties over $\C$, where it becomes
a union (sum) of the strata formed by standard modules of 
fixed deviations (without any lift).  The
latter are quasi-projective varieties.  Let us mention
without going into detail that our definition of Jacobian
factors in any ranks is connected with 
the theory of affine Springer fibers for {\em non-reduced\,}
spectral curves (their germs, to be precise).

\subsubsection{\sf Rank one modules}
Let $k$ be any field; $\overline{k}$
its algebraic closure, $\o=k[[z]]$,
$\k=k((z))$, $\r$
be any $k$\~subalgebra of $\o$
containing some principle ideal $(z^m)=z^m\o\subset \o$.
We will {\em not\,} assume for now that $\r$ is a ring
of a plane curve singularity, which is by definition an algebra
with two generators over $k$ 
and with $k((z))$ as its field of fractions.  %%%!!!

The corresponding {\em semigroup\,} 
is $\Ga=\Ga_{\r}\equal\{\nu(x),
x\in \r \}$ for 
the {\em valuation\,}
$\nu(x)$ that is the minimal $z$\~degree. The
{\em $\de$\~invariant\,} (the arithmetic genus) of $\r$ is
then $\de=\de_{\r}\equal
\dim(\o/\r)=|\Z_+\setminus \Ga_{\r}|$, where $|\cdot|$
is the cardinality. We will call
$\,\Z_+\setminus \Ga_{\r}\,$ the {\em set of gaps\,} and
denote it $G_{\r}$. Let
$\cc=\cc_\r$ be the {\em conductor\,} of $\r$, the smallest
integer such that $(z^{\cc})=z^{\cc}\o\subset \r$.
One has, $\cc\le 2\de$, where the equality holds exactly
for Gorenstein $\r$ (including plane curve
singularities); this equality is equivalent to
$\de=|\Ga_{\r}\setminus (\cc+\Z_+)|$.
%Does \N include zero?

\subsubsection{\sf Standard modules}
For any nonzero $\r$\~submodule $M\subset \o$,
its {\em $\o$\~degree\,}, the degree with respect
to $\o$ is $\deg(M)\equal$dim$_{k}(\o/ M)$.
One has deg$_{\o}R=\de$. For
arbitrary submodules $M\subset \k$, we set:
\begin{align}\label{deg-def}
\deg(M)\equal\dim(\o/(\o\cap M))-
\dim(M/(\o\cap M)).
\end{align}
This definition is a natural counterpart
of the degree of a divisor
at a given point (here at $z=0$) in the smooth situation.
Actually, we will mostly need the {\em deviation\,}
of $M$ from $\r$, which is
\begin{align}\label{deg-dev}
dev(M)\equal -\hbox{deg}_{\r}(M)=\de-\deg(M).
\end{align}
Let $\De_M=\De(M)\equal\nu(M)$; it is a $\Ga$\~module, i.e.
by definition $\Ga+\De\subset \De$.
The conductor $\cc_M$ of $M$ is naturally the
smallest $c\in \Z_+$ such that $z^{c}\o\subset M$.
Equivalently, it is the smallest $c$ such that
$\De_M\supset c+\Z_+$.

We call
it {\em standard\,} if $\Ga\subset \De_M$, i.e. $M$ contains
some $x\in 1+z\o$. For standard $M$:
$\deg(M)\!=\!|\Z_+\!\setminus\! \De_M|$,
$dev(M)\!=\!|\De_M\!\setminus\! \Ga|$.
We will use the notations $\deg(\De)$ and $dev(\De)$
for the right-hand sides here when $\De=\De_M$ and
 for any standard $\Gamma$-module $\De$ (not
necessarily the set of valuations of a module over $\r$).

{\em Invertible
modules\,} are standard ones with one generator (thus
the generator
must be an element of $1+z\o$).
Accordingly, $1)\ \De_M=\Ga$, $2)\ dev(M)=0$
and $3)\ \cc_M=\cc$ for invertible $M$; vice versa,
any standard $M$
satisfying one of these three conditions is invertible.
The latter implication
is obvious for the first two equalities and the
invertibility follows from $\cc_M=\cc$ due to the
following general lemma.

\begin{lemma}\label{LEM:cond}
For any standard $\r$\~module $M$,
 $\cc_M\le \cc-dev(M)$.
\end{lemma}
{\it Proof.} We follow \cite{PS}, but need a bit different form
of their condition $(ii)$ from Section 1. Also, we provide
the proof since this lemma is an important particular step
of the generalization (below) to arbitrary ranks.
Since $\cc_M-1$ is  missing in $\De_M$ (so missing in $\Ga$
too), all $\,\cc_M\!-\!1-\!\ga\,$ for any $\ga\in \Ga$ are not in
$\De_M$, and one has:
\begin{align*}
&\de-dev(M)\,=\,|\Z_+\setminus \De_M|\,\ge\,
|\{\ga\in \Ga\,\mid\, \nu(\ga)<\cc_M\}|\\
&=|\{\ga\in \Ga\mid \nu(\ga)\!<\!\cc\}|-
|\{\ga\in \Ga\mid \cc_M\le \nu(\ga)\!<\!\cc\}|\,\ge\,
\de-(\cc-\cc_M).
\end{align*}
\vskip -0.5cm
\sq

Any module $M$ can be canonically transformed to a standard
one $M^\circ$ as follows: $M^\circ=z^{-n} M$ for
$n=$min\,$\De_M$. Indeed,
$\deg(z^n M)\!\!=\,\deg(M)+n$ and
$dev(z^n M)=dev(M)-n$ for any $M,n$.
We set $\De^\circ_M=
\De(M^\circ)=\De-n$ and write
$dev^\circ(M)=dev(M^\circ)$,
$\cc_M^\circ=\cc_{M^\circ}$.

Using the map
$\pi: M \mapsto z^{dev(M)}M$,
the set of standard $\r$\~modules $M$ can be identified
with the set of all $\r$\~modules of degree $\de=\de_\r$
(equivalently, of deviation $0$).
 The latter modules have conductors at
least $\cc$ (and therefore at least $2\de$).
Indeed, $\cc_{\pi(M)}=\cc^\circ_M+dev(M)\le \cc$ due
to the lemma. Respectively, we set $\pi(\De)\equal
\De+dev(\De)$ for any standard $\De$.

This observation readily identifies the
set of standard modules considered under $\pi$
with the {\em projective\,}  variety ${}^\pi\!\! J_{\r}$
of {\em all\,} $\r$\~invariant subspaces
of co-dimension $\de$ in $\o/z^{\cc}\o$. It is 
called $\m$ in \cite{PS}.
They actually consider $\o/z^{2\de}\o$, but this results in
the same variety. Recall that $\cc=2\de$ for Gorenstein $\r$.
\smallskip

We will 
identify ${}^\pi\!\! J_{\r}$ when necessary with the set of all
standard modules denoted by $J=J_{\r}$. We also write
$M=M^\circ$ instead of saying that $M$ is standard.
For standard $\De$, let
\begin{align}\label{pi-delta-div}
&J(\De)\equal\{M=M^\circ \,\mid\, \De(M)=\De\};
\hbox{ \, also we represent:} \\
&\pi(\De)=\De+dev(\De)=
\{d_1< d_2<\ldots<d_\de\}\cup \{2\de+\Z_+\}.
\end{align}

Using the theory of Schubert varieties and following
\cite{PS},
\begin{align}\label{piclose}
\hbox{Closure\,}({}^\pi\!\! J(\De))\subset
\medcup_{\,\De'}\,\,{}^\pi\!\! J(\De')
\hbox{\, for \,}
\De' \hbox{\, such that\, } d_i'\ge d_i,
\end{align}
where $1\le i\le \de$.
\smallskip

By the {\em closed fiber\,} of an $\r$\~module $M$, we mean
$\overline{M}=M/\mathfrak{m}M$ for the maximal ideal
$\mathfrak{m}\subset \r$. Its dimension over $k$ will be called
the {\em Nakayama rank\,}. The
{\em $\Ga$\~rank\,}
of $\De\!=\!\De_M$ is the number of generators $b_i$
of $\De$ over $\Ga$. These generators are exactly
all primitive (indecomposable)
elements, i.e. those not in the form $a+\ga$ for
$a\in \De$ and $0<\ga\in\Ga$. Equivalently, they can be
introduced by the following conditions:\
$\De\setminus\{b_i\}$ are $\Ga$\~modules.
The $\Ga$\~rank of $\De$ is no greater than
$\varpi\equal\min\{\ga>0\,\mid\, \ga\in \Ga\}$,
the {\em multiplicity\,}
of the singularity associated with $\r$.

\begin{lemma}\label{NAKIM}
The Nakayama rank {\sl dim}$_k\overline{M}$
of any $\,M$ is no greater than
the $\Ga$\~rank of $\De_M$, and no greater
than $\varpi$. If the primitive elements of $\De=\De_M$
are all smaller than the first non-primitive element in $\De$,
then the Nakayama rank coincides with the $\Ga$\~rank.
\end{lemma}

{\it Proof.} The images of the elements $m_i\in M$ such that
$\nu(m_i)$ are primitive generate
$\overline{M}$. If the primitive elements of $\De_M$
are all smaller than the first non-primitive element in $\De$
then the images of the $m_i$ will form a basis for $\overline{M}$.
\sq

\begin{example}\label{EX}
The simplest
example of a module $M$ with maximal possible Nakayama rank
$\varpi$ is
$M\!=\!\o$. Then $\De_M=\Z_+$, $\{0,1,\ldots,\varpi\!-\!1\}$ are
primitive elements. For any $0\le i< \varpi$, one can diminish
this module $M=\o$ as follows. Let $\o\lan i\ran=
\bigl(\oplus_{j=0}^{i-1} k e_j'\bigr)
\oplus z^{i+1}\o$ for
$e_j'=z^j+\la_j z^i$, where  $\la_j\in k\,
(0\!\le\! j\!\le\! i-1)$ are any (free) parameters.
Then $\De(\o\lan i\ran)=\Z_+\!\setminus\! \{i\}$ and
the primitive elements are $\{0,1,\ldots,\varpi\!-\!1\}
\setminus \{i\};$ so the Nakayama rank of $\o\lan i\ran$
is $\varpi-1$ and it depends on $i\,$ (free) parameters.
\end{example}
\vskip -1.0cm
\sq

\medskip
{\sf Admissibility.}
A $\Ga$\~module $\De$ is called {\em admissible\,} if there
exists $M$ such that $\De_M=\De$. All $\Ga$\~ modules are
admissible for {\em monomial\, } $\r$ (quasi-homogeneous
singularities); see \cite{PS}. Monomial $\r$ are
$k$\~subalgebras in $\o$ generated by pure monomials $z^a$.
This is a significant simplification of the general theory of
$\De$ and $M$; non-admissible $\De$ occur even for the
simplest planar non-monomial $\r$ \cite{Pi}.

\subsubsection{\sf Standard flags} These are the key objects
in \cite{ChP}. A flag
$\vect{M}=\{M_0\subset M_1\subset\ldots\subset M_\ell\}$ of
$\r$\~modules $M_i\subset \o$ as above is called
{\em full $\nu$\~increasing\,} if
the corresponding {\em $\Delta$\~flag}
$\vect{\De}=\De(\vect{M})=
\{\De_0\subset \De_1\subset\ldots\subset \De_\ell\}$
for $\De_i=\De(M_i)$ is full increasing, i.e.
satisfies the following:
\begin{align} \label{deflags}
\De_{i}=\De_{i-1}\cup \{g_i\},\
 g_1<g_2<\ldots<g_\ell, \ \,1\le i\le \ell.
\end{align}

It is called {\em standard\,}
if it is full $\nu$\~increasing
and $M_0$ is standard (i.e. $\De_0\supset \Ga$),
which automatically implies  that all $M_i$ are standard.
We will call such (standard) flags simply $\ell$\~flags and
define the length of $\vect{M}$ by $l(\vect{M})=\ell$.

The set $J^{\ell}$ of all standard $\ell$\~flags $\vect{M}$
is called the
{\em flagged Jacobian factor\,} in \cite{ChP}. It can
be supplied with a natural structure of a
{\em quasi-projective\,} variety due to
following proposition, which is actually Proposition 2.3
from \cite{ChP} (somewhat related to \cite{ORS}, Section 2.1).
Its proof there did not use that $\r$ is of planar (with $2$
generators), so we will omit it (and the corresponding
equations of $J^{\ell}$).

For an arbitrary module $M$, we set  $M^{\{i\}}=M\cap
(z^i)$, which is obviously an $\r$\~module,
$\overline{M}^{\{i\}}$
the image of $M^{\{i\}}$ in $\overline{M}$. One has:
dim\,$M^{\{g\}}/M^{\{g\!+\!1\}}\!=\!1$ for
$g\!\in\! \De_M$\, and $0$ otherwise.
These dimensions can drop to $0$ upon taking the
closed fiber, which is the
reduction modulo $\mathfrak{m}$.

%dim\,$\{M}^{\{g\}}/\bar{M}^{\{g\!+\!1\}}\le 1$.
%\newcommand*{\vect}[1]{\overrightarrow{\mkern0mu#1}}

%What does the ' mean in \overline{M}_i' ?
\begin{proposition}\label{NESTED}
(i) For a full increasing $\De$\~flag
$\vect{\De}=\bigl\{\De_0\!\subset
\De_1\!=\!\De_0\cup\{g_1\}\subset\!
\cdots\!\subset \De_\ell\!=
\!\De_0\cup\{g_1,\ldots,g_\ell\}\bigr\}$,
all $\De_0\cup \{g_i\}$ for $1\!\le\! i\!\le\! \ell$
are $\Ga$\~modules. This implies that for any standard flag
$\vect{M}$ of length $\ell$, one has:
$\mathfrak{m}M_{\ell}\subset M_i$ for $0\le i\le \ell$.
Accordingly,
these $\r$\~modules $M_i$ are uniquely determined by their images
$\overline{M}'_i\equal M_i/\mathfrak{m}M_{\ell}$ in
$\overline{M}_{\ell}$, which follows
from Nakayama Lemma.

% k-subspaces of what? I think it should be \overline{M}_{\ell}
(ii) An arbitrary standard module $M$
of Nakayama rank\, {\sl dim}$_k\overline{M} \ge \ell$
appears as the top module $M_\ell=M$ in at least one
$\ell$\~flag for some $\vect{\De}$.
Such $M_\ell$ form a
projective subvariety in $J_{\r}$ due to the
upper semicontinuity of Nakayama rank. Given $\vect{\De}$
and $M_{\ell}$,
all possible $M_0$  are in one-to-one correspondence with
$k$\~subspaces
$\overline{M}'_{0}\not\subset \overline{M}_{\ell}^{\{1\}} $
of codimension $\ell$ such that
\begin{align}\label{mocond}
\hbox{dim\,}_k\,(\overline{M}'_0\!+
\!\overline{M}_{\ell}^{\{g_i\}})/
(\overline{M}'_0\!+\!\overline{M}_{\ell}^{\{g_i\!+\!1\}})=1
\for 1\le i\le \ell.
\end{align}

(iii) Given $\vect{\De}$ and a standard $M_{\ell}$, the
corresponding
$\ell$\~flags $\vect{M}$ exist if and only if
\begin{align}\label{mocondd}
\hbox{dim\,}_k\,\!\overline{M}_{\ell}^{\{g_i\}}/
\overline{M}_{\ell}^{\{g_i\!+\!1\}}=1
\for 1\le i\le \ell.
\end{align}
These flags can be identified via Nakayama Lemma with full
$\ell$\~flags of $k$\~subspaces
$\overline{M}'_{0}\!\subset\ldots\!\subset
\overline{M}'_{i}\!\subset\ldots
\!\subset\overline{M}_{\ell}$
such that
\begin{align}\label{moconddd}
\overline{M}'_i\!+\!\overline{M}_{\ell}^{\{g_i\!+\!1\}}=
\!\overline{M}_{\ell} \for 1\le i\le \ell-1.
\end{align}
Provided (\ref{mocondd}) and given $M_0$,
this space is biregular to $\mathbb{A}^{\ell(\ell-1)/2}$.
\sq
\end{proposition}

Example \ref{EX} can be interpreted
as follows. All $1$-flags (for $\ell=1$) with $M_1=\o$
are $\{M_0=\o\lan i\ran\subset \o=M_1\}$. Given $i$, the
corresponding
subspaces of $\o/z^\varpi\o\cong k^\varpi$ of codimension $1$
must contain $(z^{i+1})$ but not $(z^{i})$. Thus the space
of such flags is biregular to
$\mathbb A^{i}$, which matches Example \ref{EX}.
Here $i\ge 1$; note that the removal of $i=0$ from $\Z_+$ makes
the corresponding module
non-standard.

Proposition \ref{NESTED} identifies $J^\ell$ with the space of
pairs
$$\bigl\{M_\ell\in J\lan\ell\ran,
\{\overline{M}'_i,0\le i\le \ell\}
\bigr\},$$
where (a)\, $J\lan\ell\ran\subset J$ is a projective subvariety
${}^\pi\!\! J$ of
modules of Nakayama rank no smaller than
$\ell$ considered under the map $\pi$ above,
(b)\, $\overline{M}'_{0}\subset
\overline{M}_{\ell}$ is a $k$\~subspace of codimension $\ell$
such that $\overline{M}_{0} \not\subset
\overline{M}_{\ell}^{\{1\}}$, and
(c)\, $\{\overline{M}'\}$ is {\em any\,} full flag in
$\overline{M}_{\ell}/\overline{M}'_{0}$.
Following the
notation for $\rk\!=\!1$ above, we  will denote this variety
by ${}^\pi\!\! J^{\ell}$ or simply by $J^{\ell}$.

Note that allowing
$\overline{M}_{0} \subset \overline{M}_{\ell}^{\{1\}}$
means that we consider {\em all\,}
flags, not only standard; $M_\ell$ will be always assumed
standard when we apply the map $\pi$ above. This gives the
following {\em projectivization\,} $\hat{J}^\ell$ of
$J^\ell$ (of ${}^\pi\!\! J^\ell$ to be exact).
 It is a bundle
over $J\lan\ell\ran$ with the fibers that are Grassmannians
$Gr(n,n-\ell)$ for $n=\hbox{dim}_k \overline{M}_\ell$
equipped with (any) full $\ell$\~flags in the corresponding
quotient spaces.

%The corresponding boundary is for $\overline{M}_{0}\subset
%\overline{M}_{\ell}^{\{1\}}$.
The following general identification is helpful to
describe the boundary $\hat{J}^\ell\setminus\! J^\ell$:
\begin{align}\label{mbound}
&\{\,\vect{M}\ \mid\ l(\vect{M})=\ell,\ M_{\ell}=M,\
M_{\ell-1}\subset zM\,\}\\
&\cong
\{\,\vect{M}'\ \mid\ l(\vect{M}')=\ell-1,\ M_{\ell-1}'=
z^{-\omega}M^{\{1\}}\,\}.\notag
\end{align}
Here $\omega=\om_M=$\,min\,$\{\De_M\setminus\! \{0\}\}$, so
$z^{-\omega}M^{\{1\}}$ is a standard module
corresponding to $(\De_M\setminus\{0\})-\om$
for the pointwise subtraction of $\omega$.
The Nakayama rank of this module is no smaller than that of
$M$ minus $1$ by construction.
Thus the boundary can be generally interpreted via the
map $M\mapsto z^{-\varpi}M^{\{1\}}$
and its counterparts for higher $M^{\{i\}}$.

We will also note that given $\vect{\De}$,
conditions (\ref{mocondd})
are necessary and sufficient for a standard
$M_\ell$ to have an extension to a standard flag $\vect{M}$
with such $\vect{\De}$. This of course results in the
admissibility of $\De_0$ and all $\De_i$ for  $1\le i\le \ell$
and also in admissibility of $\De_0\cup\{g_i\}$.

Using the standard facts from the theory of Schubert
varieties (see \cite{PS}), we arrive at the following
corollary.

\begin{corollary}\label{BOUNDARY}
Given $\vect{\De}$, the closure  in $J^\ell$
of the
corresponding subvariety
$$
J^\ell(\vect{\De})\equal \{
\vect{M}\,\mid\, \De(\vect{M})=\vect{\De}\}
$$
belongs to the union
of $J^\ell(\vect{\De}^\flat)$ for $\{g_i^\flat\}$,
which are not now assumed increasing in the sense of
(\ref{deflags}), such that :

(a)\, $\pi_\ell(\De_\ell)$
and $\pi_\ell(\De^\flat_\ell)$ satisfy conditions
(\ref{piclose}) for $\pi_\ell(\De)\!=\!\De\!+\!dev(\De_\ell)$,
(b)\, $\pi_\ell(\De_0)$ and $\pi_\ell(\De^\flat_0)$
satisfy the same condition, and (c)\, $g_i^\flat\ge g_i$
for all $1\le i\le \ell$.

If $M_\ell$ is fixed, then the
union above is over all flags $\vect{\De}^\flat$ with
standard $\De^\flat_0$ and $\{g_i^\flat\}$ satisfying
conditions $(b,c)$; it coincides with
the corresponding closure. If $M_\ell$ and $M_0$ are fixed,
then this closure is the union over all permutations
$\{g_i^\flat\}$ of the initial (increasing) set $\{g_i\}$, where
relations (\ref{moconddd}) are omitted. Also, the embedding
of the closures of the subvarieties corresponding to such
$\{g_i^\flat\}$ are with respect to the Bruhat ordering of
permutations
(generated by consecutive transpositions of the pairs $g_i<g_j$).
%in Proposition \ref{NESTED},$(iii)$.
\sq
\end{corollary}

The simplest example is as follows. Let
$\{g^\flat_i=g_{\ell-i}\}$, i.e. it is
the {\em decreasing\,} ordering of the initial (increasing)
sequence. Then for fixed $M_\ell,M_0$,
a unique flag of vector spaces  exists
for such $\{g^\flat_i\}$ in $M_\ell/M_0$
from Proposition \ref{NESTED}, $(iii)$, which
is formed by the images of $M_\ell^{\{g_i\}}$.
So this is one point in this case and it
belongs to all closures corresponding to
any intermediate $\{g_i^\flat\}$. When $\ell=2$ (for $2$\~flags),
this point is the whole boundary.

\comment{
It was conjectured that the
admissibility of the latter and $\De_0$ is sufficient for the
admissibility of $\vect{\De}$. Let us prove it.
\begin{proposition}\label{ADMFLAG}
Let as assume that standard $\Ga$\~modules $\De_0$ and
$\De=\De_0\cup\{g\}$ are admissible for some
$g\in \Z_+\setminus\! \De_0$.
Then
the flag $\{\De_0\subset\De\}$ is admissible. Generally,
standard (increasing) $\De$\~flags are admissible if and only
if their constituents $\De_i$ are admissible.
\end{proposition}
{\it Proof.} Let us assume that   $\{\De_0\subset\De_1\}$
is the pair with the greatest $dev_0=|\De_0\setminus\Ga|$
such that $\De_i$ are admissible but this flag is not.
Let  $g'=\hbox{nax\,}(\Z_+\setminus \De_0)$. Then
for any $M'$ with $\De(M')=\De'$, its
the intersection with
 can be lifted to some
$M'$ for $\De'$.
}

\subsection{\bf Higher ranks}
Let us consider now $\r$\~submodules $\m\subset \o^{\rk}$
for $\rk\ge 1$ such that $\k\cdot \m=\k^{\rk}$.
%assuming that they generate the
%whole $\k^{\rk}$ over the field $\k$ of fractions of $\o$.
From now on, we fix an  $\o$\~basis
$\{\ep_i,\,0\le i\le \rk-1\}\subset \o^{\rk}$
and the standard form $(\ep_i,\ep_j)=\de_{ij}.$

\subsubsection{\sf Basic definitions}
As in the case of $\rk=1$, the $\o$\~degree of
$\m\subset\o^{\rk}$ is
naturally
 $\deg(\m)\equal$dim$_k(\o^{\rk}/ \m)$; for
arbitrary submodules $\m\subset \k^{\rk}$ of rank $\rk$
over $\k$, we set:
\begin{align}\label{deg-def-rk}
\deg(\m)\equal\hbox{dim}_k\,(\o^{\rk}/(\o^{\rk}\cap \m))-
\hbox{dim}_k\,(\m/(\o^{\rk}\cap \m)).
\end{align}
The {\em deviation\,}
of $\m$ from $\r^{\rk}$ is
$dev(\m)\equal \rk\cdot\de\!-\!\hbox{deg}_\o(\m)$.

The {\em projections\,} of $\m$ onto $\ep_i$,
their $\De$\~modules and the corresponding
deviations are:
\begin{align}\label{de-proj}
&\m^{\lan i\ran}\equal(\m,\ep_i),\,
\De^{\lan i\ran}\!=\!\De^{\lan i\ran}(\m)\!=\!
\De(\m^{\lan i\ran}),\notag\\\
&dev^{\lan i\ran}(\m)=
dev(\De^{\lan i\ran})
\for 0\le i\le \rk-1.
\end{align}
These definitions significantly depend on the choice of the
basis $\{\ep_i\}$.
\smallskip

{\sf $\De$\~Components.}
There is another way to associate a sequence of $\Ga$\~modules
to $\m$; it will be the
main tool for us. It requires
a choice of {\em cumulative valuation\,}
$\upsilon:\k^{\rk}\to \R$, extending the valuation $\nu$ above:
\begin{align}\label{de-cumv}
\upsilon(\sum_i m_i\ep_i)\!=\!
\min_i \{\nu(m_i)\!+\!\upsilon(\ep_i)\}
\for m_i\!\in\! \k, \ 0\!\le\! i\!\le\! \rk\!-\!1,
\end{align}
where $\upsilon_i\equal\upsilon(\ep_i)$ can generally be
arbitrary real numbers. The conditions
$0\le \upsilon_i\neq \upsilon_j<1$
for $i\neq j$ are sufficient below and will be
imposed unless stated otherwise. Generally,
one needs the following weaker relations:\,
$\upsilon_i-\upsilon_j\not\in \Z$ for
$i\neq j$.  The default will be:
$\upsilon_i\equal i/\rk$; though
{\em they will be multiplied by $\rk$ in the last two
sections of the paper}.

For such $\upsilon$ and $0\le i\le \rk-1$,
the {\em $\De$\~components\,} are
$\Ga$\~modules
defined as follows. Setting here and below
$\De(\m)\equal\upsilon(\m)$,
\begin{align}\label{de-comp}
\De^{(i)}(\m)\!=\!\{a\in \Z_+\mid a\!+\!\upsilon_i\in
\De(\m)\},\
dev^{(i)}\!=\!dev(\De^{(i)}(\m)).
\end{align}
We will also use the notation $\De^{(i)}_{\m}=\De^{(i)}(\m)$
or simply $\De^{(i)}$ when there is no fear of confusion.
It is immediate to check that $\De^{(i)}\subset
\De^{\lan i\ran}$.
If $\upsilon_i-\upsilon_j\in \Z$ for some $i\neq j$, then
we face {\em wall-crossing\,}, which will not be studied in
this work. Generally speaking, the $\Z$\~commensurable
components of $\upsilon$ must be considered, which
are $\Ga$\~modules. We will need their number to be $\rk$
(maximal possible) in this paper.

The stratification of $\j_{\r,\rk}$ in
Section   \ref{SEC:JACF}   will be based on 
%Our definition of {\,\em Piontkowski-type cells\,}
$\{\De^{(i)}\}$; the modules $\{\De^{\lan i\ran}\}$
will be mostly needed in this paper to define the map $s\!t$
below, in (\ref{pre-st}).
Importantly,
$\sum_{i=0}^{\rk-1} dev^{(i)}= dev(\m)$, as well as the
analogous formula $\sum_{i=0}^{\rk-1} \hbox{deg}^{(i)}
= \hbox{deg}_\o(\m)$, where
deg$^{(i)}\!=$\,\,deg\,$\De^{(i)}\!=|\Z_+\setminus\!\De^{(i)}|$
and deg$_\o(\m)$.
These relations generally do not hold for $\De^{\lan i\ran}$.

It is interesting to examine the dependence of
$\{\De^{(i)}\}$ upon the change of basis
$\{\ep_i\}\mapsto \{\hat{\ep}_i\}$ without
varying the numbers $\upsilon_i$. This is actually close
to {\em tropical geometry\,}. Let us provide
the following lemma.

\begin{lemma}\label{SORTDEL}
We fix $\upsilon$ such that $0\le \upsilon_0<\upsilon_1<\ldots
<\upsilon_{\rk-1}<1$ and an algebraic closure $\overline{k}$ of
$k$.
Let $\m\subset\o^{\rk}$ be an $\r$-module such that
$\k\cdot\m=\k^{\rk}$. Define an upper-triangular linear
substitution on $\m$ as follows:
$\ep_i\mapsto\hat{\ep}_i\!=
\!\ep_i\!+\!\sum_{j<i}c_i^j\ep_j$,
where $c_i^j\in \overline{k}$ (apart from some finite set of
affine hyperplanes). If
 $\hat{\m}\equal\{\sum_{i=0}^{\rk\!-\!1}
(y,\ep_i)\hat{\ep}_i\, \mid\, y\in \m\}$,
then $\nu(\hat{\m})$ has the following components:
\begin{align*}
&\hat{\De}^{(0)}=\De^{(0)}\cup\ldots\cup\De^{(\rk\!-\!1)},\ \,
\hat{\De}^{(1)}=\medcup\,\,_{i\neq j}\,
\{\De^{(i)}\cap\De^{(j)}\},\,\ldots\,,\\
&\hat{\De}^{(\rk\!-\!2)}=\medcup\,\,_{i}
\bigl\{\cap_{j\neq i}\De^{(j)}\bigr\},\ \,
\hat{\De}^{(\rk\!-\!1)}=\De^{(0)}\cap
\ldots \cap \De^{(\rk\!-\!1)}.
\end{align*}
If $k$ is infinite, $\overline{k}$
is not needed here: $c_i^j\in k$. Also,
lower-triangular
substitutions $\ep_i\mapsto\hat{\ep}_i\!=
\!\ep_i\!+\!\sum_{j>i}c_i^j\ep_j$,
do not influence $\De^{(i)}$.
\end{lemma}
\smallskip
{\it Proof.}
The case of upper-triangular transformations involving two
neighboring $\ep_i,\ep_{i+1}$ is sufficient to consider for
both claims.
Suppose $y=t^m\ep_i+c t^m\ep_{i+1}+ t^m\sum_{j>i+1}a_j\ep_j
\mod (t^{m+1})\o^2$, and
$y'=t^n\ep_{i+1}+t^n\sum_{j>i+1}b_j\ep_j \mod (t^{n+1})\o^2$
where $c,a_j,b_j$ are constants.
The corresponding valuations are then
$m+\upsilon_i$ and $n+\upsilon_{i+1}$. It suffices to consider
elements only in the form of $y$ and $y'$.

Obviously
lower-triangular substitutions will not change the valuations
$y$ and $y'$. Let us justify the first claim.
Given $y$, let's first assume that
elements of the form of $y'$
appear only for $n\not=m$.
Then a single change
$\ep_{i+1}\mapsto\ep_{i+1}+\al\ep_{i}$ for $\al\in k^*$
will preserve the valuation of $y$ if
$\al+c\neq 0$ and will change the valuation of $y'$
 to $n+\upsilon_i$.
Thus  $n$ will be added to $\De^{(i)}$
and will disappear from $\De^{(i+1)}$ in this case.

Note that if $\al+c=0$ here,  then $m\in \De^{(i)}$ and
$n\in\De^{(i+1)}$ will be transposed, so a different
transformation will occur. Therefore for $k=\mathbb F_2$,
moving $n$  from $\De^{(i)}$
to $\De^{(i+1)}$ is not possible over $k$ if $c=1$ and
more general substitutions (not only upper-triangular)
must be considered or, as in the lemma,
one may use $\overline{k}$.

Now given $y$, suppose $y'$ exists with $n=m$,
which means that $m\in \De^{(i)}\cap \De^{(i+1)}$.
Then $m$ remains in both, $\De^{(i)}$ and $\De^{(i+1)}$,
upon the substitution above for any $\al$.

We see that the resulting components for $\al+c\neq 0$ become
$\De^{(i)}\cup\De^{(i+1)},\De^{(i)}\cap\De^{(i+1)}$.
Then we will iterate by considering another pair of
neighboring indices and so on.
A sufficiently large extension of $k$ is
generally necessary here for finite $k$.
\sq

Note that the sum of cardinalities of $\hat{\De}^{(i)}\setminus
(N+\Z_+)$ for sufficiently large $N$
coincides with that for $\De^{(i)}\setminus
(N+\Z_+)$ calculated for  $\De^{(i)}$. It is clear
combinatorially and because this sum is
the dimension of $\m/(z^N)\o^{\rk}$.

\medskip
{\sf Conductors.} Accordingly, we have two
different sequences
of conductors $\{\cc\!\lan i\ran\}$ and $\{\cc(i)\}$,
those for $\m^{\lan i\ran}$ and, correspondingly,
for $\m^{(i)}$.

Then $\cc(i)$ is no smaller than
$\cc\!\lan i\ran$:
\begin{align}\label{CC-ups}
\hbox{\sf C}(\{\De^{(i)}\})\equal
\sum_{i=0}^{\rk\!-\!1}(z^{\cc(i)}\o)\ep_i\subset
\hbox{\sf C}(\{\De^{\lan i\ran}\})\equal
\sum_{i=0}^{\rk\!-\!1}(z^{\cc\lan i\ran}\o)\ep_i.
\end{align}

The simplest examples of modules
$\m$ are direct sums. Let $\De^{(i)}=\De(M_i)$
for some $M_i\subset \o$ (i.e. $\De^{(i)}$ are admissible).
Then
$\De^{(i)}(\m)\!=\!\De^{\lan i\ran}(\m)$ for
$\m\!=\!\oplus_{i=0}^{\rk\!-\!1} M_i\equal
\sum_{i=0}^{\rk\!-\!1} M_i \ep_i$, and
$\sum_{i=0}^{\rk\!-\!1} (z^{\cc(i)}\o)\ep_i\!=\!
\cc(\m).$ The coincidence
$\De^{(i)}\!=\!\De^{\lan i\ran}$ for all $i$ is a special feature
of direct sums.
\vskip 0.2cm

Let $\cc(\upsilon)$ be the smallest
number $\cc$ such that $(\cc+\R_+) \cap \upsilon( \o^{\rk})$
belongs to $\upsilon(\m)$ and
$\hbox{\sf C}(\upsilon)\equal
\{v\!\in\! \o^{\rk}\mid \upsilon(v)\!\ge\! \cc(\upsilon)\}$.
The {\em usual conductor\,} $\hbox{\sf C}(\m)$
of a module $\m$ extends {\sf C}$(\upsilon)$:
\begin{align}\label{condum}
&\hbox{\sf C}(\upsilon)\subset
\hbox{\sf C}(\m)\equal
\{y\in \m\, \mid\, y\o \subset \m\}.
%\subset\hbox{\sf C}(\{\De^{(i)}\}).
\end{align}
Indeed, $\hbox{\sf C}(\upsilon)$ is an $\o$\~module
inside $\m$ and therefore
inside $\hbox{\sf C}(\m)$.

\begin{definition}\label{DEF:PURE}
We call an $\r$\~submodule $\m$ {\sf pure\,} with respect to
$\{\ep_i\}$ if
$\hbox{\sf C}(\m)\!=\!\oplus_{i=0}^{\rk-1}(z^{n_i}\o)\ep_i$
for $n_0\le n_1\le \ldots \le n_{\rk-1}$. Any submodule $\m$
becomes
pure upon a proper substitution $\{\ep_i\mapsto \hat{\ep}_i\}$
in $y\in\m$ for some basis $\{\hat{\ep}_i\}$ of $\o^{\rk}$. The
numbers $n_i$  are uniquely defined by $\m$. \sq
\end{definition}

Upper-triangular substitutions
from Lemma \ref{SORTDEL} do not alter $\hbox{\sf C}(\m)$
and transform pure modules to pure ones (for a fixed basis
$\{\ep_i\}$); the numbers $n_0,\ldots,n_{\rk-1}$
remain unchanged.
Recall that upon sufficiently general substitutions, we have:
\begin{align} \label{De-incr}
\hat{\De}^{(0)}\supset\hat{\De}^{(1)}\supset
\cdots\supset \hat{\De}^{(\rk-1)}.
\end{align}
Note that $\hat{\De}^{(i)}$ are standard for standard $\De^{(i)}$.

These embeddings readily give that
$\cc(i)\le \cc(j)$ for $i<j$, $\cc(i)\le n_i$;
also,
$\cc(\rk-1)=n_{\rk-1}$. Combining (\ref{condum})
with Lemma \ref{LEM:cond}, we
conclude that for standard $\{\De^{(i)},i=0,\ldots,\rk-1\}$
one has:
\begin{align}\label{nilecc}
n_i\le n_{\rk\!-\!1}=\cc(\rk\!-\!1)=\cc(\upsilon)\le
\cc-|\hat{\De}^{(\rk\!-\!1)}\setminus \Ga|\le \cc=\cc_\r.
\end{align}

Note that
$\upsilon(\hbox{\sf C}(\m))$ is a $\Z_+$\~module inside
$\upsilon(\m)=\cup_{i=0}^{\rk\!-\!1} \De^{(i)}.$ Thus
it belongs to
$\cup_{i=0}^{\rk-1} (\cc(i)+\Z_+)$, but they do not
generally coincide.
\smallskip

\subsubsection{\sf Standard modules}
We call $\m$ {\em standard\,} (or $\rk$\~standard) if
$\m$ generates $\o^{\rk}$ over $\o$, equivalently,
the image of $\m$ in $\o^{\rk}/\,z\o^{\rk}$
generates $\o^{\rk}/\,z\o^{\rk}$.
In this case we can find
$y_i\in \m\cap (\ep_i+z\o^{\rk})$ for
$0\le i\le \rk-1$.
This gives that $\m$ is standard if and only if
all $\De^{(i)}(\m)$ are standard (i.e. contain $0$)
for any choice of the valuation
$\upsilon$. It is always assumed to satisfy the
conditions $\upsilon_i\!-\!\upsilon_j\not\in\Z$ for $i\neq j$,
and we impose relations $0\!\le\!\upsilon_0\!<\!\ldots
\!<\!\upsilon_{\rk-1}\!<\!1$ unless stated otherwise.

Recall that
$\deg(M)\!=\!\deg(\De_M)\!=\!|\Z_+\!\setminus\! \De_M|$,
$dev(M)\!=\!dev(\De_M)\!=\!|\De_M\!\setminus\! \Ga|$ for standard
modules $M\subset \o$. We set for standard $\m$:
$$
dev(\m)\equal \rk\cdot\de-\hbox{dim}\o^{\rk}/\m=
\sum_{i=0}^{\rk-1} dev(\De_{\m}^{(i)}).
$$

A generalization of the inequality from
Lemma \ref{LEM:cond}, the key for making $J_\r$ a projective
variety, is as follows:

\begin{lemma}\label{LEM:ineq}
For a standard $\r$-module $\m$, $\rk\cdot\cc\ge dev(\m)
+\sum_{i=0}^
{\rk-1}n_i$,
where
$\hbox{\sf C}(\m)\!=\!\oplus_{i=0}^{\rk-1}(z^{n_i}\o)\ep'_i$
for  $n_0\le n_1\le\cdots\le n_{\rk-1}$ and some $\o$\~basis
$\{\ep_i'\}$.
 Also, $n_i\le \cc=\cc_\r$, which is
due to (\ref{nilecc}).
\end{lemma}
{\it Proof.} We can assume that $\m$ is pure
and relations (\ref{De-incr}) hold. Observe that
$z^{-1}\hbox{\sf C}(\m)\cap \m={\sf C}(\m)$ because
 if $y\in z^{-1}\hbox{\sf C}(\m)\cap \m$, then
$y(1+z\o)\subset \m$ and this means $y\in{\sf C}(\m)$.
The reverse inclusion is clear.

For every $j\in \Ga$, we choose
an arbitrary $x_j\in\r$ such that $\nu(x_j)=j$. Then
$
\bigl(\oplus_{i,j}\, k (z^{n_i-1}/x_j)\,\ep_i\bigr)
\cap \m=\{0\}, \hbox{\,\, where\, } 0\le j\le n_i,\ j\in \Ga.
$

Indeed, for any $y$ in this intersection, let $x_j$ be
the element with maximal $j$ in its decomposition;
it can appear with several $\ep_i$.
Then $x_j y \in z^{-1}\hbox{\sf C}(\m)\cap \m={\sf C}(\m)$ and
dim\,$\o^{\rk}/\m\ge \sum_{i=0}^{\rk-1}
|\Ga\!\ni\! \ga\! <\! n_i|$.
However, $|\{\Ga\!\ni\! \ga\! <\!n\}|\ge \de\! -\!\cc+n$ as
in the proof of Lemma \ref{LEM:cond}, and
$$ \rk\cdot\de-dev(\m)=\hbox{dim\,}\o^{\rk}/\m\ge
\sum_{i=0}^{\rk-1}(\de-\cc+n_i),$$
which gives the required inequality.\sq
\vskip 0.2cm

Generalizing {\em invertible
modules\,} in $\rk=1$, we now
introduce {\em minimal modules $\m$\,},
which by definition
are standard and have exactly $\rk$ generators over $\r$.
Equivalently, $\De^{(i)}(M)=\Ga$ for all $0\le i\le \rk-1$.
These modules form a ``big cell" in the variety of all
standard ones, to be introduced below.

There is a natural
action of $A\in GL(\rk,\k)$ on any modules $\m$ defined
via the corresponding substitutions:
\begin{align}\label{GLaction}
A: \m\ni y=\sum_{i}c_i\ep_i\mapsto
\sum_{i}c_i a_i,\ A_i=(a_0,\ldots,a_{\rk-1}),
\end{align}
i.e. $\ep_i$ are replaced by the corresponding columns of $A$.
If $A\in  GL(\rk,\o)$, then $A(\m)\subset \o^n$ (we always
assume in this paper that $\m\subset \o^{\rk}$),
standard modules go to standard and minimal ones go to
minimal.

For any submodule $\m'\in \o^{\rk}$, we set
\begin{align}\label{pre-st}
&d'\lan i\ran\equal\hbox{min\,}\De^{\lan i\ran}\!=\!
\hbox{min\,}\Delta((\m')^{\lan i\ran}),\\
&s\!t:\,\, \m'\ni y'\,=\,\sum_{i}c_i\ep_i\,\,
\mapsto\,\, y=s\!t(y')\,=\,
\sum_{i}z^{-d'\lan i\ran}c_i\ep_i\in \m,\notag\\
&\pi_{\{d'\lan i\ran\}}: \m\ni y\!=\!\sum_{i}c_i\ep_i\mapsto
\sum_{i}z^{d'\lan i\ran}c_i\ep_i\in \m',\,
s\!t\hbox{\tiny$\circ$}\pi_{\{d'\lan i\ran\}}=id.\notag
\end{align}
Here the image of $s\!t$ is an $\r$\~submodule
$\m=\hbox{diag}(z^{-d'\lan i\ran})(\m')$ in $\o^{\rk}$; it is
generally not standard.
%However for any
%$\m$, there exists $A\in GL(\rk,\o)$ such that
%$A(\m)$ is pre-standard.

\begin{definition}\label{DEF:LIFT}
(A) A submodule $\tilde{\m}\subset\o^{\rk}$ is called
\,{\sf liftable\,}
if there exists $A\in GL(\rk,\o)$ such that
$\m'=A(\tilde{\m})$ and conditions  $(a,b)$ hold:

(a) $\m'$ is {\sf pure\,} with respect to the basis $\{\ep_i\}$,
which means that its conductor is in the form
{\sf C}$(\m')=\oplus_{i=0}^{\rk\!-\!1} (z^{n_i'})\ep_i$
and also the inequalities
$0\!\le\!n_0'\!\le\ldots\le\! n_{\rk-1}'$ are imposed;

(b) $\m=s\!t(\m')$ is {\em standard\,}
and pure for $\{\ep_i\}$, which includes the inequalities
$n_0\!\le\ldots\le\! n_{\rk-1}$, and also
$d'\lan 0\ran\ge d'\lan 1\ran\ge\ldots\ge
d'\lan \rk\!-\!1\ran\!\ge\!0$, where
$\{d'\lan i\ran\}$ for $\m'$ are from $(\ref{pre-st})$.
\vskip 0.2cm

(B)
One has that $d'\lan i\ran\!=\!n'_i\!-\!n_i$ due to
$\m=s\!t(\m')$ in $(b)$ for all\, $i$.
If there exists $\si\in \N$\,
such that $d'\lan i\ran\!=\!\si\!-\!n_i$ for
$i\!<\! i_{\ast}$, \
$0\!\le\! d'\lan i_{\ast}\ran+n_{i_{\ast}}\!<\!\si$,\ and
$d'\lan i\ran\!=\!0$ for $i\!>\!i_{\ast}$ for some
$0\!\le\! i_{\ast}\!\le\! \rk\!-\!1$, then $\m'$
and the
corresponding $\tilde{\m}$ from its $GL(\rk,\o)$\~orbit
are called {\sf $\si$\~minimal
liftable}.
\vskip 0.2cm

(C) Vice versa, given a standard pure $\m$,
some $\partial\ge dev(\m)$
and $\si\ge \cc$, there exists a unique  $\si$\~minimal
sequence $\{d_0 \ldots d_{\rk\!-\!1}\}$ such
that $n_i'=n_i+d_i$ are non-decreasing (as for  $\{n_i\}$),
and $d_0+ \ldots+ d_{\rk\!-\!1}=\partial$.
Then pure $\m'$ of deviation $\,dev(\m)-\partial\,$
is defined as the image of the map
$\pi_d: \m\ni y=\sum c_i \ep_i\mapsto
y'=\sum c_i z^{d_i}\ep_i$. The module $\m$ is then
$s\!t(\m')$ from $(b)$.
\sq
\end{definition}

Actually, $\si=\cc$ is sufficient in this paper.
An important particular case is
when $n_0'=\cdots=n_{\rk-1}'=\cc$ in $(a)$, i.e.
$\tilde{\m}$ and $\m'$ are with the conductor
$(z^\cc)\o^{\rk}$. Then $\cc$\~minimal
liftable modules for these $\{n_i'\}$
are simply all liftable ones. Indeed,
$\m=s\!t(\m')$ is standard and
$d'\lan i\ran\!=\cc\!-\!n_i$. We need only to check
that $d'\lan \rk\!-\!1\ran\!\ge\!0$, which holds automatically
due to (\ref{nilecc}). The inverse of $s\!t$ is
$\pi: \m\ni y\!=\!\sum c_i \ep_i\mapsto y'\!=
\!\sum c_i z^{\cc-n_i}\ep_i$
in this case.

Generally $s\!t^{-1}$ is
$\pi_d: \m\ni y=\sum c_i \ep_i\mapsto
y'=\sum c_i z^{d_i}\ep_i$ for $d_i=d'\lan i\ran$;
the inequalities
for $\{d_i\}$
mean that it {\em flattens\,} the conductor of $\m$
toward making it totally flat, i.e. toward
$z^{\hbox{\tiny const}}\o^{\rk}$ for some {\small const}.

Concerning the $\si$\~minimality,
the inequalities for $\{d_i\}$ alone may be insufficient
to ensure that there is only one such set
for a given $\m$ subject to
 $dev(\tilde{\m})=dev(\m)-\partial=dev(\m')$ from $(C)$.
The minimality condition there means that we are supposed to
do full {\em flattening\,} up to $\si$ of the
conductor consecutively for
$i=0,1,\ldots$\ .
\medskip

\subsubsection{\sf Jacobian factors and cells}\label{SEC:JACF}
Let us fix a valuation $\upsilon$. Recall that
$\De(\m)\equal\upsilon(\m)=\cup_i\{\upsilon_i+\De^{(i)}\}$
for $\upsilon_i=\upsilon(\ep_i)$,
and  deg$_\o(\m)=\hbox{deg\,}\{\De^{(i)}\}=
\sum_{i=0}^{\rk\!-\!1}$deg\,$(\De^{(i)})$, which is
the codimension of $\m$ in $\o^{\rk}$.
Here deg\,$\De=|\Z_+\setminus \De|$ for $\Ga$\~modules
$\De\subset \Z_+$.

Given a collection $\{\De^{(i)}, 0\le i\le \rk-1\}$ of
standard $\Ga$\~modules, the corresponding $\De$\~{\em cell\,}
$\j\{\De^{(i)}\}=\j_{\r,\rk}\{\De^{(i)}\}$
is a set of all {\em standard\,} modules
$\m\subset \o^{\rk}$ with
$\De(\m)=\{\upsilon_i+\De^{(i)}\}$. If
$\De=\cup_i\{\upsilon_i+\De^{(i)}\}$, then we will sometimes use
the notation $\j_{\rk}(\De)$ instead of $\j\{\De^{(i)}\}$.
It has a natural structure of a quasi-projective
variety. Indeed, such $\m$ can be identified with
$\r$\~invariant subspaces in
$\o^{\rk}/z^{\cc}\o^{\rk}$ of co-dimension
$\hbox{deg\,}\{\De^{(i)}\}$
satisfying the following ``jump-conditions":

We set $\m^{\{r\}}\equal\{y\in \m\,\mid\,
\upsilon(y)\ge r\}$ for $r\in \R_+$, which are $\r$\~submodules
of $\m$. Let $\vep>0$ be a real number smaller than
any  $|\upsilon_i-\upsilon_j|$ for  $i\neq j$. Then the spaces
$\m^{\{r\}}/\m^{\{r+\vep\}}$ modulo $z^{\cc}\o^{\rk}$
must be one-dimensional exactly
at $r\in \De(\m)\setminus \upsilon(z^{\cc}\o^{\rk})$ and $\{0\}$
otherwise.

The following spaces do not depend on the choice
of the valuation:
\begin{align}\label{jrk-def}
&\j[d]\!=\!\j_{\r,\rk}[d]\equal\{\m \mid dev(\m)=d\},\,
\j\!=\!\j_{\r,\rk}\equal\cup_d \j[d],\\
&%\where
dev(\m)\!=\!\sum_{i=0}^{\rk\!-\!1}dev(\De^{(i)}),\,
\De(\m)\!=\!\{\upsilon_i+\De^{(i)}, 0\!\le\! i\!\le\!
 \rk\!-\!1\}.\notag
\end{align}

For our Main Conjecture the
definition of $\j_{\r,\rk}$ as a disjoint union of
quasi-projective
varieties is sufficient. However it is
important to extend the interpretation of the rank-one
$J_\r$ above
as a ``single variety" to any ranks, including
the description of the boundaries of the $\De$\~cells.
This is what we proceed to do now.
\vskip 0.2cm

{\sf Jacobian factor as a single variety.}
For a submodule $\m\subset \o^{\rk}$, let
$\m_\o\!=\!\m\!\otimes_\r\!\o\!=\!\{\sum_j \o
y_j, y_j\!\in\! \m\}$ and
$ \hbox{Aut}_\o(\m)\equal\{ A\in GL(\rk,\o)\,\mid\,
A(\m)=\m\}$. Then
\begin{align}\label{autmc}
A(\m_\o)=\m_\o, \
A(\hbox{\sf C}(\m))=\hbox{\sf C}(\m)\for A\in \hbox{Aut}_\o(\m).
\end{align}
See e.g. \cite{Ar} for constructible sets and algebraic spaces.
\begin{proposition}\label{PROP:LIFT}
(i) For $\si=\cc=\cc_\r$,  let
$\j^{cond}_{\r,\rk}$ be a set of all
$\cc$\~liftable modules
$\tilde{\m}$ of conductor $(z^\cc)\o^{\rk}$. Then
the set of corresponding pure $\m'$ can be identified with
all pure standard $\m$ via the map $s\!t$
and its inverse
$\pi:
\m\ni y=\sum c_i \ep_i\mapsto y'=\sum c_i z^{\cc-n_i}\ep_i$.
\vskip 0.1cm

Let $\m=s\!t(\m')$ for such $\m'$. Then Aut$_\o(\m')=\pi
\hbox{Aut}_\o(\m)\pi^{-1}$ and
the union of $GL(\rk,\o)$\~orbits of pure standard $\m$, which
is $\j_{\r,\rk}$,
can be identified with the union of $GL(\rk,\o)$\~orbits
of their $\pi$\~images, which
is $\j^{cond}_{\r,\rk}$. The latter
is a union of connected constructible sets
corresponding to $\tilde{\m}$ with the same
conductor $(z^\cc)\o^{\rk}$ and fixed deviations
$dev(\tilde{\m})$,
which is compatible with the natural quasi-projective 
structure of the spaces $\j_{\r,\rk}[d]$ for $d=dev(\m)$.
\vskip 0.2cm

(ii)  For $\si=\cc=\cc_\r$,  let
$\j^{\de\cdot\rk}_{\r,\rk}$ be a set of all
$\cc$\~minimal
liftable modules $\tilde{\m}$ such that
$dev(\tilde{\m})\!=\!0$, equivalently,
dim$_\o(\o^{\rk}/\tilde{\m})\!=\!\de\cdot\rk$.
Then $\sum_{i=0}^{\rk-1} d'\lan i\ran=dev(\m)$
for the corresponding
pure $\m=s\!t(\m')$ from Definition \ref{DEF:LIFT},
where the inverse of $s\!t$ (for pure $\m,\m'$)  is now
$\pi_d: \m\ni y=\sum c_i \ep_i\mapsto y'=
\sum c_i z^{d'\lan i\ran}\ep_i$.
\vskip 0.1cm

As above,  Aut$_\o(\m')=\pi_d \hbox{Aut}_\o(\m)\pi_d^{-1}$ and 
the union of $GL(\rk,\o)$\~orbits of pure standard $\m$, which 
is $\j_{\r,\rk}$, can be identified as a set  with the 
union orbits of the $GL(\rk,\o)$\~orbits of their 
$\pi_d$\~images, which is $\j^{\de\cdot\rk}_{\r,\rk}$. The 
latter is a constructible subset in the projective 
variety of all $\r$\~invariant $k$\~vector subspaces in 
$(\o/(z^{\cc}))^{\rk}$ of co-dimension $\de\cdot\rk$, which 
is compatible with the quasi-projective structure of 
$\j_{\r,\rk}[d]$. 
\end{proposition}

{\it Proof.} For $(i)$, we need to check that
Aut$_\o(\m')=\pi\, \hbox{Aut}_\o(\m)\,\pi^{-1}$,
where  $\m=s\!t(\m')$, $\m'=\pi(\m)$, the first is
standard and both modules are pure. The problem is that
$\pi\, GL(\rk,\o)\,\pi^{-1}$ does not generally belong to
$GL(\rk,\o)$. It is the same problem for $(ii)$ with $\pi_d$
instead of $\pi$. However, if the stabilizers are identified,
then the corresponding $GL(\rk,\o)$\~orbits will be identified
too. Let us show that the image of
$\j_{\r,\rk}$ under  $\pi$ (or $\pi_d$) is naturally 
a constructible sets, which is  
compatible with the
natural quasi-projective structure above  of 
 $\j_{\r,\rk}[d]$.

Let us begin with $\pi$. Using (\ref{autmc}),
$\hbox{\sf C}(\m)=z^{\cc}\pi^{-1}(\o^{\rk})$,
which implies that
$A\,\pi^{-1}\,(\o^{\rk})=\pi^{-1}\,(\o^{\rk})$ and
$\pi\,A\,\pi^{-1}\,(\o^{\rk})=\o^{\rk}$.
This gives the required identification. The set
$\j^{cond}_{\r,\rk}$ is naturally a union of projective spaces
minus
some their quasi-projective subspaces formed by
non-liftable submodules; thus they are constructible sets.
This gives $(i)$.
\vskip 0.2cm

This argument (the usage of the conductors)
becomes a bit more technical for $(ii)$.
As above, we need to check that
$\pi_d\, A\, \pi_d^{-1}\in GL(\rk,\o)$, explicitly,
$d'\lan i\ran-d'\lan j\ran+\nu(a_{ij})\ge 0$ for $A=(a_{ij})$.
The indices are arbitrary here, but it suffices to
assume that $i>j$ because $a_{ij}\in \o$. Since
$d'\lan i\ran-d'\lan j\ran+\nu(a_{ij})=
(n_j-n_i)+(n'_i-n'_j)+\nu(a_{ij})$,  the following
inequalities are needed for the proof:
$\nu(a_{ij})\ge (n_i-n_j)-(n_i'-n_j').$ The invariance
of {\sf C}$(\m)$ gives that
$\nu(a_{ij})\ge n_i-n_j$ for $A\in$ Aut$_\o(\m)$, and
the inequalities $n'_i-n'_j\ge 0$ are imposed in $(ii)$ for
$i>j$.
\vskip 0.2cm

We did not use the {\em minimality\,} in this argument. It is
needed only for the {\em uniqueness\,}
of the passage from standard
pure $\m$ to liftable pure $\m'$ of deviation
$\de\cdot\rk\,$, i.e. for making this in terms 
$dev(\m)$ only. Concerning 
the identification of $\j_{\r,\rk}$ with
$\j^{\de\cdot\rk}_{\r,\rk}$,
the latter makes the former
a constructible subset in the projective
variety of $\r$\~invariant vector subspaces
in $(\o/(z^{\cc}))^{\rk}$ of co-dimension $\de\cdot\rk$.
Here we need to remove all non-$\cc$-liftable ones,
which form a certain quasi-projective subvariety (described
by algebraic (in)equalities).
In the case of $\rk=1$, all
modules of $dev=0$ are liftable to standard ones. 
\sq
\vskip 0.2cm
%\vfil

Part $(ii)$ provides a counterpart of the projective
variety $\pi(J_{\r})={}^\pi\!\! J_{\r}$
defined after Lemma \ref{LEM:cond}, though this
is now less canonical.
%\smallskip

First, we need to impose
the {\em purity\,} of the conductors with respect to the
basis $\{\ep_i\}$, and
the resulting identification is only due to the isomorphism
of the corresponding $GL(\rk,\o)$\~orbits.

Second, we need to impose the {\em $\cc$\~minimality condition\,}
to make the passage from $\m$ to $\m'$
depending only on $\m$.  Third, 
$\j^{\de\cdot\rk}_{\r,\rk}$ is a
{\em constructible set\,} now, not a projective or
quasi-projective variety.

As far as the count of points over
finite fields is concerned, $\j_{\r,\rk}$ and
$\j^{\de\cdot\rk}_{\r,\rk}$ give the same. However generally,
we have two different definitions, which can be
identified only ``orbit-wise". In contrast to the former (the main
in this particular paper), $\j^{\de\cdot\rk}_{\r,\rk}$ opens
a road to an important theory of (relative) boundaries
of the natural cells there (for instance, $\De$\~cells),
similar to the Schubert calculus.

\vskip 0.2cm
The space
$\j^{cond}_{\r,\rk}$ is quite interesting too. By construction,
it is a union of pieces (not a connected
variety)  in contrast to
$\j^{\de\cdot\rk}_{\r,\rk}$, even for $\rk=1$.
The conductor of modules is fixed there, which is
$z^{\cc} \o^n$, but the deviations are not fixed.
Correspondingly, the closure here is {\em relative\,},
restricted only to the modules with coinciding conductors.
In $\j^{\de\cdot\rk}_{\r,\rk}$, conductors generally change
upon taking the limits, which cannot be seen in
$\j^{cond}_{\r,\rk}$. However the map $\pi:$
$\j_{\r,\rk}\mapsto \j^{cond}_{\r,\rk}$
is simpler than $\pi_d$ and is
given entirely in terms of the (pure) conductors of $\m$;
$\pi_d$ requires some
special choice of $\{d_i\}$.

%%%!!! ADD AND CHANGE THE FOLLOWING TWO and one above:$
The simplest cases of stable modules $\m$ are {\em minimal\,}
ones and $\o$\~invariant modules. There is only one
standard $\o$\~invariant module, which is $\o^{\rk}$.
Accordingly, its image in
$\j^{cond}_{\r,\rk}$ is $z^{\cc}\o^{\rk}$ and its
 $\pi_d$\~image in $\j^{\de\cdot\rk}_{\r,\rk}$ is
 $z^{\de}\o^{\rk}$. Minimal modules remain unchanged under
 both projections.

%\vskip 0.2cm
\vfil
Following \cite{PS}, let us outline the proof
of the {\em connectivity\,} of the closure
$\overline{\j}^{\de\cdot\rk}_{\r,\rk}$
of $\j^{\de\cdot\rk}_{\r,\rk}$ in
$Gr(\de\cdot\rk, (\o/z^{\cc}\o)^{\rk})$ over the field $k$
of zero characteristic (though the latter is not strictly
necessary).
Any $GL(\rk,\o)$\~invariant
{\em closed\,} subset in $\overline{\j}^{\de\cdot\rk}_{\r,\rk}$
contains
at least one $\r$\~module fixed by the action
of the upper triangle unipotent group
$$U\!=\hbox{exp}\{a\!=\!(a_{i,j})\!\in\! \hbox{Mat}(\rk,\o)
\mid a_{i,j}\!=\!0\hbox{ for } i\!>\!j,\,
a_{i,i}\!\in\! (z)\}.
$$
We use that the orbits of $U$ are affine spaces.
Note that this action preserves conductors of
{\em pure\,} modules because they are of the form
$\oplus z^{n'_i}\o^{\rk}$ for {\em non-decreasing\,} $\{n_i'\}$.
So it preserves the
space of pure modules and its closure in
$\overline{\j}^{\de\cdot\rk}_{\r,\rk}$.
We consider pure  $\m'\in \j^{\de\cdot\rk}_{\r,\rk}$
only with {\em non-increasing\,} $\{d'\lan i\ran\}$ in
Lemma \ref{DEF:LIFT}, $(b)$. Sufficiently
general elements $A\in U$ transform $\m'$ to
pure $\m''$ with pairwise equal $d''\lan i\ran$.
Recall that $d\lan i\ran$ is the minimal power
of $z$ in the projection of $\m$ onto $\ep_i$.
Therefore, if such an orbit is a single
$U$\~invariant module, then
it can be only $z^\de\o^{\rk}$.  This gives the required
connectivity.
\vskip 0.2cm

{\sf The simplest example.} Let
$\r=k[[z^2,z^3]],\, \rk=2$. Then $\de=1,\cc=2$ and
$\r$\~invariant $2$\~dimensional subspaces
in $4$\~dimensional $\o^2/z^2\o^2$ form $Gr(2,k^4)$,
since the $\r$\~invariance here holds automatically.
It is a quadric in $\mathbb P^5$ over $k$, and is a disjoint union of
the cells
$\mathbb A^4$, $\mathbb A^3$,  $2\mathbb A^2$, $\mathbb A^1$ and
$\mathbb A^0\!=\!pt$. Using these cells,
$\j^{\de\cdot\rk}_{\r,\rk}$ is a
constructible set obtained by
removing $\mathbb A^2\cup\mathbb A^1$ from it,
i.e. $\j^{\de\cdot\rk}_{\r,\rk}=\mathbb P^4\setminus
\mathbb A^1$, which
is not a quasi-projective variety.

Indeed, we need to remove one  $GL(2,\o)$\~orbit
of non-liftable $V'_0=\ep_0\o+z\ep_0\o\mod z^2\o^2$, which
is naturally
a union of the following two families   %%%%!!!!
of modules for $\la_0^{1,2}\in k$:
$\tilde{V}_1=(\ep_1+\la_0^1 z\ep_0)\o+z\ep_1\o$ and
$\tilde{V}_2=(\ep_0+\la_0^1\ep_1+\la_0^2 z\ep_1)\o+
(z\ep_0+\la_0^1 z\ep_1)\o$. For $V'_0$, which is pure:
$$
n_0'=0,n_1'=2,\ d'\lan 0\ran=0,
d'\lan 1\ran=2,\ s\!t(V_0)=\ep_0\o+\ep_1\o.
$$
It is non-liftable, because the inequality
$d'\lan 0\ran\ge d'\lan 1\ran$ required for the liftability
is violated. Recall that $\pi_d$ must {\em flatten\,}
the conductor $(z^{n_0})e_0+(z^{n_1})e_1$, which
is not the case. Indeed, $V_0'$ coincides with its
conductor, $n_0=0=n_1$ and therefore $d'\lan 0\ran$ must
coincide with $d'\lan 1\ran$ for the liftability.

The modules $\tilde{V}_{1,2}$ are $z$\~invariant too, so
they coincide with their conductors. However they are not pure
and must be first transformed to $V'_0$ before applying $s\!t$.
Recall that $\{d'\lan i\ran\}$ and $s\!t$ were defined for
{\em any\,}
modules. For instance, $d'\lan 0\ran=1, d'\lan 1\ran=0$ for
$\tilde{V}_1$ and
$s\!t(\tilde{V}_1)=(\ep_1+\la_0^1 \ep_0)\o+ %%%%%!!!!!
z\ep_0\o+ z\ep_1\o$.
This module is not standard and cannot be standard because
the initial pure $V'_0$ is not. Generally, applying $s\!t$
does not make much sense before the passage to pure $\m'$.

Let us give an example of a liftable family.
Let $\tilde{V}_3=(\ep_0+\la_0^1\ep_1+\la_0^2 z\ep_1)\o+
(z\ep_0+\la_1^1 z\ep_1)\o$ with $\la_1^1\neq\la_0^1$.
If $\la_1^1=\la_0^1$, then we arrive at $\tilde{V}_2$.
The substitution is $\ep'_0=\ep_0+\la_0^1\ep_1,
\ep'_1=\ep_0+\la_1^1\ep_1$. The corresponding pure $V'_3$
is $(\ep'_0+\la_0^2 z\frac{\ep'_0-\ep'_1}{\la_0^1-\la_1^1})\o+
z\ep'_1\o$ with $n_0'\!=\!0,n_1'\!=\!1$,\ $d'\lan 0\ran\!=\!0,
d\lan 1\ran\!=\!1$. Its $s\!t$\~lift to a standard one is
$(\ep'_0+\la_0^2 z\frac{\ep'_0-\ep'_1/z}{\la_0^1-\la_1^1})\o+
\ep'_1\o$.

\medskip

\subsubsection{\sf Flagged Jacobian factors}
The definitions closely follow those in the   %%%!!! definitions
case of $\rk=1$. We need to choose a valuation $\upsilon$,
but will later see that this actually does not influence
flagged Jacobian factors up to some canonical identifications.
As above, we assume that
$$0\le \upsilon(\ep_0)=\upsilon_0<\ldots<
\upsilon(\ep_{\rk-1})=\upsilon_{\rk-1}<1.
$$
Recall that $\De(\m)\equal\upsilon(\m)$. A flag
$\vect{\m}=\{\m_0\subset \m_1\subset\ldots\subset \m_\ell\}$ of
standard $\r$\~modules $\m_i\subset \o^{\rk}$  is called
{\em full $\upsilon$\~increasing\,} if
the corresponding {\em $\Delta$\~flag}
$\vect{\De}=\De(\vect{\m})=
\{\De_0\subset \De_1\subset\ldots\subset \De_\ell\}$
 is full increasing for $\De_i\equal\De(\m_i)$, i.e.
by definition satisfies the following:
\begin{align} \label{deflagsrk}
\De_{i}=\De_{i-1}\cup \{g_i\},\
 g_1<g_2<\ldots<g_\ell, \,\,1\le i\le \ell,\, \,
g_i\in \upsilon(\o^{\rk}).
\end{align}

It is called a {\em standard flag\,}
if it is full $\upsilon$\~increasing
and $\m_0$ is standard,
which automatically implies  that all $\m_i$ are standard.
We will call such standard flags simply $\ell$\~flags and write
$l(\vect{\m})=\ell$.
\smallskip

{\sf Jacobian $\rk$\~factors.}
The set $\j^{\ell}\!=\!\j_{\r,\rk}^{\ell}$ of all standard
$\ell$\~flags $\vect{\m}$
will be called {\em rank-$\rk$ flagged Jacobian factor\,}.
One has:
\begin{align}\label{defdeflag}
&\j^{\ell}=\cup_{d\,} \j^{\ell}[d\,],\
\j^{\ell}[d\,]\equal\{\vect{\m}\in
\j^{\ell}\,\mid\, dev(\m_0)\!=\!d\},\\
\j^{\ell}[d\,]&=\cup_{\,\vect{\De},\,dev(\De_0)\!=d\,}
\j^{\ell}(\vect{\De}),\ \,
\j^{\ell}(\vect{\De})\equal\{\vect{\m}\,\mid\, \De(\vect{\m})=
\vect{\De}\}.\notag
\end{align}

For fixed $d$ and/or $\vect{\De}$, these spaces are naturally 
quasi-projective varieties. One can extend Proposition
\ref{PROP:LIFT} to introduce a counterpart of
$\j^{\de\cdot\rk}_{\r,\rk}$ \, for arbitrary $\ell>0$ as
a ``connected space", which is set-theoretically
isomorphic to $\j_{\r,\rk}^{\ell}$. We will omit details;
a proper generalization of Proposition \ref{NESTED} is
needed here.

\subsection{\bf Main Conjecture}
We will now state the Main Conjecture.
The definition of
$\j^{\ell}_{\r,\rk}$ and all constructions above were
for any $\r$, however the geometric interpretation
of DAHA superpolynomials is conjectured (and known in
quite a few examples) only
for plane curve singularities. Though see \cite{ChD1}
concerning {\em pseudo-algebraic\,} knots, where the
geometric superpolynomials can be expected too.

\subsubsection{\sf Plane curve singularities}
Let us restrict ourselves to subrings $\r\subset\o$ with two
(algebraic) generators over the base field $k$.
For such {\em plane curve singularities\,},
$\r$ is Gorenstein, $\cc=2\de$ and also there is an isomorphism
$\Z_+\setminus\Ga\ni g\mapsto \cc-g\in
 \Ga\setminus \{\cc+\Z_+\} $.
We mention that this relation provides an explicit formula for
$\De(M^\ast)$, where
$M^\ast\equal\{ y\in \k\,\mid\, yM\in \r\}$. It is given in
terms of $\De(M)$. The simplest example is for invertible
modules; then
$M^\ast$ is invertible and ``$\ast$"\, is the involution
$x\mapsto -x$
of the corresponding generalized Jacobian.

We will not use the dual modules in this paper; see e.g.
\cite{GM1}. However this is an important feature and let
us provide here at least their definition
in the case of any ranks.

We need a non-degenerate form in $k^{\rk}$; the natural
choice is $(\ep_i,\ep_j)=\de_{i,j}$. It was already used in
(\ref{de-proj}) in a similar context:
$\m^{\lan i\ran}=(\m,\ep_i)$.
Then we extend $(\cdot,\cdot)$ to $\k^{\rk}$ and set
$\m^{\ast}
=\{\tilde{y}\in\k^{\rk}\,\mid\,
(\tilde{y},\m)\in \r\}$. This module contains
$\oplus_{i=0}^{\rk-1}(\m^{\lan i\ran})^\ast\ep_i$
and its conductor is exactly $\o^{\rk}$ for any
standard $\m$.
Assuming that $\m$ is pure with $\mathbb C(\m)=
\oplus_{i=0}^{\rk-1}(z^{n_i})\ep_i$, one has:
 $\m^\ast\subset\oplus_{i=0}^{\rk-1}z^{-n_i}\o\ep_i$.

We find it convenient to
use below the realization of $\o^{\rk}$
as the ring extension $\tilde{\o}=\o[z^{1/\rk}]$
(and  $\tilde{\k}=\k(z^{1/\rk})$). Then naturally,
\begin{align}\label{ramz}
\ep_i=z^{i/\rk},\ \upsilon(\ep_i)=i/\rk,\  0\le i\le \rk-1.
\end{align}
This interpretation is not really necessary in this
particular paper. However it clarifies sometimes the
nature of our considerations. For the $\ast$\~duality,
one can take  $(\tilde{x},\tilde{y})=
tr_{\k}(\tilde{x} \tilde{y})$, the trace of $\tilde{\k}/\k$,
where $\tilde{x},\tilde{y}\in \tilde{\k}$. For such a choice
of the form $(\cdot,\cdot)$, the definition of $\m^{\ast}$
generalizes the classical definition of the
{\em different ideal\,} of an extension of local fields.
\vskip 0.2cm

\subsubsection{\sf Changing the valuation}
As above, $k$ will be an arbitrary field (unless stated
otherwise). We fix a basis $\{\ep_i, 0\le i\le \rk-1\}$
in $k^{\rk}$ and the evaluation $\upsilon$ extending
$\nu$ ($z$\~valuation) in $\k$, uniquely determined by
$\upsilon_i=\upsilon(\ep_i)$. As in the previous sections,
we always assume that
$0\le \upsilon_0<\ldots<\upsilon_{\rk-1}<1$. Obviously,
the modules $\De^{\lan i\ran}$ (projections onto $e_i$)
are simply permuted accordingly. However
the transformations of the corresponding $\Ga$\~modules
$\De^{(i)}$ upon such permutations can be nontrivial.
This actually does not influence rank-$\rk$ flagged Jacobian
factors too much; they can be identified for
different choices of the valuation. We will provide
the following lemma, which explains how this can be done.

\begin{lemma}\label{LEM:perm}
The construction of the flagged Jacobian factors
does not depend on the particular choice of $\{\upsilon_i\}$
chosen as above.
Furthermore, any flag $\vect{\m}$
of modules satisfying inequalities
(\ref{deflagsrk}) for $g_i$ only within $\Z$\~orbits
(i.e. only when $g_i-g_j\in \Z$) can be
canonically transformed
to an {\sf increasing flag\,}, where $g_i<g_j$
for any $i<j$.
\end{lemma}
{\it Proof.} The first claim is obvious.
Let us justify the second. As above, $\vect{\De}$ will be
the $\De$\~flag of $\vect{\m}$:
$\De_{i}=\De_{i-1}\cup \{g_i\}$. However now we allow
$g_i$ to be added to the corresponding $\De$ with respect to
some permutation $w=\{i_1,i_2,\ldots\}$ of their indices,
assuming that the
initial (increasing) order is preserved in $w$
within all $\Z$\~orbits.

Let $I$ be a set $I=\{i_\circ,\ldots,i_{\bullet}\}$
of consecutive  indices in $w$
such that $g_i-g_j\not\in \Z$ for $i\neq j\in I$
and  maximal in the following
sense. Let $i_{in}=i_\circ-\!1,\, i_{out}=i_{\bullet}+\!1$,\
$g_{in}=g_{i_{in}},\, g_{out}=g_{i_{out}}$. We assume
that $g_{i}-g_{in}$ and
$g_{out}-g_{j}$ are in $\N$ for some $i,j\in I$. Given
such $I$, let us transform $\vect{\m}$ to make $g_i$
increasing.

We change $\upsilon_i$
(they can become beyond $[0,1]$) to ensure the inequalities
$|g_i-g_j|<1$ for all $i,j\in I$. Then automatically
$g_{in}< I<g_{out}$.

Let us now sort the indices of $I$:
$$
i'_\circ=\hbox{min}(I), \
(i_{\circ}+1)'=\hbox{min}(I\setminus\{i'_\circ\})
\hbox{\,  and so on till \,} i'_{\bullet}=\hbox{max}(I).
$$
Accordingly,
$\De'_{i'_\circ}\!=\!\De_{in}\cup\{g_{i'_\circ}\}$,\,
$\De'_{i'}\!=\!\De'_{(i-1)'}\cup\{g_{i'}\}$,\,
$\De'_{i'_\bullet}\!=\!\De_{i_\bullet}$.

We are going to modify $\vect{\m}$ to obtain $\vect{\m}'$
with $\vect{\De}'=\{\De'_{i'}\}$.
 The modules $\m_{in}$ and $\m_{out}$ for $i_{in},i_{out}$
(and those before the former and after the latter)
will remain unchanged. Using that $|g_i-g_j|<1$ in $I$, let us
switch from $\{\upsilon_i=\upsilon(\ep_i)\}$ to the
opposite one
$\{\upsilon_{i}^{op}=\upsilon_{\rk-i-1}\}$. I.e.
the valuation $\upsilon$ is changed by the opposite one; now
$\upsilon^{op}_{i}>\upsilon^{op}_{j}$ if $i<j$. We will use
$\upsilon^{op}_i$ only within $I$.

Then filtration in $\m_{out}$ corresponding to
$\upsilon^{op}$ is as follows:
$\m_{out}^{\{g\}_{op}}=\{y\in \m_{out}\mid
\upsilon^{op}(y)\ge g\}$. Using it, we set:
\begin{align*}
&\m'_{i'}\equal\m_{out}^{\{g_{i'}\}_{op}}\for i'_\circ\le i'\le
i'_\bullet\,.
\end{align*}

This first module here is $\m'_{i'_\circ}$.
It is a one-dimensional $\r$\~submodule in
$\m_{out}/\m_{in}$ with the $\De$\~set
$\De'_{i'_\circ}=\De_{in}\cup\{g_{i'_\circ}\}$.
The existence of such modules for any $\De_0\cup \{g_i\}$
is a general fact; cf.
$\De_0\cup \{g_i\}$ from Proposition
\ref{NESTED}, $(i)$.
However generally they are not unique such.
Using $\upsilon^{op}$, we define
them {\em canonically\,}, which
is possible only because $g_i$ here are all
from different $\Z$\~orbit.

Then we go to the resulting flag $\vect{\m}'$,
find there another $I$ and proceed by induction.
Strictly speaking, this procedure depends on
the way we pick $I$\~sets, but it can be made
canonical upon some combinatorial considerations,
which we omit.
%\sq
%\vskip -0.5cm

\subsubsection{\sf Bad reduction}
One of the key corollaries of our Main Conjecture
\ref{CONJ:MAIN} (below)
is that the geometric superpolynomials are actually
{\em topological\,} invariants of plane curve singularities.
Analytic classification of such singularities is much more
involved than the topological classification (and is not
finished). This triggers important questions concerning
the topological invariance of various constructions here
and in the related theory of affine Springer fibers.
The theory of bad reductions of plane curve singularities
modulo prime numbers up to topological invariance is an
important example.

Considering plane curve
singularities topologically, i.e. up to isotopy
of their links,  one can always find the corresponding
ring $\r$ to be defined over $\Z$. Furthermore,
for {\em any\,} prime $p$, there exists
$\r$ over $\Z$ in the same topological class
such that the corresponding $\Ga$ of
$\r\otimes_{\Z}\mathbf F_p$  coincides with that over $\C$.
Classically, the places of good reduction are prime numbers $p$
where a given manifold remains smooth. This definition
is of course not applicable as such, but at least $\Ga$
must remain unchanged modulo $p$.
%It is not impossible that it is sufficient without flags,
%i.e. for $\ell=0$.

This is the weakest possible definition of good reduction.
We call a prime $p\,$ a place of {\em good reduction
in the strong sense\, }
if there exists $\r'$ over $\Z$ topologically equivalent to
$\r\otimes_{\r}\C$
such that all Piontkowski-type cells
$\j^{\ell}_{\r',\rk}\{\vect{\De}\}$ do not change their
bi-regular type when going from $\C$ to $\mathbb F_p$
(via $\Z$). This is what we use (and check) practically.
%\vskip 0.2cm

We note that
\cite{ChRH} hints that the weak definition can be
insufficient for flags. Heuristic (indirect)
arguments there indicate  that $p=2$ can be a place of bad
reduction for $\r=\Z[[z^4,z^6+z^7]]$ for $1$\~flags ($\ell=1$)
and $\rk=1$, but we have not checked this so far. Actually,
the ring $\r$ has bad reduction modulo $p=2$, but
$\r'=\Z[[z^4+z^5,z^6]]$
of the same topological type has a good reduction
at $p=2$.  Presumably,
$p=2$ is always a place of bad reduction for {\em $1$\~flags\,}
for any $\r'$ over $\Z$ of the same topological type as $\r$.
Actually,
this example is more important for \cite{ChRH} than for our
present paper; we simply stick to the strong
understanding of bad reduction in the following conjecture.

\subsubsection{\sf Geometric superpolynomials}
The notation is from this and previous sections;
$\r\subset \o=\C[[z]]$ is a ring of a given unibranch
plane curve singularity $\c$, which will be
considered over $\Z$ in $(ii)$ below.

The corresponding link of $\c$
is given by the sequence of Newton pairs
$\vec{\rr},\vec{\ss}$. In terms of the latter,
$\h_{\,\vec\rr,\,\vec\ss}\,(\om_{\rk};q,t,a)$
is the DAHA superpolynomial for the Young diagram
corresponding to $\om_{\rk}$, the $\rk$\~column,
which is a topological invariant of the link of $\c$.
For standard submodules,
$\rk$ is their rank and
$\ell=0,1,\cdots$ is the length of the flags
of such modules, which are
$\vect{\m}=\{\m_0\subset\ldots\subset\m_{\ell}\}$
satisfying conditions (\ref{deflagsrk}).

\begin{conjecture}\label{CONJ:MAIN}
Let $\r\subset \o$ be defined over $\Z$.
The corresponding rank-$\rk$ flagged Jacobian factor
$\j^{\ell}_{\r,\rk}$ is a union of quasi-projective
varieties $\j^{\ell}_{\r,\rk}[d]$, where $d$ is $dev(\m_0)$
for standard flags $\vect{\m}$.

(i) Considering   $\j^{\ell}_{\r,\rk}[d]$ over $\C$, by
$H^{\!BM}_{j}(\j^{\ell}_{\r,\rk}[d])$ we mean
{\sf Borel-Moore homology\,}, i.e.
relative singular homology of
$\j^{\ell}_{\r,\rk}[d]$ compactified by one point $pt$
(relative with respect to $pt$).

We conjecture that the odd homology
$H^{\!BM}_{2i+1}(\j^\ell_{\r,\rk}[d])$
vanishes for all $i,d\ge 0$ and that
$\h_{\,\vec\rr,\,\vec\ss}\,(\om_{\rk};q,t,a)$
coincides with the {\sf singular
superpolynomial\,} defined as follows:
\begin{align}\label{conjcoh}
\h^{sin\!g}_{\r,\rk}(q,t,a)\equal
\sum_{d,i,\ell}\hbox{dim}_{{}_\R\!}
\bigl(H^{\!BM}_{2i}(\j^\ell_{\r,\rk}[d])\bigr)
q^{d+\ell}\, t^{\delta\cdot\rk^2-i} a^\ell.
\end{align}
\vskip 0.2cm

(ii) Let prime $p$ be a place of good reduction
in the strong sense for $\r\subset \o$ over $\Z$.
Accordingly,  the varieties $\j^{\ell}_{\r,\rk}[d]$
will be now considered as schemes (and then reduced) over
$\mathbb F=\mathbb F_{p^m}$ for $m\in \N$.  The number of
$\mathbb{F}$\~points of a scheme $X$ over $\mathbb{F}$
will be denoted by  $|X(\mathbb{F})|$. The
{\sf motivic superpolynomial\,} is:
\begin{align}\label{conjmot}
\h^{mot}_{\r,\rk}(q,t,a)\equal
t^{\delta\cdot\rk^2}\sum_{d,\ell}\,\mid\!
\j^\ell_{\r}[d](\mathbb{F})\!\mid\,
q^{d+\ell} a^\ell.
\end{align}
Letting $1/t=p^{m}$, we conjecture that
$\h^{mot}_{\r,\rk}(q,t,a)=
\h_{\,\vec\rr,\,\vec\ss}\,(\om_{\rk};q,t,a).$
\end{conjecture}
\vskip -0.5cm
\sq

The range of the indices in $(i)$  is
 $\,0\le d\le \de\cdot \rk$,\,
$0\le i \le \de\cdot \rk^2\,$ and  $\,\ell\,$  %%%!!!\rk^2
is from $0$ to
$\,\rk\cdot(\ss_1\rr_2\cdots \rr_{\ell}-1)$, where
we assume that $\ss_1<\rr_1$. It is the same in $(ii)$
for $d,\ell$. This can be readily seen for
$\h^{sin\!g}_{\r,\rk}, \h^{mot}_{\r,\rk}$, and is not
difficult to check for $\h_{\,\vec\rr,\,\vec\ss}$.
The product $\ss_1\rr_2\cdots \rr_{\ell}$ here is
the multiplicity of the singularity.

If the varieties $\j^{\ell}_{\r,\rk}[d]$ are paved
by affine spaces over $\C$ and over $\mathbb F$, then
conjectures in $(i)$ and $(ii)$ are essentially
equivalent. Many Piontkowski-type cells
$\j^{\ell}_{\r,\rk}(\vect{\De})$ are affine spaces.
They are all affine spaces for torus knots. However in $\rk>1$,
all non-torus $\r$ (non-quasi-homogeneous
singularities) starting with $\r=\C[[z^4,z^6+z^7]]$
always have such cells that are not affine spaces
in all examples we considered.
This does not mean that $\j^{\ell}_{\r,\rk}[d]$
cannot be paved by affine spaces; we do not have
counterexamples and the structure of the Piontkowski-type cells
indicates that this is not impossible.

%%%NEW:

\setcounter{equation}{0}
\section{\sc Torus knots}
The purpose of this section is to prove that
$\j^{\ell}_{\r,\rk }(\vect{\De})$ are affine spaces
and compute their dimension in the case $\r=k[[z^p,z^q]]$,
where $p$ and $q$ are coprime and $p<q$. This can be used
to verify Conjecture \ref{CONJ:MAIN}, but this is
not done at the moment for $\rk>1$ (see \cite{Mel}
for $\rk=1$).
We note that
the affineness of $\j^{\ell}_{\r,\rk }(\vect{\De})$
can be justified using the theorem of 
Bialynicki-Birula for varieties with a $\C^*$\~action.

\subsection{\bf Main claims}
We will use
the conventions from (\ref{ramz}) with minor
modifications. For instance, we use $p,q$ instead of 
$\rr,\ss$ above and multiply the valuation
and all $\De$ by $\rk$ to simplify notations.

\subsubsection{\sf Basic definitions}
We will consider rank-$\rk$ $\r$-modules to be
embedded in the ring
$\tilde{\o}=k[[z]](z^{1/\rk})$ and 
{\it the valuation
$\upsilon(z^{j/\rk})=j$ will be used\,}, 
so our $\Delta$ will be rank-$\rk$
$\tilde{\Gamma}$\~modules over
$\tilde{\Gamma}\equal\rk\cdot \Gamma=\upsilon(\r)$. This
valuation is the
one from (\ref{ramz}) multiplied by $\rk$.

All modules $\m$ will be assumed standard of rank $\rk$.
Recall that $\m$ is {\em standard\,}
if $\tilde{\o}\m=\tilde{\o}$, equivalently,
$\De=\De(\m)=\upsilon(\m)$ is standard, i.e. $\De$
contains $0,1,\cdots,\rk-1$.
For a standard $\tilde{\Gamma}$\~module $\De$:
\begin{align}\label{newdelta}
&\Delta=\medcup\ _{i=0}^{\rk -1}\,\bigl(\tilde{\Delta}^{(i)}
\bigr),\ \tilde{\Delta}^{(i)}=\{x\in
\De: x= i \hbox{\, mod\, } \rk \}.
\end{align}

Recall from (\ref{de-comp}), that $\De^{(i)}$ are
the components of $\De$, which are standard rank-one
modules over $\Ga$. One has:
$$ \De^{(i)}=\{(x-i)/\rk\mid x\in  \tilde{\Delta}^{(i)}\}.
$$

Each $\Delta^{(i)}$ has a $p$-\emph{basis} in
the sense
of \cite{Pi}. Explicitly, there exists a unique $p$-tuple
$(a_{i,0},...,a_{i,p-1})$ such that
\[
\Delta^{(i)}=\medcup\,_{j=0}^{p-1}\,(a_{i,j}+p\N),
\hspace{1mm} a_{i,j}= jq\hbox{\, mod\,} p,\, \where
0\leq j\leq p-1.
\]
We call the set
\begin{align}\label{pbas}
\{\rk\cdot a_{i,j}+i\,\mid\,
0\leq i\leq \rk -1,0\leq j\leq p-1\}
\end{align}
the $p$-\emph{basis of} $\De$; notice the multiplication
by $\rk$.

We will use the definition of the quasi-projective varieties
$\j_{\r,\rk}^{\ell}(\De)$ from Section \ref{SEC:JACF}.
Since the ring remains the same through this section, we mostly
omit $\r$ here; also we set $\j_{\rk}(\De)=$
$\j_{\r,\rk}^{\ell=0}(\De)$.

\subsubsection{\sf Theorems \ref{unflag},\ref{flag}}
Let
$
\gamma_{\De}(b)\equal|[b,\infty)\backslash\De|
$
be the \emph{gap counting function}. Here the interval
$[b,\infty)$ is a subset of $\N$. We need
$\gamma_{\De}$ in order
to state the following theorems, to be proven later.

\begin{theorem}\label{unflag}
Let $\De$ be a
rank-$\rk$ module over $\tilde{\Gamma}=
\upsilon(k[[z^p,z^q]])=\rk\cdot\bigl(\Z_+p+\Z_+q\bigr)$.
Then $\j_{\r}(\De)$
$=\j_{\r,\rk}^{\ell=0}(\De)$ is biregular 
to the affine space $\mathbb{A}^{\sum_{i=0}^{\rk -1} N_i}$,
where
\[
N_i=\sum_{j=0}^{p-1}(\gamma_{\De}(\rk \cdot a_{i,j}+i)-
\gamma_{\De}(\rk \cdot a_{i,j}+i+\rk\cdot q)).
\]
\end{theorem}
\vskip -1.0cm\sq

\comment{
\rmk
Note that the formula in Theorem \ref{unflag} is a
``summed" version of Theorem 12 from \cite{Pi}, where the
sum is over the $\tilde{\Delta}^{(i)}$. This is almost true
except that the
function $\gamma_{\De}$ counts gaps that may be in other
$\tilde{\Delta}^{(j)}$ so the dimensions here are generally
larger than just
naively summing the dimensions in the rank-$1$ case. When
$\rk =1$ we recover Piontkowski's theorem.\sq
\vskip 0.2cm
}

If $g\in \N\setminus \De$ and $\De\cup\{g\}$ is
$\tilde{\Gamma}$-module, we define
\[
\mu_{\De,g}\equal\ \text{dim}\left(\j_{\rk }^1
\{\De,\De\cup\{g\}\}\right)-
\text{dim}\left(\j_{\rk }(\De)\right)
\]
to be the \emph{dimension change}. The following theorem
is the flagged generalization  of Theorem
\ref{unflag}.

\begin{theorem}\label{flag}
Let $\vect{\De}=\{\De_i\}$ be a rank-$\rk $ standard flag of
length $\ell$ defined in
(\ref{deflagsrk}) but now with respect to $\upsilon$ above and
$\tilde{\Ga}$.
Then $\j^{\ell}_{\r,\rk }(\vect{\De})$ is biregular
isomorphic to an affine space $\mathbb{A}^N$ with
\[
N=\text{\emph{dim}}\left(\j_{rk}(\De_0)\right)
+\sum^{\ell-1}_{i=0}\mu_{\De_i,g_{i+1}}, \hbox{ where:\,\, }
\]
\[
\mu_{\De,g}=\gamma_{\De'}(g)-\gamma_{\De'}(g+\rk \cdot p)-
(\gamma_{\De'}
(g+\rk \cdot q)-\gamma_{\De'}(g+\rk \cdot p+\rk \cdot q)).
\]
Here $\De=\De_i,g=g_{i+1}$; more generally,
$\De$ is any rank-$\rk $ $\tilde{\Gamma}$-module
such that $\De'=\De\cup\{g\}$ is also a
$\tilde{\Gamma}$-module.
\end{theorem}

\subsubsection{\sf Reduction algorithm}
The following algorithm plays an important role in 
our proof of the theorems.
Given the $p$-basis for $\De$ as
in formula (\ref{pbas}), we define the following elements of
$k[[z]](z^{1/\rk})$:
\begin{align}\label{gen}
m_{i,j}=z^{\frac{\rk \cdot a_{i,j}+i}{\rk}}
+\sum_{\substack{\rk \cdot a_{i,j}+i+k\not\in \De
}}\lambda_k^{i,j}z^{\frac{\rk \cdot a_{i,j}+i+k}{\rk}},
\end{align}
where $ 0\leq i\leq \rk-1$, $0\leq j\leq p-1$
and the $\lambda_k^{i,j}$ are variables (the cardinality 
of the set of $\lambda_k^{i,j}$ is
$\sum_{i,j}\gamma_{\De}(\rk \cdot a_{i,j}+i)$).

One can show that for all standard modules there
exists a set
of generators of the same form as the $m_{i,j}$.
This results from the following reduction algorithm.

Suppose $y\!\in\! k[[z]](z^{1/\rk })$
and let $y_0\!=\!y$.
Then define inductively,
$y_{k+1}\!=\!y_k$ if $k\not\in \De$ and set
$s_k\!=\!0$.
If $k\in \De$, then find the monomial, $c_kz^{\frac{k}{\rk }}$,
with power $k/\rk$ in
$y_k$ and the element
$s_k\in \r$ such that $s_km_{i,j}\!=\!c_kz^{\frac{k}{\rk }}\!+\!
\ldots$
for one of the
generators $m_{i,j}$ (it may not be unique). Then we
let $y_{k+1}\!=\!y_k\!-\!s_km_{i,j}$. Finally, the sequence of
elements $y_k\in \o$ converges to an element $y_{\infty}$,
which has the form
$\sum_{l\in \N\setminus\De}d_\ell z^{l/\rk }$
for some coefficients $d_k$.

In general, the result of the reduction algorithm will depend
on the choices made.
We will make it canonical as follows.
For $k\in \De$, there is a unique choice of $i,j$ and $l$ such that
$k=\rk\cdot p\cdot l+\rk\cdot a_{i,j}+i$.
Letting $s_k=(z^p)^l$,
the reduction algorithm becomes canonical. Only this version will
be used below.

\begin{definition}
For any element $y\in k[[z]](z^{1/\rk })$, let $R_{\De}(y)$ be
the result of the reduction algorithm above
applied to $y$ with the conventions from the previous
paragraph. In this case $R_{\De}$ is a
$\mathbb{C}$-linear projection,
i.e. $R_{\De}(R_{\De}(y))=R_{\De}(y)$.
\end{definition}

Later, we will need $R_{\De}$ for length-one flags (pairs)
$\{\De_0\subset\De_1\}$.
In this case, we will use the notation
$R_i\equal R_{\De_i}$, $i=0,1$.

\subsection{\bf Proof of Theorem \ref{unflag}}
It is similar to the proof
of the rank-one case in \cite{Pi}. Let $\m$ be the module
generated by the elements $(\ref{gen})$. It is
necessary to characterize when $\upsilon(\m)=\De$.
To this end we need
to study the syzygies of the leading terms of the $m_{i,j}$.

\subsubsection{\sf Syzygies of
\texorpdfstring{{\mathversion{bold}$k[{\Gamma}]$}}{Gamma}
-modules}
Let $k[{\Gamma}]$ be the $k$-span of the set
$\{z^{\gamma}\mid \gamma\in \Gamma\}$. Similarly,
we define $k[\De]$ as the $k$-span of
$\{z^{\frac{l}{\rk}}\mid l\in \De\}$; it is a module over
$k[\Gamma]$ via multiplication.
The following Lemma extends Lemma 4 in \cite{Pi} to
rank-$\rk $ modules.

\begin{lemma}\label{syz}
For a rank-$\rk$ standard $\tilde{\Ga}$\~module $\De$, let
\[
A=\bigl(z^\frac{{\rk \cdot a_{0,0}}}{\rk},...,
z^\frac{{\rk \cdot a_{\rk-1,p-1}+p-1}}{\rk}\bigr)
\in k[\De]^{\rk\cdot p}
\]
be a
generating set of a $k[{\Gamma}]$\~module $k[\De]$ ordered
such that the valuations are increasing.
Then there exists
a minimal generating set of syzygies of $A$
consisting of vectors of the form
\[
v=(0,...,0,z^{\gamma_i},0,...,0,-z^{\gamma_j},0,...,0)\in
k[\Gamma]^{\rk \cdot p}
\]
with $A\cdot v=0$, where $\cdot$ is the usual dot product.
\end{lemma}
{\it Proof.}
Since $\De=\cup\tilde{\Delta}^{(i)}$ we have that
\[
k[\De]=\oplus_{i=0}^{\rk-1} k[\tilde{\De}^{(i)}].
\]
Using the connection
$\tilde{\Delta}^{(i)}=\rk\cdot\Delta^{(i)}+i$, we apply
Lemma 4 of \cite{Pi} to
$k[\Delta^{(i)}]$
to obtain a minimal basis of syzygies and then
modify the powers accordingly to
receive a minimal basis of syzygies
of $k[\tilde{\Delta}^{(i)}]$. From the direct sum decomposition
above, we then
have a minimal basis of syzygies of $k[\De]$. \sq
\vskip 0.2cm

\begin{definition}\label{DEFINITX}
For $\vec{r}=(r_{i,j})\in \r^{\rk \cdot p}$, let
 $\sigma=\min\{\upsilon(r_{i,j})+\rk \cdot a_{i,j}+i\}$. The
{\sf initial vector} $in(\vec{r})$ is as follows:
$in(\vec{r})=(\ze_{i,j})$
with $\ze_{i,j}$ equal to the monomial of lowest degree in
$\,r_{i,j}\,$ if $\,\nu(r_{i,j})+\rk \cdot a_{i,j}+i=
\sigma\,$ and
$\,0\,$ otherwise. \sq
\end{definition}
The proof of the following lemma is essentially from \cite{Pi}.

\begin{lemma}\label{PIO}
Let $\m$ be an $\r$-module generated by
$\{m_{i,j}\}$,
%$V\subset \bigoplus_j \r(a_j)$ be
$V$ be a subset of $\oplus_{k=0}^{\rk \cdot p} \r$
such that the initial vectors
$\{\text{in}(\vec{r})\mid \vec{r}\in V\}$ of $V$ linearly generate 
the {\sf syzygies}
of the set $\{z^\frac{{\rk \cdot a_{i,j}+i}}{\rk}\}$ of $k[\Delta]$.

Then  $\upsilon(\m)=\Delta$ if and only if for each
$\vec{r}=(r_{i,j})\in V$ the initial terms $\,in\,$
$\sum r_{i,j}m_{i,j}$ cancel, i.e.,
$\upsilon(\sum r_{i,j}m_{i,j})>\sigma$ and for every $i$ and $j$,
there exists $p_{i,j}\in \r$ such that
$\upsilon(p_{i,j}m_{i,j})>\sigma$ and $\sum r_{i,j}m_{i,j}=
\sum p_{i,j}m_{i,j}$.
If such $p_{i,j}$ exist, then the element
$\sum p_{i,j}m_{i,j}\in \m$
is called a {\sf higher order expression} for
$\sum r_{i,j}m_{i,j}$.\sq
\end{lemma}

The following properties of $R_{\Delta}$ from \cite{Pi}
can be readily extended to the rank-$\rk$ case. First,
 $R_{\Delta}(y)=0$ if and only if $y\in \m$. Second, if
$R_{\Delta}(y)=0$, we may
write $y+f=0$ since $R_{\Delta}(y)=y+f$ for some $f\in \m$. In
this case $-f$ is a higher order expression for $y$.

\subsubsection{\sf
\texorpdfstring{{\mathversion{bold}$\j_{\rk}(\De)$}}{Jacobian
factor}\ is affine}
By combining Lemma \ref{syz}, Lemma \ref{PIO}, and Proposition 5
from \cite{Pi}, the coincidence $\upsilon(\m)=\Delta$ holds
if and only if
\begin{align}
T^{i,j}&\equal z^qm_{i,j}-z^{(\alpha_{i,j+1}-\alpha_{i,j})p}
m_{i,j+1}=
\sum_{k=1}^{\infty}c^{i,j}_kz^
{\frac{\rk \cdot a_{i,j}+i+\rk \cdot q+k}{\rk }}, \\
T^{i,p-1}&\equal z^qm_{i,p-1}-z^{(q-\alpha_{i,p-1})p}m_{i,0}=
\sum_{k=1}^{\infty}c^{i,p-1}_kz^
{\frac{\rk \cdot a_{i,p-1}+p-1+\rk \cdot q+k}{\rk }}\notag
\end{align}
for $j<p-1$ have higher order expressions.
Here $\alpha_{i,j}$ is defined as the unique integer satisfying
$a_{i,j}=jq-\alpha_{i,j}p$. Note that we have
\[
c^{i,j}_k=\chi(\rk \cdot a_{i,j}+i+k)\lambda^{i,j}_k-
\chi(\rk \cdot a_{i,j+1}+i+k)\lambda^{i,j+1}_k,
\]
where $\chi$ is the characteristic function of
$\mathbb{N}\backslash \Delta$.

As previously mentioned, we may use $R_{\Delta}$ to obtain
higher order expressions for $y\in\m$. In particular, we must have 
$R_{\Delta}(T^{i,j})=0$, since $T^{i,j}\in \m$. Following
Piontkowski, let $\tilde{c}^{i,j}_k$ be the coefficients of
$R_{\Delta}(T^{i,j})$. The reduction
algorithm changes $c^{i,j}_k$ only by adding a
polynomial in the variables $\lambda^{i,j}_{l}$ for $l<k$.
If either $\chi(\rk \cdot a_{i,j}+i+k)=0$ or
$\chi(\rk \cdot a_{i,j+1}+i+k)=0$, then we have that
$z^{(\rk \cdot a_{i,j}+i+\rk \cdot q+k)/\rk }$
is eliminated during the reduction algorithm. Indeed,
\begin{align*}
\rk \cdot a_{i,j}+i+k+\rk \cdot q&=\rk \cdot
a_{i,j+1}+i+k+\rk \cdot(\alpha_{i,j+1}-\alpha_{i,j})p, \\
\rk \cdot a_{i,p-1}+i+k+\rk \cdot q&=
\rk \cdot a_{i,0}+i+k+\rk \cdot (q-\alpha_{i,p-1})p.
\end{align*}
We conclude then that the coefficients of $R_{\Delta}(T^{i,j})$
have the form
\begin{align*}
\tilde{c}^{i,j}_k=\lambda^{i,j}_k-
\lambda^{i,j+1}_k+\mbox{ polynomial in }
\lambda_{l}^{i,j}\mbox{ for }l<k.
\end{align*}

Equations $\tilde{c}^{i,j}_k=0$ are necessary and sufficient
for $\upsilon(\m)=\Delta$; their number is
$\sum_{i,j}\gamma_{\Delta}(\rk \cdot a_{i,j}+i+\rk \cdot q)$.
Since $c(\Delta)=\max_{i,j}\{\rk \cdot a_{i,j}+i-
\rk \cdot (p-1)\}$, one has
$R_{\Delta}(T^{i,j})= 0$ for at least for one pair of indices. 
Therefore some $\lambda^{i,j}_k$ can be explicitly solved for  in
terms of a certain collection
of free $\lambda^{i,j}_{l}$. We obtain that $\j_{\rk}(\De)$ is the 
graph of
a regular function on an affine space and
therefore it is isomorphic to affine space. Its dimension is
\[
\sum_{i=0}^{\rk -1}\sum_{j=0}^{p-1}
(\gamma_{\Delta}(\rk \cdot a_{i,j}+i)-
\gamma_{\Delta}(\rk \cdot a_{i,j}+i+\rk \cdot q)).
\]\vskip -0.5cm\sq

\subsection{\bf Proof of Theorem \ref{flag}}
Let $\vect{\De}$ be a rank-$\rk$ standard flag 
of length $\ell$ for $\tilde{\Gamma}$.
The key is the justification of this theorem for
$\vect{\De}$ of length $1$ i.e. for $\vect{\De}=
\{\Delta_0\subset \Delta_1\}$,
where $\Delta_1=\Delta_0\cup\{g\}$. Then iterations will
be used to obtain
the general formula.
Our proof follows \cite{ChP} (in the rank-$1$ case).

We begin by noting that there exist integers $u,v$ such that
$g+\rk \cdot p=\rk \cdot a_{u,v}+u$; indeed,
$\Delta_1$ is a module over $\tilde{\Gamma}$.
Therefore the $p$-basis for $\Delta_1$ can be obtained from the
$p$-basis for $\Delta_0$ by replacing $\rk \cdot a_{u,v}+u$
with $g$.

In addition to the $m_{i,j}$ from (\ref{gen}), we introduce
the following
element of $k[[z]](z^{1/\rk })$ corresponding to $g$:
\begin{align*}
h\equal z^{g/\rk }+\sum_{\substack{g+k\not\in \Delta}}
\lambda_{k}^hz^{\frac{g+k}{\rk }}.
\end{align*}
As above, $\lambda^h_k$ are considered variables.

In the case of flags, the module generated by
\[
S=\{m_{i,j}:0\leq i\leq \rk -1,0\leq j\leq p-1\}
\]
will be denoted by $\m_0$. Let $\m_1$ be the module
generated by
\[
(S\backslash\{m_{u,v}\})\cup\{h\}.
\]

The conditions required for these two modules $\m_0, \m_1$
to form an element of $\j^1_{\rk}(\vect{\De})$ are as follows:
$\upsilon(\m_0)=\Delta_0$, $\upsilon(\m_1)=\Delta_1$, and
$\m_0\subset \m_1$. Here
$\upsilon(\m)=\Delta_0$ and $\upsilon(\n)=\Delta_1$
are, by the last section, equivalent to $R_0(T^{i,j})=0$ and
$R_1(Q^{i,j})=0$, where
\begin{align*}
Q^{i,j} =
\begin{cases}
T^{i,j} &
\parbox{10em}{if $g+\rk \cdot p$ is not equal to
$\rk \cdot a_{i,j}+i$
or $\rk \cdot a_{i,j+1}+i$, } \\
z^qh-z^{(\alpha_{i,j+1}-\alpha_{h})p}m_{i,j+1} &\mbox{if }g+
\rk \cdot p=\rk \cdot a_{i,j}+i, \\
z^qm_{i,j}-z^{(\alpha_{h}-\alpha_{i,j})p}h &\mbox{if }g+
\rk \cdot p=\rk \cdot a_{i,j+1}+i, \\
\end{cases}
\end{align*}
where $j<p-1$ and $\alpha_h=\alpha_{u,v}+1$.
For $j=p-1$, $Q^{i,p-1}$
is defined similarly (see $T^{i,p-1}$ above).

We will show below
that when $\m_0\subset\m_1$ the equations resulting from the
coefficients of
$R_0(T^{i,j})=0$ are actually redundant once we 
consider the equations from coefficients of $R_1(Q^{i,j})=0$.
Hence we just need to find the conditions ensuring
$\m_0\subset \m_1$. For this, we introduce
\begin{align}\notag
F\equal m_{u,v}-(z^p)h=\sum_{k=1}^{\infty}
d_kz^{\frac{\rk \cdot a_{u,v}+u+k}{\rk }}.
\end{align}

\begin{proposition}\label{containment}
Suppose $\m_0$ and $\m_1$ are rank-$\rk$
$\r$-modules defined as above.
Then $\m_0\subset \m_1$ if and only if $R_1(F)=0$.
\end{proposition}
{\it Proof.}
We have $\m_0\subset \m_1$ if and only if $m_{u,v}\in \m_1$.
But this is true if and only if $R_1(m_{u,v})=0$, which holds
if and only if $R_1(F)=0$. Indeed, $F$ is the first step in the
reduction algorithm for $m_{u,v}$. \sq

Let $\tilde{d_k}$ be the coefficients of $R_1(m_{u,v})$.
It follows from \ref{containment} that
$\tilde{d_k}=0$ is necessary for $\m_0\subset \m_1$.
Similar to the
analysis of the $\tilde{c}^{i,j}_k$
in the previous section, we have
\[
\tilde{d_k}=\lambda^{u,v}_k-\lambda^h_k+
\mbox{polynomial in }\lambda^*_{\ell}\mbox{, }\ell<k.
\]
Combining these equations with the equations from
$R_1(Q^{i,j})=0$ (identical to those from
the last section), we solve them  for some $\lambda$ variables 
in terms of a certain
fixed collection of free $\lambda$\~variables.
This gives the required; cf. the end of the previous theorem.
%Therefore $\j^1_{\rk}(\vec{\De})$ is the graph of a regular
%function on
%an affine space and is hence isomorphic to an affine space.

\subsubsection{\sf Computing
\texorpdfstring{{\mathversion{bold}$\mu_{\De,g}$}}{mu(Delta)}}
To compute $\mu_{\De,g}$ we need to show that
the equations resulting from
the coefficients of $R_{\Delta_0}(T^{i,j})=0$ are redundant
once the equations from $R_{\Delta_1}(Q^{i,j})=0$ are
imposed.

\begin{definition}
Let $P$ be an element of a rank-$\rk$ $\r$-module
$\m\subset k[[z]](z^{1/\rk})$ such that $P$ has polynomial
coefficients in terms of the $\lambda$\~variables.
The \emph{coefficient ideal}
$\mathfrak{I}(P)$ of $P$ is the ideal generated by
the coefficients of $P$
in the polynomial ring over the $\lambda$\~variables. \sq
\end{definition}

The basic property of $\mathfrak{I}$ we will need
is given in the following lemma.

\begin{lemma}\label{ideals}
One has: $\mathfrak{I}(R_{\upsilon(\m)}(rP))\subset
\mathfrak{I}(R_{\upsilon(\m)}(P))$,
where $r\in \r$ and $P$ is as above.
\end{lemma}
{\it Proof.}
For any rank-$\rk $ module,
there exists $\cc\in \mathbb{N}$ such that for all $a>\cc$
we have $z^{a/\rk }\in \m$. Therefore a certain finite truncation of $rP$
is sufficient to determine $R_{\upsilon(\m)}(rP)$. Hence
we may assume that
$r$ is a polynomial in $\r$. Furthermore, since
$R_{\upsilon(\m)}$ is $\mathbb{C}$-linear
we may assume that $r=z^{k}$ for some $k\in \Gamma$.

The reduction algorithm for $rP$ is exactly parallel to that for
$P$ with valuations shifted by $\rk \cdot k$ until the
first
$i\not\in\upsilon(\m)$ such that $i+\rk \cdot k \in \upsilon(\m)$. 
The remaining coefficients
%will have a multiple of the final
%coefficient of $z^{i/\rk }$ in $R_{\upsilon(\m)}(P)$.
%Due to such $i$,
%the final coefficients
of $R_{\upsilon(\m)}(rP)$ will be sums of
multiples of the coefficients of $R_{\upsilon(\m)}(P)$.
This implies
that $\mathfrak{I}(R_{\upsilon(\m)}(rP))\subset
\mathfrak{I}(R_{\upsilon(\m)}(P))$.\sq
\smallskip

Continuing with the proof of the theorem,
note that $T^{i,j}=Q^{i,j}$ except when $a_{i,j}=a_{u,v}$ and
$a_{i,j+1}=a_{u,v}$. Recall that $u,v$ are integers such that
$g+\rk \cdot p=\rk \cdot a_{u,v}+u$.
Hence we only need to check that
\begin{align*}
\mathfrak{I}(R_0(T^{u,v}))&\subset \mathfrak{I}(R_1(Q^{u,v})), \\
\mathfrak{I}(R_0(T^{u,v-1}))&\subset \mathfrak{I}(R_1(Q^{u,v-1}))
\end{align*}
to show that the equations from $R_0(T^{u,v})=0$ and
$R_0(T^{u,v-1})=0$ are redundant.
Thus we need to relate $R_0$ and $R_1$.

\begin{lemma}\label{R01}
If $\m_0\subset \m_1$, then
$R_0(T^{u,v})\!=\!R_1(T^{u,v})$ and
$R_0(T^{u,v-1})\!=\!R_1(T^{u,v-1})$.
\end{lemma}
{\it Proof.}
Observe that
$
R_0(T^{u,v})-R_1(T^{u,v})=\sum_{i,j}'r_{i,j}m_{i,j}+sF,
$
where $r_{i,j},s\in \r$ and
by $'$ we indicate that the pair $i=u, j=v$ is omitted.

 Applying $R_1$ to each side we obtain
\[
R_0(T^{u,v})-R_1(T^{u,v})=R_1(sF);
\]
here the left-hand side is an eigenvector
for $R_1$ and all the $r_{i,j}m_{i,j}$ belongs to $\m_1$. 
By the previous lemma,
$\mathfrak{I}(R_1(sF))\subset \mathfrak{I}(R_1(F))$. By Lemma
\ref{containment}, $R_1(F)=0$. Hence $R_1(sF)=0$.
The proof for $T^{u,v-1}$ is identical. \sq
\smallskip

\subsubsection{\sf The end of the proof}
Now we will show that
$\mathfrak{I}(R_0(T^{u,v-1}))= \mathfrak{I}(R_1(Q^{u,v-1}))$.
Consider the difference
\begin{align*}
T^{u,v-1}-Q^{u,v-1}=-z^{(\alpha_{u,v}-\alpha_{u,v-1})p}F.
\end{align*}
Applying $R_1$ to both sides and using Lemmas
\ref{containment},\ref{ideals},\ref{R01}, we obtain
$R_0(T^{u,v-1})=R_1(Q^{u,v-1})$, which implies
$\mathfrak{I}(R_0(T^{u,v-1}))= \mathfrak{I}(R_1(Q^{u,v-1}))$.

Let us check that
$\mathfrak{I}(R_0(T^{u,v}))\subset \mathfrak{I}(R_1(Q^{u,v}))$.
Now we consider:
\[
T^{u,v}-z^pQ^{u,v}=z^qF.
\]
Applying $R_1$ to both sides and Lemmas
\ref{containment}, \ref{ideals} and \ref{R01}, we arrive at
$R_0(T^{u,v})=R_1(z^pQ^{u,v})$. However
$\mathfrak{I}(R_1(z^pQ^{u,v}))\subset
\mathfrak{I}(R_1(Q^{u,v}))$ due to
Lemma \ref{ideals}. The case for $j=p-1$ is similar.
Finally, we may calculate:
\begin{align*}
&\mu_{\De,g}=\gamma_{\Delta_1}(g)-\gamma_{\Delta_1}
(g+\rk \cdot p)-
\gamma_{\Delta_1}(\rk \cdot a_{u,v-1}+u+
\rk \cdot q)\\
-\gamma&_{\Delta_1}(g+\rk \cdot q)
+\gamma_{\Delta_1}(\rk \cdot a_{u,v-1}\!+\!u\!+\!\rk \cdot q)
+\gamma_{\Delta_1}(g\!+\!\rk \cdot p\!+\!\rk \cdot q) \\
=\gamma&_{\Delta_1}(g)-\gamma_{\Delta_1}(g\!+\!\rk \cdot p)-
(\gamma_{\Delta_1}(g\!+\!\rk \cdot q)-\gamma_{\Delta_1}
(g\!+\!\rk \cdot p\!+\!\rk \cdot q)).
\end{align*}
%\sq

%\setcounter{section}{-1}
\setcounter{equation}{0}
\section{\sc Examples for
\texorpdfstring{{\mathversion{bold}$\Ga\!=\!
\lan 4,6,6\!+\!v\ran$}}
{<4,6,6+v>}}
This section is devoted to numerical confirmations
of our Main Conjecture in the case of the (simplest)
series of non-torus plane curve singularities with the rings
$\r=\C[[z^4,z^6+z^{v}]]$ for
odd $v>6$. Then  $\de\!=\!6\! +\! \frac{v - 1}{2}\! -\! 1$
and $\Ga=\lan 4,6,6+v\ran$ for such $\r$.

\subsection {\bf Numerical simulations}

\subsubsection{\sf General remarks}
We will use the realization of $\o^{\rk}$
as the ring extension $\tilde{\o}=\o[z^{1/\rk}]$; see
(\ref{ramz}). Namely, we set:
\begin{align}\label{ramzz}
\ep_i=z^{i/\rk},\ \upsilon_i=\upsilon(\ep_i)=i/\rk,\
0\le i\le \rk-1.
\end{align}
Recall from (\ref{de-comp}) that for an arbitrary
$\Ga$\~module $\De\in \tilde{o}$ or for any
$\De=\De(\m)=\upsilon(\m)$,
\begin{align}
&\De^{(i)}\!=\!\{a\in \Z_+\mid a\!+\!\upsilon_i\!\in\!
\De\},\, dev^{(i)}\!=\!dev(\De^{(i)}),\,
dev(\De)\!=\!\sum_{i=0}^{\rk-1}dev^{(i)}.
\notag
\end{align}

Thus standard $\De$\~modules in $\cup_{i=1}^{\rk-1}
\{\upsilon_i+\Z\}$
are determined by  collections
$\{\De^{(i)}\}$ of usual standard $\Ga$\~modules in $\Z_+$ and
their number is $\kap^{\rk}$, where $\kap$ is the number
of usual standard $\Ga$\~modules.
\vskip 0.2cm

{\sf $D$\~sets and $D^\dag$\~sets.}
We define the $D$\~sets
$D(M)$ for standard rank-one $M$, {\em the sets of added gaps,\,}
 as $\De(M)\setminus \Ga$;
then $dev(\m)=|D(M)|$.

Accordingly, $D^{(i)}=\De^{(i)}\setminus \Ga$,
but we mostly use in this section
\comment{
$$
D(\m)\equal\rk\cdot\upsilon(\m)
\setminus \bigl(\Ga+ \{0,1,2,\ldots,\rk-1\}\bigr)
=\cup_{i=0}^{\rk\!-\!1}(i+D^{(i)}).
$$
 Notice multiplication by
$\rk$, which is convenient to avoid denominators,
and that we will always removed $\{0,\cdots,\rk\!-\!1\}$.
Accordingly,
}
\begin{align}\label{devDe}
D(\m)\,\equal&\,\,\rk\cdot\Bigl(\upsilon(\m)\setminus
\bigl(\cup_{i=0}^{\rk-1}\{\frac{i}{\rk}+\Ga\}\bigr)\Bigr)\,=\,
\cup_{i=0}^{\rk-1}\, \{i+\rk\cdot D^{(i)}\},\notag\\
&dev(\De(M))=\sum_{i=0}^{\rk-1} |D^{(i)}|=|D(\m)|,\
dev(\m)=|D(\m)|.
\end{align}
 Notice multiplication by
$\rk$, which is convenient to avoid denominators.
The deviation is zero for {\em minimal\,} standard modules
(minimal are inverse for $\rk=1$), and it
is $\rk\cdot\de$ for the greatest standard $\m=\tilde{\o}$.

For standard flags, we will set
$\d=\d(\vect{\m})=\{D_j,0\le i\le \ell\}$, where
$D_j=D(\m_j)$.
Accordingly, $\j_{\r,\rk}^\ell[\d]$ (notice using $[\cdot]$)
is a subvariety of flags of modules with a given $D$\~flag
$\d=\d(\vect{\m})$.

The elements of $D$ from (\ref{devDe})
are called {\em primitive\,}
if they cannot be represented in the form $\ga+g$
for $0\neq\ga\in \rk\cdot\Ga$ and $g\in D\cap\{0,\ldots,\rk-1\}$.
The notation
will be $D^\dag$; $D$ can be obviously uniquely recovered
from $D^\dag$.
We set $\d^\dag=\{D_j^\dag\}$ for flags.
\vfil

{\sf The tables.}
We discuss the rings $\r=\C[[z^4,z^6+z^{v}]]$
only in the case $\rk=2$; the total number
of $\Ga$\~modules is then from $25^2=625$ for $v=7$ to
 $41^2=1681$ for $v=15$ (when $\Ga\!=\! \lan 4,6,21\ran$),
the greatest example we reached. Recall that
{\em we multiply by $\rk=2$ the valuation $\upsilon$
in this section\,}
($\upsilon$ extends $\nu$ in $\r$) to avoid the denominators
in the outputs, i.e. the ``natural"
$D(\m)$ and $D(\m)^\dag$ defined in terms of initial
$\upsilon$ (with $\upsilon_i=i/\rk$)
are those below divided by $\rk=2$.
\vfil

The number of cells grows significantly
for flags, so we mostly analyzed numerically only
$0,1$\~flags ($\ell=0,1$); the $a$\~degree is $6$ for the whole
series. Some partial confirmations
of our conjecture were obtained beyond $\ell=0,1$
for sufficiently small or large powers of $q$.
We note that the corresponding
equations of the cells can have thousands of monomials
in the examples considered below even without flags,
and the monomials themselves can be long due to the large number
of $\la$\~variables.

If $v\le 15$ (i.e. till $\Ga\!\!=\!\lan 4,6,21\ran$),
then our conjecture is checked numerically for $\ell=0,1$
(i.e. for $a^0,a^1$) and also when $t=1$ for {\em any\,}
powers of $a$.
The DAHA superpolynomials can be calculated smoothly
for these and even more involved
examples, but they are long; for instance, there
are $24649$\, $q,t,a$\~monomials in the DAHA superpolynomial
for $\rk=2,v=15$.
\vfil

Due to such a size,  we provide only some portions of our
calculations in this paper, mostly focusing on an important
phenomenon of {\em non-affine cells\,}, which is of
independent
interest. Such cells always occur for non-torus singularities
in $\rk>1$ (in all examples we reached); they influence
the formulas of geometric superpolynomials in quite nontrivial
ways. The validity of our conjecture in the presence
of non-affine Piontkowski-type cells is a very good test;
needless to say that our conjecture
holds in all example we calculated.
\vskip 0.2cm

For the following, the maximal dimension
of the cells is reached for $D=\varnothing$; it
equals $\partial\equal 4\de$ for $\rk=2$. %$\varrho=2$.
Also, $2\de$
is the maximal deviation $|D|$ of
standard modules  (reached for $\m=\tilde{\o}$);
we will set
\begin{align}\label{defdev}
\kap\equal|D_\ell|=|D_0|+\ell \hbox{\, for a\, }
D\hbox{-flag\,\,  }
\d=\{D_0,\ldots,D_\ell\}.
\end{align}
%\smallskip

\subsubsection{\sf Our approach}
We begin our analysis of the cells $\j_{\r,\rk}^\ell[\d]$
with the straightforward elimination of the
$\la$\~variables, essentially due to Piontkowski.
Namely our computer program produces equations
for  (generalized) $\la$\~variables and then
eliminates all variables, occurring linearly in
at least one of the equations with a constant coefficient.
We express then such $\la$\~variables in terms of the
remaining ones in these equations and eliminate them in the
remaining equations.
We call such elimination/reduction {\em straightforward}.

\vfil
In this process, the program does not look for
possible linear changes of variables, i.e. for
linear combinations
of $\la$\~variables that can be eliminated this way
(by a straightforward substitution); finding them
is generally quite a challenge numerically and even
theoretically. If such a procedure eliminates all variables
involved, then the corresponding cell is affine\
(biregular isomorphic to some $\mathbb{A}^N$). If not,
i.e. when no further reduction of $\la$\~variables (or even
those in terms of their linear combinations) is possible,
the cell can be still affine; we then continue using more
advanced algebraic tools.
 Statistically, the instances
where the straightforward reduction of $\la$\~variables does not
eliminate (solve) all equations are actually relatively rare.
For singularities of torus type, such a reduction is
always sufficient (and all cells are affine).
\vfil

In the cases when straightforward substitutions
are insufficient,  called below
{\em potentially non-affine\,}, the program
switches to exact algebraic analysis of the remaining equations.
Fortunately the number of remaining equations is generally small
and they are not very involved in the examples we considered.
The (biregular) {\em types\,} of such cells are exactly
determined, which is combined (to double-check) with the count
of points modulo $p=3$. This prime number is a place of
{\em good reduction\,}
for the whole family $\r=\C[[z^4,z^6+z^{v}]]$. We provide all
potentially non-affine cells for $v=7,9,15$  If the resulting
cell is affine, then the symbol of the corresponding type is
$A$ in the tables below.
\vskip 0.2cm

The number of the variables remaining after the
straightforward elimination of $\la$\~variables minus
the number
of remaining equations will be called
{\em potential dimension\,} and denoted by $\partial -d$;
recall that $\partial\equal 4\de$ for $\rk=2$.
The {\em potential codimension\,} will be then denoted by
$d$; this is mostly important practically than theoretically,
but certainly reflects well the complexity of the cell.
This $d$ coincides with the actual codimension of resulting cells
for types $A,X,Y,Z,W$ in Table \ref{Table4-6-types} ($A$
stands for {\em affine\,} cells).

Providing all {\em potentially\,}
non-affine cells 
is quite reasonable
for us; this includes all non-affine cells
and affine ones where the affineness is
not straightforward,
which can be helpful in (future)  general theory.
There is extensive literature devoted to algebraic
varieties which are affine because of non-trivial reasons;
we face similar problems.
%\vskip 0.2cm\vfil

The table below 
for $\Ga\!=\!\lan 4,6,15\ran$ and $\ell=0,1$
provides the pairs of
{\em primitive  presentations\,} $\d^\dag$ (with primitive
elements only) for $0,1$\~flags
$\d^\dag=\{D_0^\dag, D_1^\dag\}$ in {\em all\,} cases
where the straightforward elimination of $\la$-variables
is insufficient. When $\ell=0$ (the
case of single modules), i.e. we have a single $D^\dag$,
we put $\{D^\dag,D^\dag\}$.
The next table gives {\em all\,} such primitive $D^\dag$
for $\ell\!=\!0$ (i.e. without flags)
for $\Ga\!=\!\lan 4,6,21\ran$. We omit the intermediate
cases $6+v=17,19$. There are no new
types there vs. those for $\Ga\!=\!\lan 4,6,17\ran$.
Also, only type $N$ there is missing for  
$\Ga\!=\!\lan 4,6,15\ran$. {\em We conjecture that
Table \ref{Table4-6-types} contains all possible
types of cells (for any flags) that can occur for
the whole series $\lan 4,6,6+v\ran$}.
%\smallskip

\subsubsection{\sf Field with 1 element}
This is the case of $t\!=\!1$. We note that if
$t\!=\!1$ and $a\!=-t^{n+1}\!=\!-1$, then 
the DAHA-Jones polynomials and superpolynomials become
trivial. However for $\rk\!=\!1$,
one can switch here to $q\!=1\!$ due to the super-duality.
Importantly, the substitution $t\mapsto q$ leads to 
{\em another\,}  
case of ``field with one element"; see \cite{ChRH}.

The following general observation 
can be proved for the series $\r=\C[[z^4,z^6+z^{v}]]$.
There exists no greater than one equation with
nonzero constant term where no single $\la$\~variable
occurs only linearly. The constant term is actually
$\pm 2$ (if nonzero) in the normalization when the
coefficients of all monomials are relatively prime
integers (all of them); $p=2$ is the place of bad reduction
for this series. This claim formally
includes non-admissible flags; then this equation
is (always) of the form $0=2$.
The presence of such
an equation for any module in a flag (in this family of
rings) gives that there will be no other such equations of
this type for the whole flag.

Making now $t=1$, we obtain that only the flags
with such equations contribute $0$ to the
superpolynomial evaluated at $t=1$ (as well as non-admissible
ones do for any $t$). The corresponding (symbols of)
types from Table \ref{Table4-6-types}
are $X,Z,L$.
Let us call them 
{\em $C^*$\~types}, adding the non-admissible cells to this
category too. In the
cases $A,Y,W,M,N$, called below
{\em $\C$\~types\,}, the contribution
to the geometric superpolynomial at $t=1$ is always
$1$ ($a^m q^{\kap}$, to be exact). Thus the contribution
at $t=1$ is expected to be always $0$ (for types $X,Y,Z$) or $1$
otherwise
due to our conjecture on the types for the whole family
$\r=\C[[z^4,z^6+z^{v}]]$.

\smallskip
This gives that
the superpolynomial (any powers of $a$) at $t=1$ can
be found if one knows the types $\C,\C^*$ for
{\em single\,} $\d=D_0$ (no flags). Namely, a flag $\d$ is of type
$\C^*$ (and therefore contributes $0$) if and only if at least one
$D_i$ is of such type; otherwise it is of type $\C$
(and then contributes $1$). We have
a draft of justification of this claim for
$\r=\C[[z^4,z^6+z^{v}]]$. This  
seems more general than just for this family.

Also, for any $\r$, our Main Conjecture gives that
if at least one rank-one component  $D_i^{(j)}$ is
non-admissible (the whole flag can be still admissible),
then its contribution to the geometric superpolynomial
at $t=1$ is zero, i.e. it is a non-affine cell
of $\C^*$\~type. Indeed, a DAHA superpolynomial at $t=1$
for $\rk>1$
is that of rank one raised to
the $\rk${\tiny th} power, an important symmetry of DAHA
superpolynomials.
\vskip 0.2cm

Finally,  the passage from single
modules to flags is essentially combinatorial
at $t=1$. One only needs to know the list of  $D$\~flags,
which is an interesting combinatorial problem in its own
right, and the corresponding $\C^*$\~types for
{\em single\,} $D$ (which can be further reduced
to the case $\rk=1$). For instance, the knowledge of
non-admissible modules is not needed here ($\C^*$\~types
include them); let us add some details here.

%\vskip 0.2cm
%We note that we substitute $t\mapsto 1$ in the DAHA
%superpolynomials,
%{\em not in the DAHA itself} (which would result in trivial
%superpolynomials).

\vskip 0.2cm

{\sf Unusual non-admissibility.}
It is important to note that {\em it is not
true for $\rk>1$
that a $D$\~flag $\d=\{D_i\}$ is non-admissible, i.e.
does not correspond to any flag $\vect{\m}$,
(if and) only if one of $D_i$ is non-admissible\,}, which
always holds in the examples we reached for $\rk=1$. The
$\C^*$\~types (non-admissible $\d$ are included in this
class)
have this ``only if" property, at least for the families
we calculated, but this is not always true for
non-admissibility. Let us provide examples.

There are (only) $2$ such nonadmissible
$1$\~flags for $\r\!=\!\C[[z^4,z^6+z^9]]$
with admissible components (both). Their (full) $\{D_0,D_1\}$
are:
\begin{align*}
&\{[5, 18, 22, 26, 34, 35], \ [5, 14, 18, 22, 26, 34, 35] \},\\
&\{[4, 5, 18, 22, 26, 34, 35],\
[4, 5, 14, 18, 22, 26, 34, 35]\}.
\end{align*}
Recall that we {\em always multiply $D(\m)$ by $2$ in any
outputs\,} in this section.
The corresponding {\em primitive\,} $D_i^\dag$ are:
\begin{align*}
&\{\{5, 18, 22\}, \{5, 14, 18\} \},\
\{\{4, 5, 18, 22\},\{4, 5, 14, 18\}\}.
\end{align*}
The cell-types for all four
$D$\~modules above are $X$. There are $6,10$
nonadmissible $1$\~flags with admissible components
correspondingly for $\r=\C[[z^4,z^6+z^{11,13}]]$ and no such
pairs for $\r=\C[[z^4,z^6+z^{7}]].$ The types for the components
$D_i (i=0,1)$ of such non-admissible
pairs are always $X$ in the examples we
calculated. Generally, type $\C^*$ is expected here due
to the discussion above, not only $X$.
\medskip

\subsection {\bf Non-affine cells for
\texorpdfstring{{\mathversion{bold}$v\!=\!9,15$}}
{v=9,15}}
\subsubsection {\sf The case of
\texorpdfstring{{\mathversion{bold}$v\!=\!9$}}
{v=9}}
The first table gives all {\em potentially non-affine
cells\,} for $\r=\C[[z^4,z^6+z^{v}]]$ and
$\Ga=\lan 4,6,15\ran$, where one table entry contains
1) $\d^\dag=\{D_0^\dag, D_1^\dag\}$ (they are
coinciding for $\ell=0$),
2) the corresponding {\em potential dimensions} of
$\j_{\r,\rk}^{\ell=1}[\d]$ and 3) the types.
One needs the table of types (we provide) to calculate
the actual dimension. Recall that the potential
dimension is $\partial -d$, which is the number of variables
minus the number of equations upon the
{\em straightforward elimination\,}
of $\la$\~variables.

In all tables
below we sort the entries with respect
to the $\kap$-deviations
$\kap$ from $0$ to $2\de$; the order is from left to right
and then down. Recall that
$\kap$ is $|D_1|$; it cannot be seen from $|D_1^\dag|$\,,
but can be readily calculated in terms of $D_1^\dag$.

\makeatletter
\newcommand{\srcsize}{\@setfontsize{\srcsize}{7pt}{7pt}}
\makeatother

{\srcsize
\hskip -0.7cm
\begin{longtable}{|l|l|}
\caption*
{\small
Types of cells for $\r=\C[[z^4,z^6+z^{9}]]$ when $\ell=0,1$:}\\
\hline
\{18, 22, 35\},\{4, 18, 22, 35\}:\ 34,A \ $\rightsquigarrow$&
\{18, 22, 35\},\{5, 18, 22\}:\ 33,X\ $\rightsquigarrow$\\
\{19, 23, 34\},\{5, 19, 23, 34\}:\ 32,A &
\{4, 18, 22, 35\},\{4, 18, 22, 35\}:\ 32,A\\
\{4, 18, 22, 35\},\{4, 5, 18, 22\}:\ 34,X &
\{4, 18, 22, 35\},\{4, 14, 18, 35\}:\ 33,A\\
\{4, 18, 22, 35\},\{4, 18, 22, 23\}:\ 33,A &
\{4, 18, 22, 35\},\{4, 18, 22, 27\}:\ 32,A\\
\{5, 18, 22\},\{5, 18, 22\}:\ 31,X &
\{5, 18, 22\},\{4, 5, 18, 22\}:\ 33,X\\
\{5, 18, 22\},\{5, 18, 22, 23\}:\ 32,X &
\{5, 18, 22\},\{5, 18, 22, 27\}:\ 31,X\\
\{5, 19, 23, 34\},\{5, 19, 23, 34\}:\ 30,A &
\{5, 19, 23, 34\},\{5, 15, 19, 34\}:\ 32,A\\
\{5, 19, 23, 34\},\{5, 19, 22, 23\}:\ 31,A &
\{5, 19, 23, 34\},\{5, 19, 23, 26\}:\ 30,A\\
\{14, 18, 35\},\{4, 14, 18, 35\}:\ 33,A &
\{14, 18, 35\},\{5, 14, 18\}:\ 32,X\\
\{15, 19, 34\},\{4, 15, 19\}:\ 32,X &
\{15, 19, 34\},\{5, 15, 19, 34\}:\ 31,Y\\
\{18, 22, 27\},\{4, 18, 22, 27\}:\ 31,A &
\{18, 22, 27\},\{5, 18, 22, 27\}:\ 30,X\\
\{4, 5, 18, 22\},\{4, 5, 18, 22\}:\ 31,X &
\{4, 5, 18, 22\},\{4, 5, 18, 22, 23\}:\ 32,X\\
\{4, 5, 18, 22\},\{4, 5, 18, 22, 27\}:\ 31,X &
\{4, 14, 18, 35\},\{4, 14, 18, 35\}:\ 31,A\\
\{4, 14, 18, 35\},\{4, 5, 14, 18\}:\ 33,X &
\{4, 14, 18, 35\},\{4, 14, 18, 23\}:\ 32,A\\
\{4, 14, 18, 35\},\{4, 14, 18, 27\}:\ 31,A &
\{4, 15, 19\},\{4, 15, 19\}:\ 30,X\\
\{4, 15, 19\},\{4, 5, 15, 19\}:\ 32,X &
\{4, 15, 19\},\{4, 15, 19, 22\}:\ 31,X\\
\{4, 15, 19\},\{4, 15, 19, 26\}:\ 30,X &
\{4, 18, 22, 27\},\{4, 18, 22, 27\}:\ 29,A\\
\{4, 18, 22, 27\},\{4, 5, 18, 22, 27\}:\ 32,X &
\{4, 18, 22, 27\},\{4, 14, 18, 27\}:\ 31,A\\
\{4, 18, 22, 27\},\{4, 18, 19, 22\}:\ 30,A &
\{4, 18, 22, 27\},\{4, 18, 22, 23, 27\}:\ 29,A\\
\{5, 14, 18\},\{5, 14, 18\}:\ 30,X &
\{5, 14, 18\},\{4, 5, 14, 18\}:\ 32,X\\
\{5, 14, 18\},\{5, 14, 18, 23\}:\ 31,X &
\{5, 14, 18\},\{5, 14, 18, 27\}:\ 30,X\\
\{5, 15, 19, 34\},\{5, 15, 19, 34\}:\ 29,Y &
\{5, 15, 19, 34\},\{4, 5, 15, 19\}:\ 31,X\\
\{5, 15, 19, 34\},\{5, 15, 19, 22\}:\ 30,Y &
\{5, 15, 19, 34\},\{5, 15, 19, 26\}:\ 29,Y\\
\{5, 18, 22, 27\},\{5, 18, 22, 27\}:\ 28,X &
\{5, 18, 22, 27\},\{4, 5, 18, 22, 27\}:\ 31,X\\
\{5, 18, 22, 27\},\{5, 14, 18, 27\}:\ 30,X &
\{5, 18, 22, 27\},\{5, 18, 19, 22\}:\ 29,X\\
\{5, 18, 22, 27\},\{5, 18, 22, 23, 27\}:\ 28,X &
\{14, 18, 27\},\{4, 14, 18, 27\}:\ 30,Y\\
\{14, 18, 27\},\{5, 14, 18, 27\}:\ 29,Z &
\{15, 19, 26\},\{4, 15, 19, 26\}:\ 29,X\\
\{15, 19, 26\},\{5, 15, 19, 26\}:\ 28,Y &
\{4, 5, 14, 18\},\{4, 5, 14, 18\}:\ 30,X\\
\{4, 5, 14, 18\},\{4, 5, 14, 18, 23\}:\ 31,X &
\{4, 5, 14, 18\},\{4, 5, 14, 18, 27\}:\ 30,X\\
\{4, 5, 15, 19\},\{4, 5, 15, 19\}:\ 29,X &
\{4, 5, 15, 19\},\{4, 5, 15, 19, 22\}:\ 30,X\\
\{4, 5, 15, 19\},\{4, 5, 15, 19, 26\}:\ 29,X &
\{4, 5, 18, 22, 27\},\{4, 5, 18, 22, 27\}:\ 29,X\\
\{4, 5, 18, 22, 27\},\{4, 5, 14, 18, 27\}:\ 31,X &
\{4, 5, 18, 22, 27\},\{4, 5, 18, 19, 22\}:\ 30,X\\
\{4, 5, 18, 22, 27\},\{4, 5, 18, 22, 23, 27\}:\ 29,X &
\{4, 14, 18, 27\},\{4, 14, 18, 27\}:\ 28,Y\\
\{4, 14, 18, 27\},\{4, 5, 14, 18, 27\}:\ 30,L &
\{4, 14, 18, 27\},\{4, 10, 14, 27\}:\ 30,A\\
\{4, 14, 18, 27\},\{4, 14, 18, 19\}:\ 29,Y &
\{4, 14, 18, 27\},\{4, 14, 18, 3, 27\}:\ 28,Y\\
\{4, 15, 19, 26\},\{4, 15, 19, 26\}:\ 27,X &
\{4, 15, 19, 26\},\{4, 5, 15, 19, 26\}:\ 30,X\\
\{4, 15, 19, 26\},\{4, 11, 15, 26\}:\ 29,X &
\{4, 15, 19, 26\},\{4, 15, 18, 19\}:\ 28,X\\
\{4, 15, 19, 26\},\{4, 15, 19, 22, 26\}:\ 27,X &
\{5, 14, 18, 27\},\{5, 14, 18, 27\}:\ 27,Z\\
\{5, 14, 18, 27\},\{4, 5, 14, 18, 27\}:\ 29,L &
\{5, 14, 18, 27\},\{5, 10, 14, 27\}:\ 29,X\\
\{5, 14, 18, 27\},\{5, 14, 18, 19\}:\ 28,Z &
\{5, 14, 18, 27\},\{5, 14, 18, 23, 27\}:\ 27,Z\\
\{5, 15, 19, 26\},\{5, 15, 19, 26\}:\ 26,Y &
\{5, 15, 19, 26\},\{4, 5, 15, 19, 26\}:\ 29,X\\
\{5, 15, 19, 26\},\{5, 11, 15, 26\}:\ 28,Y &
\{5, 15, 19, 26\},\{5, 15, 18, 19\}:\ 27,Y\\
\{5, 15, 19, 26\},\{5, 15, 19, 22, 26\}:\ 26,Y &
\{10, 14, 27\},\{5, 10, 14, 27\}:\ 28,X\\
\{11, 15, 26\},\{4, 11, 15, 26\}:\ 28,X &
\{11, 15, 26\},\{5, 11, 15, 26\}:\ 27,Y\\
\{14, 18, 23, 27\},\{4, 14, 18, 23, 27\}:\ 27,Y &
\{14, 18, 23, 27\},\{5, 14, 18, 23, 27\}:\ 26,Y\\
\{4, 5, 14, 18, 27\},\{4, 5, 14, 18, 27\}:\ 27,L &
\{4, 5, 14, 18, 27\},\{4, 5, 10, 14, 27\}:\ 30,X\\
\{4, 5, 14, 18, 27\},\{4, 5, 14, 18, 19\}:\ 28,L &
\{4, 5, 14, 18, 27\},\{4, 5, 14, 18, 23, 27\}:\ 27,L\\
\{4, 5, 15, 19, 26\},\{4, 5, 15, 19, 26\}:\ 27,X &
\{4, 5, 15, 19, 26\},\{4, 5, 11, 15, 26\}:\ 29,X\\
\{4, 5, 15, 19, 26\},\{4, 5, 15, 18, 19\}:\ 28,X &
\{4, 5, 15, 19, 26\},\{4, 5, 15, 19, 22, 26\}:\ 27,X\\
\{4, 10, 14, 27\},\{4, 5, 10, 14, 27\}:\ 29,X &
\{4, 11, 15, 26\},\{4, 11, 15, 26\}:\ 26,X\\
\{4, 11, 15, 26\},\{4, 5, 11, 15, 26\}:\ 28,X &
\{4, 11, 15, 26\},\{4, 11, 15, 18\}:\ 27,X\\
\{4, 11, 15, 26\},\{4, 11, 15, 22, 26\}:\ 26,X &
\{4, 14, 18, 23, 27\},\{4, 14, 18, 23, 27\}:\ 25,Y\\
\{4, 14, 18, 23, 27\},\{4, 5, 14, 18, 23, 27\}:\ 28,M &
\{4, 14, 18, 23, 27\},\{4, 10, 14, 23, 27\}:\ 27,Y\\
\{4, 14, 18, 23, 27\},\{4, 14, 15, 18\}:\ 26,Y &
\{4, 14, 18, 23, 27\},\{4, 14, 18, 19, 23\}:\ 25,Y\\
\{5, 10, 14, 27\},\{5, 10, 14, 27\}:\ 26,X &
\{5, 10, 14, 27\},\{4, 5, 10, 14, 27\}:\ 28,X\\
\{5, 10, 14, 27\},\{5, 10, 14, 19\}:\ 27,X &
\{5, 10, 14, 27\},\{5, 10, 14, 23, 27\}:\ 26,X\\
\{5, 11, 15, 26\},\{5, 11, 15, 26\}:\ 25,Y &
\{5, 11, 15, 26\},\{4, 5, 11, 15, 26\}:\ 27,X\\
\{5, 11, 15, 26\},\{5, 11, 15, 18\}:\ 26,Y &
\{5, 11, 15, 26\},\{5, 11, 15, 22, 26\}:\ 25,Y\\
\{5, 14, 18, 23, 27\},\{5, 14, 18, 23, 27\}:\ 24,Y &
\{5, 14, 18, 23, 27\},\{4, 5, 14, 18, 23, 27\}:\ 27,M\\
\{5, 14, 18, 23, 27\},\{5, 10, 14, 23, 27\}:\ 26,Y &
\{5, 14, 18, 23, 27\},\{5, 14, 15, 18\}:\ 25,Y\\
\{5, 14, 18, 23, 27\},\{5, 14, 18, 19, 23\}:\ 24,Y &
\{10, 14, 23, 27\},\{4, 10, 14, 23, 27\}:\ 26,Y\\
\{10, 14, 23, 27\},\{5, 10, 14, 23, 27\}:\ 25,Y &
\{4, 5, 10, 14, 27\},\{4, 5, 10, 14, 27\}:\ 26,X\\
\{4, 5, 10, 14, 27\},\{4, 5, 10, 14, 19\}:\ 27,X &
\{4, 5, 10, 14, 27\},\{4, 5, 10, 14, 23, 27\}:\ 26,X\\
\{4, 5, 11, 15, 26\},\{4, 5, 11, 15, 26\}:\ 25,X &
\{4, 5, 11, 15, 26\},\{4, 5, 11, 15, 18\}:\ 26,X\\
\{4, 5, 11, 15, 26\},\{4, 5, 11, 15, 22, 26\}:\ 25,X &
\{4, 5, 14, 18, 23, 27\},\{4, 5, 14, 18, 23, 27\}:\ 25,M\\
\{4, 5, 14, 18, 23, 27\},\{4, 5, 10, 14, 23, 27\}:\ 27,M &
\{4, 5, 14, 18, 23, 27\},\{4, 5, 14, 15, 18\}:\ 26,M\\
\{4, 5, 14, 18, 23, 27\},\{4, 5, 14, 18, 19, 23\}:\ 25,M &
\{4, 10, 14, 23, 27\},\{4, 10, 14, 23, 27\}:\ 24,Y\\
\{4, 10, 14, 23, 27\},\{4, 5, 10, 14, 23, 27\}:\ 26,M &
\{4, 10, 14, 23, 27\},\{4, 10, 14, 15\}:\ 25,Y\\
\{4, 10, 14, 23, 27\},\{4, 10, 14, 19, 23\}:\ 24,Y &
\{5, 10, 14, 23, 27\},\{5, 10, 14, 23, 27\}:\ 23,Y\\
\{5, 10, 14, 23, 27\},\{4, 5, 10, 14, 23, 27\}:\ 25,M &
\{5, 10, 14, 23, 27\},\{5, 10, 14, 15\}:\ 24,Y\\
\{5, 10, 14, 23, 27\},\{5, 10, 14, 19, 23\}:\ 23,Y &
\{4, 5, 10, 14, 23, 27\},\{4, 5, 10, 14, 23, 27\}:\ 23,M\\
\{4, 5, 10, 14, 23, 27\},\{4, 5, 10, 14, 15\}:\ 24,M &
\{4, 5, 10, 14, 23, 27\},\{4, 5, 10, 14, 19, 23\}:\ 23,M\\
\hline
\caption*{\small Potentially non-affine
$\{D_0^\dag, D_1^\dag\}$ for $\rk=2,
\Ga=\lan 4,6,15\ran$}
\end{longtable}
\hspace {0.1cm}
}

\vskip -1.cm
The types are explained below. Recall that the number before
the type is $\partial -d$, which is the number of variables
minus the number of equations upon the
{\em straightforward elimination\,}
of $\la$\~variables. It coincides with the dimension of a cell
unless for types $L,M,N$. In these cases
the actual {\em codimension\,} is not $d$; it is $d-1$
for $L,M$ and $d-2$ for $N$, though type $N$ does not
appear for $\Ga=\lan 4,6,15\ran$.
We note that only types $A,X$
occur for $\Ga=\lan 4,6,13\ran$ for any flags. Type $N$ appears
for the first time for $\Ga=\lan 4,6,17\ran$.
No new types occur for $\lan 4,6,19\ran$ and $\lan 4,6,21\ran$
at least as $\ell=0,1$. See the conjecture on types below.
%\vfil

We omit the tables for
$17$ and $19$ and provide dimensions only
for $13$ ($a=0$). See the online Appendix to the paper
for the case of $15$ for $\ell=0,1$. For $\ell=0,1$
and partially for any $\ell$, our Main Conjecture was checked
for $13,15,17,19,21$ (quite a test of our conjecture and
the validity of our programs).

%\smallskip

\subsubsection {\sf The conjecture on types}
We conjecture that admissible flags for
the family $\C[t^4,t^6+t^v]$ can be only those from
Table \ref{Table4-6-types}. Concerning
the notations there, ``$+$" means
a set-theoretical disjoint union. Type $Z$ is
a union of two copies of affine spaces
$\mathbb{A}^{N}\setminus\mathbb{A}^{N-1}$ as an
algebraic variety (i.e. disconnected). The cells of types $L,N$
are (geometrically) connected; the corresponding intersections
are straightforward. Note that we need only the set-theoretical
presentation to determine the contributions of cells
to the geometric $\h$\~polynomial.
%\vfil

The spaces
$\mathbb{A}^{N}\vee\mathbb{A}^{N}$ and
$\mathbb{A}^{N+1}\vee\mathbb{A}^{N}$ are the
unions of these two affine spaces with the ``natural"
intersection $\mathbb{A}^{N-1}$. Set-theoretically
(or in the Grothendieck ring of varieties over $\C$)
these cells can be presented as
$\mathbb{A}^{N+\ep}+
(\mathbb{A}^{N}\setminus \mathbb{A}^{N-1})$ for $\ep=0,1$.
%\vfil

Type $W$ requires more comment;
geometrically it is the union of three copies of
$\mathbb{A}^N$ where the
intersection of the first and the third coincides
with the total (triple) intersection. This is reflected
in the contribution of $W$ to the superpolynomial
from the table (actually the class of $W$ in the Grothendieck
ring).
%Its $t$\~part generally counts $\mathbb{F}_p$\~points
%of the corresponding cell upon $t=1/p$ for the places
%$p$ of good reduction.
%\vfill

Recall that {\em potential dimension\,} is denoted by
$\partial -d$ and $d$ is
{\em potential codimension\,};
$\partial\equal 4\de$ for $\rk=2$. The contribution
of an {\em affine cell\,} (called type $A$) to the geometric
superpolynomial is simply $q^{\kap}t^d$ (then $d$ coincides with
the actual codimension), where $\kap$ is from
(\ref{defdev}).
\vskip 0.5cm
%\vfill

%\setcounter{table}{0}
{\small
\renewcommand\thetable{1}
\begin{longtable}{|rr|c|l|}
\hline
label&& type in terms of $N=\partial-d$ & contribution to $\h$\\
\hline
A $\rightsquigarrow$& &$\mathbb{A}^{N}$
& $q^{\kap}\,\, t^d $\\
X $\rightsquigarrow$& & $\mathbb{A}^{N}
\setminus \mathbb{A}^{N-1}$
& $q^{\kap}(t^d-t^{d+1})$\\
Y $\rightsquigarrow$& & $\mathbb{A}^{N}\vee\mathbb{A}^{N}$
& $q^{\kap}(2 t^d-t^{d+1})$\\
Z\, $\rightsquigarrow$& & $(\mathbb{A}^{N}
\setminus\mathbb{A}^{N-1})+
(\mathbb{A}^{N}\setminus\mathbb{A}^{N-1})$
& $q^{\kap}(2t^d-2t^{d+1})$\\
W\! $\rightsquigarrow$& & $\mathbb{A}^{N}\vee\mathbb{A}^{N}
\vee\mathbb{A}^{N}$
& $q^{\kap}(3t^d-2t^{d+1})$\\
L $\rightsquigarrow$& & $(\mathbb{A}^{N+1}
\setminus \mathbb{A}^N)+
(\mathbb{A}^{N}\setminus \mathbb{A}^{N-1})$
& $q^{\kap}(t^{d-1}-t^{d+1})$\\
M $\rightsquigarrow$& & $\mathbb{A}^{N+1}\vee \mathbb{A}^N$
& $q^{\kap}(t^{d-1}+t^d-t^{d+1})$\\
N $\rightsquigarrow$& & $\mathbb{A}^{N+2}\!+
\!(\mathbb{A}^{N+1}\!\setminus
\mathbb{A}^N)\!+\!
(\mathbb{A}^{N}\!\setminus \mathbb{A}^{N-1})$
& $q^{\kap}(t^{d-2}+t^{d-1}\!-t^{d+1})$\\
\hline
\caption{Non-affine types for $\rk=2,\ \Ga=\lan 4,6,6+v\ran$}
\end{longtable}
}
\hspace {0.1cm}
\vspace {0.cm}
\label{Table4-6-types}

%Types of cells for $\r=\C[z^4,z^6+z^{15}]\ran$ for $m=0$
%(without flags):

\subsubsection{\sf The case of
\texorpdfstring{{\mathversion{bold}$v\!=\!15$}}
{v=15}}
Counting the dimensions of $1$\~flags becomes
very long for $\r=\C[[z^4,z^6+z^{15}]]$, so we
calculated all of them only for $0$\~flags ($\ell=0$).
Potentially non-affine cells will be provided,
they are of interest. Also, one can see some
combinatorial patterns, when going through the
tables of types we provide in this paper for
the family $\C[[z^4,z^6+z^v]]$ as well as patterns
of non-admissible modules.

{\srcsize
\hskip -0.7cm
\begin{longtable}{|l|l|l|}
\caption*
{\small
Types of cells for $\r=\C[[z^4,z^6+z^{15}]]$ when $\ell=0$:}\\
\hline
\{ 4, 30, 34, 47\}:\ 44,A \ $\rightsquigarrow$&
\{ 5, 30, 34\}:\ 43,X \ $\rightsquigarrow$&
\{ 5, 31, 35, 46\}:\ 42,A \ $\rightsquigarrow$\\
\{ 4, 5, 30, 34\}:\ 43,X &
\{ 4, 26, 30, 47\}:\ 43,A &
\{ 4, 27, 31\}:\ 42,X\\
\{ 4, 30, 34, 39\}:\ 41,A &
\{ 5, 26, 30\}:\ 42,X &
\{ 5, 27, 31, 46\}:\ 41,Y\\
\{ 5, 30, 34, 39\}:\ 40,X &
\{ 4, 5, 26, 30\}:\ 42,X &
\{ 4, 5, 27, 31\}:\ 41,X\\
\{ 4, 5, 30, 34, 39\}:\ 41,X &
\{ 4, 23, 27\}:\ 41,X &
\{ 4, 26, 30, 39\}:\ 40,Y\\
\{ 4, 27, 31, 38\}:\ 39,X &
\{ 5, 22, 26\}:\ 41,X &
\{ 5, 23, 27, 46\}:\ 40,Y\\
\{ 5, 26, 30, 39\}:\ 39,Z &
\{ 5, 27, 31, 38\}:\ 38,Y &
\{ 4, 5, 22, 26\}:\ 41,X\\
\{ 4, 5, 23, 27\}:\ 40,X &
\{ 4, 5, 26, 30, 39\}:\ 39,L &
\{ 4, 5, 27, 31, 38\}:\ 39,X\\
\{ 4, 19, 23\}:\ 40,X &
\{ 4, 22, 26, 39\}:\ 39,Y &
\{ 4, 23, 27, 38\}:\ 38,Z\\
\{ 4, 26, 30, 35, 39\}:\ 37,Y &
\{ 5, 18, 22\}:\ 40,X &
\{ 5, 19, 23, 46\}:\ 39,Y\\
\{ 5, 22, 26, 39\}:\ 38,Z &
\{ 5, 23, 27, 38\}:\ 37,W &
\{ 5, 26, 30, 35, 39\}:\ 36,Y\\
\{ 4, 5, 18, 22\}:\ 40,X &
\{ 4, 5, 19, 23\}:\ 39,X &
\{ 4, 5, 22, 26, 39\}:\ 38,L\\
\{ 4, 5, 23, 27, 38\}:\ 37,L &
\{ 4, 5, 26, 30, 35, 39\}:\ 37,M &
\{ 4, 15, 19\}:\ 39,X\\
\{ 4, 18, 22, 39\}:\ 38,Y &
\{ 4, 19, 23, 38\}:\ 37,Z &
\{ 4, 22, 26, 35, 39\}:\ 36,W\\
\{ 4, 23, 27, 34, 38\}:\ 35,Y &
\{ 5, 14, 18\}:\ 39,X &
\{ 5, 15, 19, 46\}:\ 38,Y\\
\{ 5, 18, 22, 39\}:\ 37,Z &
\{ 5, 19, 23, 38\}:\ 36,W &
\{ 5, 22, 26, 35, 39\}:\ 35,W\\
\{ 5, 23, 27, 34, 38\}:\ 34,Y &
\{ 4, 5, 14, 18\}:\ 39,X &
\{ 4, 5, 15, 19\}:\ 38,X\\
\{ 4, 5, 18, 22, 39\}:\ 37,L &
\{ 4, 5, 19, 23, 38\}:\ 36,L &
\{ 4, 5, 22, 26, 35, 39\}:\ 35,N\\
\{ 4, 5, 23, 27, 34, 38\}:\ 35,M &
\{ 4, 14, 18, 39\}:\ 37,Y &
\{ 4, 15, 19, 38\}:\ 36,Z\\
\{ 4, 18, 22, 35, 39\}:\ 35,W &
\{ 4, 19, 23, 34, 38\}:\ 34,W &
\{ 4, 22, 26, 31, 35\}:\ 33,Y\\
\{ 5, 14, 18, 39\}:\ 36,Z &
\{ 5, 15, 19, 38\}:\ 35,W &
\{ 5, 18, 22, 35, 39\}:\ 34,W\\
\{ 5, 19, 23, 34, 38\}:\ 33,W &
\{ 5, 22, 26, 31, 35\}:\ 32,Y &
\{ 4, 5, 14, 18, 39\}:\ 36,L\\
\{ 4, 5, 15, 19, 38\}:\ 35,L &
\{ 4, 5, 18, 22, 35, 39\}:\ 34,N &
\{ 4, 5, 19, 23, 34, 38\}:\ 33,N\\
\{ 4, 5, 22, 26, 31, 35\}:\ 33,M &
\{ 4, 11, 15, 38\}:\ 35,X &
\{ 4, 14, 18, 35, 39\}:\ 34,W\\
\{ 4, 15, 19, 34, 38\}:\ 33,W &
\{ 4, 18, 22, 31, 35\}:\ 32,W &
\{ 4, 19, 23, 30, 34\}:\ 31,Y\\
\{ 5, 10, 14, 39\}:\ 35,X &
\{ 5, 11, 15, 38\}:\ 34,Y &
\{ 5, 14, 18, 35, 39\}:\ 33,W\\
\{ 5, 15, 19, 34, 38\}:\ 32,W &
\{ 5, 18, 22, 31, 35\}:\ 31,W &
\{ 5, 19, 23, 30, 34\}:\ 30,Y\\
\{ 4, 5, 10, 14, 39\}:\ 35,X &
\{ 4, 5, 11, 15, 38\}:\ 34,X &
\{ 4, 5, 14, 18, 35, 39\}:\ 33,N\\
\{ 4, 5, 15, 19, 34, 38\}:\ 32,N &
\{ 4, 5, 18, 22, 31, 35\}:\ 31,N &
\{ 4, 5, 19, 23, 30, 34\}:\ 31,M\\
\{ 4, 10, 14, 35, 39\}:\ 33,Y &
\{ 4, 11, 15, 34, 38\}:\ 32,Y &
\{ 4, 14, 18, 31, 35\}:\ 31,W\\
\{ 4, 15, 19, 30, 34\}:\ 30,W &
\{ 4, 18, 22, 27, 31\}:\ 29,Y &
\{ 5, 10, 14, 35, 39\}:\ 32,Y\\
\{ 5, 11, 15, 34, 38\}:\ 31,Y &
\{ 5, 14, 18, 31, 35\}:\ 30,W &
\{ 5, 15, 19, 30, 34\}:\ 29,W\\
\{ 5, 18, 22, 27, 31\}:\ 28,Y &
\{ 4, 5, 10, 14, 35, 39\}:\ 32,M &
\{ 4, 5, 11, 15, 34, 38\}:\ 31,M\\
\{ 4, 5, 14, 18, 31, 35\}:\ 30,N &
\{ 4, 5, 15, 19, 30, 34\}:\ 29,N &
\{ 4, 5, 18, 22, 27, 31\}:\ 29,M\\
\{ 4, 10, 14, 31, 35\}:\ 30,Y &
\{ 4, 11, 15, 30, 34\}:\ 29,Y &
\{ 4, 14, 18, 27, 31\}:\ 28,W\\
\{ 4, 15, 19, 26, 30\}:\ 27,Y &
\{ 5, 10, 14, 31, 35\}:\ 29,Y &
\{ 5, 11, 15, 30, 34\}:\ 28,Y\\
\{ 5, 14, 18, 27, 31\}:\ 27,W &
\{ 5, 15, 19, 26, 30\}:\ 26,Y &
\{ 4, 5, 10, 14, 31, 35\}:\ 29,M\\
\{ 4, 5, 11, 15, 30, 34\}:\ 28,M &
\{ 4, 5, 14, 18, 27, 31\}:\ 27,N &
\{ 4, 5, 15, 19, 26, 30\}:\ 27,M\\
\{ 4, 10, 14, 27, 31\}:\ 27,Y &
\{ 4, 11, 15, 26, 30\}:\ 26,Y &
\{ 4, 14, 18, 23, 27\}:\ 25,Y\\
\{ 5, 10, 14, 27, 31\}:\ 26,Y &
\{ 5, 11, 15, 26, 30\}:\ 25,Y &
\{ 5, 14, 18, 23, 27\}:\ 24,Y\\
\{ 4, 5, 10, 14, 27, 31\}:\ 26,M &
\{ 4, 5, 11, 15, 26, 30\}:\ 25,M &
\{ 4, 5, 14, 18, 23, 27\}:\ 25,M\\
\{ 4, 10, 14, 23, 27\}:\ 24,Y &
\{ 5, 10, 14, 23, 27\}:\ 23,Y &
\{ 4, 5, 10, 14, 23, 27\}:\ 23,M\\
\hline
\caption*{\small Potentially non-affine
$D^\dag$ for $\Ga=\lan 4,6,21\ran$}
\end{longtable}
}
\hspace {0.1cm}

\subsection {\bf %The case
\texorpdfstring{{\mathversion{bold}$\r\!=\!\C[[z^4,z^6\!+\!z^7]],
\ell\!=\!0$}}
{The case v=7, length=0}}
\subsubsection {\sf Potentially non-affine cells}
Let us begin with the list of all {\em potentially non-affine
cells\,} for $\r=\C[[z^4,z^6+z^7]]$ and their dimensions,
i.e. when the affineness cannot be established
using {\em straightforward elimination\,} of
the $\la$\~variables.

\hskip -0.7cm
{\tiny
\begin{longtable}{|l|l|l|}
\caption*
{\small
Types of cells for $\r=\C[[z^4,z^6+z^{7}]]$ when $\ell=0$:}\\
\hline
\{ 4, 14, 18, 31\}:\ 30,A &
\{ 5, 14, 18\}:\ 29,X &
\{ 4, 14, 18, 31\}:\ 28,A\\
\{ 4, 5, 14, 18\}:\ 30,X &
\{ 4, 14, 18, 19\}:\ 29,A &
\{ 4, 14, 18, 23\}:\ 28,A\\
\{ 5, 14, 18\}:\ 27,X &
\{ 4, 5, 14, 18\}:\ 29,X &
\{ 5, 14, 18, 19\}:\ 28,X\\
\{ 5, 14, 18, 23\}:\ 27,X &
\{ 4, 14, 18, 23\}:\ 27,A &
\{ 5, 14, 18, 23\}:\ 26,X\\
\{ 4, 5, 14, 18\}:\ 27,X &
\{ 4, 5, 14, 18, 19\}:\ 28,X &
\{ 4, 5, 14, 18, 23\}:\ 27,X\\
\{ 4, 14, 18, 23\}:\ 25,A &
\{ 4, 5, 14, 18, 23\}:\ 28,X &
\{ 4, 10, 14, 23\}:\ 27,A\\
\{ 4, 14, 15, 18\}:\ 26,A &
\{ 4, 14, 18, 19, 23\}:\ 25,A &
\{ 5, 14, 18, 23\}:\ 24,X\\
\{ 4, 5, 14, 18, 23\}:\ 27,X &
\{ 5, 10, 14, 23\}:\ 26,X &
\{ 5, 14, 15, 18\}:\ 25,X\\
\{ 5, 14, 18, 19, 23\}:\ 24,X &
\{ 4, 10, 14, 23\}:\ 26,A &
\{ 5, 10, 14, 23\}:\ 25,X\\
\{ 4, 5, 14, 18, 23\}:\ 25,X &
\{ 4, 5, 10, 14, 23\}:\ 27,X &
\{ 4, 5, 14, 15, 18\}:\ 26,X\\
\{ 4, 5, 14, 18, 19, 23\}:\ 25,X &
\{ 4, 10, 14, 23\}:\ 24,A &
\{ 4, 5, 10, 14, 23\}:\ 26,X\\
\{ 4, 10, 14, 15\}:\ 25,A &
\{ 4, 10, 14, 19, 23\}:\ 24,A &
\{ 5, 10, 14, 23\}:\ 23,X\\
\{ 4, 5, 10, 14, 23\}:\ 25,X &
\{ 5, 10, 14, 15\}:\ 24,X &
\{ 5, 10, 14, 19, 23\}:\ 23,X\\
\{ 4, 5, 10, 14, 23\}:\ 23,X &
\{ 4, 5, 10, 14, 15\}:\ 24,X &
\{ 4, 5, 10, 14, 19, 23\}:\ 23,X\\
\hline
\caption*{\small Potentially non-affine
$D^\dag$ for $\Ga=\lan 4,6,13\ran$}
\end{longtable}
}
\hspace {0.1cm}
\label{Table4-6-13-na}

\subsubsection {\sf Dimensions and deviations}
Let us give the dimensions and deviations $\kap$
of the cells for
$\r=\C[[z^4,z^6+z^7]]$ for single modules, i.e. when $\ell=0$,
$\kap=|D|$. In
this example, the dimensions always coincide with potential
dimensions $\partial-d$, which is not true for
$\r=\C[[z^4,z^6+z^9]]$ and any further members of this family.
After {\em primitive\,} $D^\dag$, we put  the dimension of the
corresponding cell and then the deviation $|D|$ (which can not
be readily seen since with use primitive $D^\dag$). The ordering
is from left to right and then downward. The symbol $na$ means
that the flag is non-admissible. Recall that $\Ga=\lan
4,6,13\ran$, $\de=8$, $\partial=4\de=32$ and the maximal
deviation is $\kap=2\de=16$ in this case. This table must be used
together with the list of non-affine cells, which are all of
type $A$ or $X$ (see above).

{\tiny
\hskip -0.7cm
\begin{longtable}{|l|l|l|}
\caption*
{\small
Dimensions and deviations for $\r=\C[[z^4,z^6+z^{7}]]$,
$\ell=0$:}\\
\hline
\{\ \}\ :\ 32,0 \ $\rightsquigarrow$&
\{ 30\}\ :\ 31,1 \ $\rightsquigarrow$&
\{ 31\}\ :\ 30,1 \ $\rightsquigarrow$\\
\{ 4\}\ :\ na,2 &
\{ 5\}\ :\ na,2 &
\{ 18\}\ :\ 31,2\\
\{ 19\}\ :\ 30,2 &
\{ 22\}\ :\ 30,2 &
\{ 23\}\ :\ 29,2\\
\{ 30, 31\}\ :\ 28,2 &
\{ 4, 18\}\ :\ 31,3 &
\{ 4, 22\}\ :\ na,3\\
\{ 4, 31\}\ :\ na,3 &
\{ 5, 19\}\ :\ 30,3 &
\{ 5, 23\}\ :\ na,3\\
\{ 5, 30\}\ :\ na,3 &
\{ 14\}\ :\ 30,3 &
\{ 15\}\ :\ 29,3\\
\{ 18, 22\}\ :\ 29,3 &
\{ 18, 31\}\ :\ 29,3 &
\{ 19, 23\}\ :\ 28,3\\
\{ 19, 30\}\ :\ 28,3 &
\{ 22, 31\}\ :\ 27,3 &
\{ 23, 30\}\ :\ 26,3\\
\{ 4, 5\}\ :\ na,4 &
\{ 4, 14\}\ :\ 30,4 &
\{ 4, 18, 22\}\ :\ 30,4\\
\{ 4, 18, 31\}\ :\ 29,4 &
\{ 4, 19\}\ :\ na,4 &
\{ 4, 22, 31\}\ :\ na,4\\
\{ 4, 23\}\ :\ na,4 &
\{ 5, 15\}\ :\ 29,4 &
\{ 5, 18\}\ :\ na,4\\
\{ 5, 19, 23\}\ :\ 29,4 &
\{ 5, 19, 30\}\ :\ 28,4 &
\{ 5, 22\}\ :\ na,4\\
\{ 5, 23, 30\}\ :\ na,4 &
\{ 10\}\ :\ 29,4 &
\{ 11\}\ :\ 28,4\\
\{ 14, 18\}\ :\ 28,4 &
\{ 14, 31\}\ :\ 28,4 &
\{ 15, 19\}\ :\ 27,4\\
\{ 15, 30\}\ :\ 27,4 &
\{ 18, 19\}\ :\ 28,4 &
\{ 18, 22, 31\}\ :\ 27,4\\
\{ 18, 23\}\ :\ 26,4 &
\{ 19, 22\}\ :\ 26,4 &
\{ 19, 23, 30\}\ :\ 25,4\\
\{ 22, 23\}\ :\ 24,4 &
\{ 4, 5, 18\}\ :\ na,5 &
\{ 4, 5, 19\}\ :\ na,5\\
\{ 4, 5, 22\}\ :\ na,5 &
\{ 4, 5, 23\}\ :\ na,5 &
\{ 4, 10\}\ :\ 29,5\\
\{ 4, 14, 18\}\ :\ 29,5 &
\{ 4, 14, 31\}\ :\ 28,5 &
\{ 4, 15\}\ :\ na,5\\
\{ 4, 18, 19\}\ :\ 29,5 &
\{ 4, 18, 22, 31\}\ :\ 28,5 &
\{ 4, 18, 23\}\ :\ 27,5\\
\{ 4, 19, 22\}\ :\ na,5 &
\{ 4, 19, 23\}\ :\ na,5 &
\{ 4, 22, 23\}\ :\ na,5\\
\{ 5, 11\}\ :\ 28,5 &
\{ 5, 14\}\ :\ na,5 &
\{ 5, 15, 19\}\ :\ 28,5\\
\{ 5, 15, 30\}\ :\ 27,5 &
\{ 5, 18, 19\}\ :\ 28,5 &
\{ 5, 18, 22\}\ :\ na,5\\
\{ 5, 18, 23\}\ :\ na,5 &
\{ 5, 19, 22\}\ :\ 27,5 &
\{ 5, 19, 23, 30\}\ :\ 26,5\\
\{ 5, 22, 23\}\ :\ na,5 &
\{ 6\}\ :\ 28,5 &
\{ 7\}\ :\ 27,5\\
\{ 10, 14\}\ :\ 27,5 &
\{ 10, 31\}\ :\ 27,5 &
\{ 11, 15\}\ :\ 26,5\\
\{ 11, 30\}\ :\ 26,5 &
\{ 14, 18, 31\}\ :\ 27,5 &
\{ 14, 19\}\ :\ 26,5\\
\{ 14, 23\}\ :\ 25,5 &
\{ 15, 18\}\ :\ 26,5 &
\{ 15, 19, 30\}\ :\ 25,5\\
\{ 15, 22\}\ :\ 24,5 &
\{ 18, 19, 22\}\ :\ 25,5 &
\{ 18, 19, 23\}\ :\ 24,5\\
\{ 18, 22, 23\}\ :\ 23,5 &
\{ 19, 22, 23\}\ :\ 22,5 &
\{ 2\}\ :\ 28,6\\
\{ 3\}\ :\ 27,6 &
\{ 4, 5, 14\}\ :\ na,6 &
\{ 4, 5, 15\}\ :\ na,6\\
\{ 4, 5, 18, 19\}\ :\ 28,6 &
\{ 4, 5, 18, 22\}\ :\ na,6 &
\{ 4, 5, 18, 23\}\ :\ na,6\\
\{ 4, 5, 19, 22\}\ :\ na,6 &
\{ 4, 5, 19, 23\}\ :\ na,6 &
\{ 4, 5, 22, 23\}\ :\ na,6\\
\{ 4, 6\}\ :\ 28,6 &
\{ 4, 10, 14\}\ :\ 27,6 &
\{ 4, 10, 31\}\ :\ 27,6\\
\{ 4, 11\}\ :\ na,6 &
\{ 4, 14, 18, 31\}\ :\ 28,6 &
\{ 4, 14, 19\}\ :\ 27,6\\
\{ 4, 14, 23\}\ :\ 26,6 &
\{ 4, 15, 18\}\ :\ 27,6 &
\{ 4, 15, 19\}\ :\ na,6\\
\{ 4, 15, 22\}\ :\ na,6 &
\{ 4, 18, 19, 22\}\ :\ 27,6 &
\{ 4, 18, 19, 23\}\ :\ 26,6\\
\{ 4, 18, 22, 23\}\ :\ 25,6 &
\{ 4, 19, 22, 23\}\ :\ na,6 &
\{ 5, 7\}\ :\ 27,6\\
\{ 5, 10\}\ :\ na,6 &
\{ 5, 11, 15\}\ :\ 26,6 &
\{ 5, 11, 30\}\ :\ 26,6\\
\{ 5, 14, 18\}\ :\ 27,6 &
\{ 5, 14, 19\}\ :\ 26,6 &
\{ 5, 14, 23\}\ :\ na,6\\
\{ 5, 15, 18\}\ :\ 26,6 &
\{ 5, 15, 19, 30\}\ :\ 26,6 &
\{ 5, 15, 22\}\ :\ 25,6\\
\{ 5, 18, 19, 22\}\ :\ 26,6 &
\{ 5, 18, 19, 23\}\ :\ 25,6 &
\{ 5, 18, 22, 23\}\ :\ na,6\\
\{ 5, 19, 22, 23\}\ :\ 24,6 &
\{ 6, 10\}\ :\ 26,6 &
\{ 6, 31\}\ :\ 26,6\\
\{ 7, 11\}\ :\ 25,6 &
\{ 7, 30\}\ :\ 25,6 &
\{ 10, 14, 31\}\ :\ 26,6\\
\{ 10, 19\}\ :\ 25,6 &
\{ 10, 23\}\ :\ 24,6 &
\{ 11, 15, 30\}\ :\ 25,6\\
\{ 11, 18\}\ :\ 24,6 &
\{ 11, 22\}\ :\ 23,6 &
\{ 14, 15\}\ :\ 24,6\\
\{ 14, 18, 19\}\ :\ 24,6 &
\{ 14, 18, 23\}\ :\ 23,6 &
\{ 14, 19, 23\}\ :\ 22,6\\
\{ 15, 18, 19\}\ :\ 23,6 &
\{ 15, 18, 22\}\ :\ 22,6 &
\{ 15, 19, 22\}\ :\ 21,6\\
\{ 18, 19, 22, 23\}\ :\ 20,6 &
\{ 2, 4\}\ :\ 27,7 &
\{ 2, 6\}\ :\ 26,7\\
\{ 2, 31\}\ :\ 27,7 &
\{ 3, 5\}\ :\ 26,7 &
\{ 3, 7\}\ :\ 25,7\\
\{ 3, 30\}\ :\ 26,7 &
\{ 4, 5, 10\}\ :\ na,7 &
\{ 4, 5, 11\}\ :\ na,7\\
\{ 4, 5, 14, 18\}\ :\ 27,7 &
\{ 4, 5, 14, 19\}\ :\ 26,7 &
\{ 4, 5, 14, 23\}\ :\ na,7\\
\{ 4, 5, 15, 18\}\ :\ 26,7 &
\{ 4, 5, 15, 19\}\ :\ na,7 &
\{ 4, 5, 15, 22\}\ :\ na,7\\
\{ 4, 5, 18, 19, 22\}\ :\ 27,7 &
\{ 4, 5, 18, 19, 23\}\ :\ 26,7 &
\{ 4, 5, 18, 22, 23\}\ :\ na,7\\
\{ 4, 5, 19, 22, 23\}\ :\ na,7 &
\{ 4, 6, 10\}\ :\ 25,7 &
\{ 4, 6, 31\}\ :\ 26,7\\
\{ 4, 7\}\ :\ na,7 &
\{ 4, 10, 14, 31\}\ :\ 26,7 &
\{ 4, 10, 19\}\ :\ 26,7\\
\{ 4, 10, 23\}\ :\ 25,7 &
\{ 4, 11, 15\}\ :\ na,7 &
\{ 4, 11, 18\}\ :\ 25,7\\
\{ 4, 11, 22\}\ :\ na,7 &
\{ 4, 14, 15\}\ :\ 25,7 &
\{ 4, 14, 18, 19\}\ :\ 26,7\\
\{ 4, 14, 18, 23\}\ :\ 25,7 &
\{ 4, 14, 19, 23\}\ :\ 24,7 &
\{ 4, 15, 18, 19\}\ :\ 25,7\\
\{ 4, 15, 18, 22\}\ :\ 24,7 &
\{ 4, 15, 19, 22\}\ :\ na,7 &
\{ 4, 18, 19, 22, 23\}\ :\ 23,7\\
\{ 5, 6\}\ :\ na,7 &
\{ 5, 7, 11\}\ :\ 24,7 &
\{ 5, 7, 30\}\ :\ 25,7\\
\{ 5, 10, 14\}\ :\ na,7 &
\{ 5, 10, 19\}\ :\ 25,7 &
\{ 5, 10, 23\}\ :\ na,7\\
\{ 5, 11, 15, 30\}\ :\ 25,7 &
\{ 5, 11, 18\}\ :\ 24,7 &
\{ 5, 11, 22\}\ :\ 24,7\\
\{ 5, 14, 15\}\ :\ 24,7 &
\{ 5, 14, 18, 19\}\ :\ 25,7 &
\{ 5, 14, 18, 23\}\ :\ 24,7\\
\{ 5, 14, 19, 23\}\ :\ 23,7 &
\{ 5, 15, 18, 19\}\ :\ 24,7 &
\{ 5, 15, 18, 22\}\ :\ 23,7\\
\{ 5, 15, 19, 22\}\ :\ 23,7 &
\{ 5, 18, 19, 22, 23\}\ :\ 22,7 &
\{ 6, 10, 31\}\ :\ 25,7\\
\{ 6, 19\}\ :\ 24,7 &
\{ 6, 23\}\ :\ 23,7 &
\{ 7, 11, 30\}\ :\ 24,7\\
\{ 7, 18\}\ :\ 23,7 &
\{ 7, 22\}\ :\ 22,7 &
\{ 10, 14, 19\}\ :\ 23,7\\
\{ 10, 14, 23\}\ :\ 23,7 &
\{ 10, 15\}\ :\ 22,7 &
\{ 10, 19, 23\}\ :\ 21,7\\
\{ 11, 14\}\ :\ 22,7 &
\{ 11, 15, 18\}\ :\ 22,7 &
\{ 11, 15, 22\}\ :\ 21,7\\
\{ 11, 18, 22\}\ :\ 20,7 &
\{ 14, 15, 18\}\ :\ 21,7 &
\{ 14, 15, 19\}\ :\ 20,7\\
\{ 14, 18, 19, 23\}\ :\ 19,7 &
\{ 15, 18, 19, 22\}\ :\ 18,7 &
\{ 2, 4, 6\}\ :\ 24,8\\
\{ 2, 4, 31\}\ :\ 26,8 &
\{ 2, 5\}\ :\ na,8 &
\{ 2, 6, 31\}\ :\ 25,8\\
\{ 2, 19\}\ :\ 25,8 &
\{ 2, 23\}\ :\ 24,8 &
\{ 3, 4\}\ :\ na,8\\
\{ 3, 5, 7\}\ :\ 23,8 &
\{ 3, 5, 30\}\ :\ 25,8 &
\{ 3, 7, 30\}\ :\ 24,8\\
\{ 3, 18\}\ :\ 24,8 &
\{ 3, 22\}\ :\ 23,8 &
\{ 4, 5, 6\}\ :\ na,8\\
\{ 4, 5, 7\}\ :\ na,8 &
\{ 4, 5, 10, 14\}\ :\ na,8 &
\{ 4, 5, 10, 19\}\ :\ 25,8\\
\{ 4, 5, 10, 23\}\ :\ na,8 &
\{ 4, 5, 11, 15\}\ :\ na,8 &
\{ 4, 5, 11, 18\}\ :\ 24,8\\
\{ 4, 5, 11, 22\}\ :\ na,8 &
\{ 4, 5, 14, 15\}\ :\ 24,8 &
\{ 4, 5, 14, 18, 19\}\ :\ 26,8\\
\{ 4, 5, 14, 18, 23\}\ :\ 25,8 &
\{ 4, 5, 14, 19, 23\}\ :\ 24,8 &
\{ 4, 5, 15, 18, 19\}\ :\ 25,8\\
\{ 4, 5, 15, 18, 22\}\ :\ 24,8 &
\{ 4, 5, 15, 19, 22\}\ :\ na,8 &
\{ 4, 5, 18, 19, 22, 23\}\ :\ 24,8\\
\{ 4, 6, 10, 31\}\ :\ 24,8 &
\{ 4, 6, 19\}\ :\ 25,8 &
\{ 4, 6, 23\}\ :\ 24,8\\
\{ 4, 7, 11\}\ :\ na,8 &
\{ 4, 7, 18\}\ :\ 24,8 &
\{ 4, 7, 22\}\ :\ na,8\\
\{ 4, 10, 14, 19\}\ :\ 24,8 &
\{ 4, 10, 14, 23\}\ :\ 24,8 &
\{ 4, 10, 15\}\ :\ 23,8\\
\{ 4, 10, 19, 23\}\ :\ 23,8 &
\{ 4, 11, 14\}\ :\ 23,8 &
\{ 4, 11, 15, 18\}\ :\ 23,8\\
\{ 4, 11, 15, 22\}\ :\ na,8 &
\{ 4, 11, 18, 22\}\ :\ 22,8 &
\{ 4, 14, 15, 18\}\ :\ 23,8\\
\{ 4, 14, 15, 19\}\ :\ 22,8 &
\{ 4, 14, 18, 19, 23\}\ :\ 22,8 &
\{ 4, 15, 18, 19, 22\}\ :\ 21,8\\
\{ 5, 6, 10\}\ :\ na,8 &
\{ 5, 6, 19\}\ :\ 24,8 &
\{ 5, 6, 23\}\ :\ na,8\\
\{ 5, 7, 11, 30\}\ :\ 23,8 &
\{ 5, 7, 18\}\ :\ 23,8 &
\{ 5, 7, 22\}\ :\ 23,8\\
\{ 5, 10, 14, 19\}\ :\ 23,8 &
\{ 5, 10, 14, 23\}\ :\ 23,8 &
\{ 5, 10, 15\}\ :\ 22,8\\
\{ 5, 10, 19, 23\}\ :\ 22,8 &
\{ 5, 11, 14\}\ :\ 22,8 &
\{ 5, 11, 15, 18\}\ :\ 22,8\\
\{ 5, 11, 15, 22\}\ :\ 22,8 &
\{ 5, 11, 18, 22\}\ :\ 21,8 &
\{ 5, 14, 15, 18\}\ :\ 22,8\\
\{ 5, 14, 15, 19\}\ :\ 21,8 &
\{ 5, 14, 18, 19, 23\}\ :\ 21,8 &
\{ 5, 15, 18, 19, 22\}\ :\ 20,8\\
\{ 6, 10, 19\}\ :\ 22,8 &
\{ 6, 10, 23\}\ :\ 22,8 &
\{ 6, 15\}\ :\ 21,8\\
\{ 6, 19, 23\}\ :\ 20,8 &
\{ 7, 11, 18\}\ :\ 21,8 &
\{ 7, 11, 22\}\ :\ 21,8\\
\{ 7, 14\}\ :\ 20,8 &
\{ 7, 18, 22\}\ :\ 19,8 &
\{ 10, 11\}\ :\ 20,8\\
\{ 10, 14, 15\}\ :\ 20,8 &
\{ 10, 14, 19, 23\}\ :\ 19,8 &
\{ 10, 15, 19\}\ :\ 18,8\\
\{ 11, 14, 15\}\ :\ 19,8 &
\{ 11, 14, 18\}\ :\ 18,8 &
\{ 11, 15, 18, 22\}\ :\ 17,8\\
\{ 14, 15, 18, 19\}\ :\ 16,8 &
\{ 2, 4, 5\}\ :\ na,9 &
\{ 2, 4, 6, 31\}\ :\ 23,9\\
\{ 2, 4, 19\}\ :\ 25,9 &
\{ 2, 4, 23\}\ :\ 24,9 &
\{ 2, 5, 6\}\ :\ na,9\\
\{ 2, 5, 19\}\ :\ 24,9 &
\{ 2, 5, 23\}\ :\ na,9 &
\{ 2, 6, 19\}\ :\ 23,9\\
\{ 2, 6, 23\}\ :\ 23,9 &
\{ 2, 15\}\ :\ 22,9 &
\{ 2, 19, 23\}\ :\ 21,9\\
\{ 3, 4, 5\}\ :\ na,9 &
\{ 3, 4, 7\}\ :\ na,9 &
\{ 3, 4, 18\}\ :\ 24,9\\
\{ 3, 4, 22\}\ :\ na,9 &
\{ 3, 5, 7, 30\}\ :\ 22,9 &
\{ 3, 5, 18\}\ :\ 23,9\\
\{ 3, 5, 22\}\ :\ 23,9 &
\{ 3, 7, 18\}\ :\ 22,9 &
\{ 3, 7, 22\}\ :\ 22,9\\
\{ 3, 14\}\ :\ 21,9 &
\{ 3, 18, 22\}\ :\ 20,9 &
\{ 4, 5, 6, 10\}\ :\ na,9\\
\{ 4, 5, 6, 19\}\ :\ 24,9 &
\{ 4, 5, 6, 23\}\ :\ na,9 &
\{ 4, 5, 7, 11\}\ :\ na,9\\
\{ 4, 5, 7, 18\}\ :\ 23,9 &
\{ 4, 5, 7, 22\}\ :\ na,9 &
\{ 4, 5, 10, 14, 19\}\ :\ 23,9\\
\{ 4, 5, 10, 14, 23\}\ :\ 23,9 &
\{ 4, 5, 10, 15\}\ :\ 22,9 &
\{ 4, 5, 10, 19, 23\}\ :\ 23,9\\
\{ 4, 5, 11, 14\}\ :\ 22,9 &
\{ 4, 5, 11, 15, 18\}\ :\ 22,9 &
\{ 4, 5, 11, 15, 22\}\ :\ na,9\\
\{ 4, 5, 11, 18, 22\}\ :\ 22,9 &
\{ 4, 5, 14, 15, 18\}\ :\ 23,9 &
\{ 4, 5, 14, 15, 19\}\ :\ 22,9\\
\{ 4, 5, 14, 18, 19, 23\}\ :\ 23,9 &
\{ 4, 5, 15, 18, 19, 22\}\ :\ 22,9 &
\{ 4, 6, 10, 19\}\ :\ 22,9\\
\{ 4, 6, 10, 23\}\ :\ 22,9 &
\{ 4, 6, 15\}\ :\ 22,9 &
\{ 4, 6, 19, 23\}\ :\ 22,9\\
\{ 4, 7, 11, 18\}\ :\ 21,9 &
\{ 4, 7, 11, 22\}\ :\ na,9 &
\{ 4, 7, 14\}\ :\ 21,9\\
\{ 4, 7, 18, 22\}\ :\ 21,9 &
\{ 4, 10, 11\}\ :\ 21,9 &
\{ 4, 10, 14, 15\}\ :\ 21,9\\
\{ 4, 10, 14, 19, 23\}\ :\ 21,9 &
\{ 4, 10, 15, 19\}\ :\ 20,9 &
\{ 4, 11, 14, 15\}\ :\ 20,9\\
\{ 4, 11, 14, 18\}\ :\ 20,9 &
\{ 4, 11, 15, 18, 22\}\ :\ 19,9 &
\{ 4, 14, 15, 18, 19\}\ :\ 19,9\\
\{ 5, 6, 10, 19\}\ :\ 21,9 &
\{ 5, 6, 10, 23\}\ :\ na,9 &
\{ 5, 6, 15\}\ :\ 21,9\\
\{ 5, 6, 19, 23\}\ :\ 21,9 &
\{ 5, 7, 11, 18\}\ :\ 20,9 &
\{ 5, 7, 11, 22\}\ :\ 21,9\\
\{ 5, 7, 14\}\ :\ 20,9 &
\{ 5, 7, 18, 22\}\ :\ 20,9 &
\{ 5, 10, 11\}\ :\ 20,9\\
\{ 5, 10, 14, 15\}\ :\ 20,9 &
\{ 5, 10, 14, 19, 23\}\ :\ 20,9 &
\{ 5, 10, 15, 19\}\ :\ 19,9\\
\{ 5, 11, 14, 15\}\ :\ 19,9 &
\{ 5, 11, 14, 18\}\ :\ 19,9 &
\{ 5, 11, 15, 18, 22\}\ :\ 18,9\\
\{ 5, 14, 15, 18, 19\}\ :\ 18,9 &
\{ 6, 10, 15\}\ :\ 19,9 &
\{ 6, 10, 19, 23\}\ :\ 19,9\\
\{ 6, 11\}\ :\ 18,9 &
\{ 6, 15, 19\}\ :\ 17,9 &
\{ 7, 10\}\ :\ 18,9\\
\{ 7, 11, 14\}\ :\ 18,9 &
\{ 7, 11, 18, 22\}\ :\ 17,9 &
\{ 7, 14, 18\}\ :\ 16,9\\
\{ 10, 11, 14\}\ :\ 17,9 &
\{ 10, 11, 15\}\ :\ 16,9 &
\{ 10, 14, 15, 19\}\ :\ 15,9\\
\{ 11, 14, 15, 18\}\ :\ 14,9 &
\{ 2, 4, 5, 6\}\ :\ na,10 &
\{ 2, 4, 5, 19\}\ :\ 23,10\\
\{ 2, 4, 5, 23\}\ :\ na,10 &
\{ 2, 4, 6, 19\}\ :\ 22,10 &
\{ 2, 4, 6, 23\}\ :\ 22,10\\
\{ 2, 4, 15\}\ :\ 22,10 &
\{ 2, 4, 19, 23\}\ :\ 22,10 &
\{ 2, 5, 6, 19\}\ :\ 21,10\\
\{ 2, 5, 6, 23\}\ :\ na,10 &
\{ 2, 5, 15\}\ :\ 21,10 &
\{ 2, 5, 19, 23\}\ :\ 21,10\\
\{ 2, 6, 15\}\ :\ 20,10 &
\{ 2, 6, 19, 23\}\ :\ 20,10 &
\{ 2, 11\}\ :\ 19,10\\
\{ 2, 15, 19\}\ :\ 18,10 &
\{ 3, 4, 5, 7\}\ :\ na,10 &
\{ 3, 4, 5, 18\}\ :\ 22,10\\
\{ 3, 4, 5, 22\}\ :\ na,10 &
\{ 3, 4, 7, 18\}\ :\ 21,10 &
\{ 3, 4, 7, 22\}\ :\ na,10\\
\{ 3, 4, 14\}\ :\ 21,10 &
\{ 3, 4, 18, 22\}\ :\ 21,10 &
\{ 3, 5, 7, 18\}\ :\ 20,10\\
\{ 3, 5, 7, 22\}\ :\ 21,10 &
\{ 3, 5, 14\}\ :\ 20,10 &
\{ 3, 5, 18, 22\}\ :\ 20,10\\
\{ 3, 7, 14\}\ :\ 19,10 &
\{ 3, 7, 18, 22\}\ :\ 19,10 &
\{ 3, 10\}\ :\ 18,10\\
\{ 3, 14, 18\}\ :\ 17,10 &
\{ 4, 5, 6, 10, 19\}\ :\ 20,10 &
\{ 4, 5, 6, 10, 23\}\ :\ na,10\\
\{ 4, 5, 6, 15\}\ :\ 21,10 &
\{ 4, 5, 6, 19, 23\}\ :\ 22,10 &
\{ 4, 5, 7, 11, 18\}\ :\ 19,10\\
\{ 4, 5, 7, 11, 22\}\ :\ na,10 &
\{ 4, 5, 7, 14\}\ :\ 20,10 &
\{ 4, 5, 7, 18, 22\}\ :\ 21,10\\
\{ 4, 5, 10, 11\}\ :\ 20,10 &
\{ 4, 5, 10, 14, 15\}\ :\ 20,10 &
\{ 4, 5, 10, 14, 19, 23\}\ :\ 21,10\\
\{ 4, 5, 10, 15, 19\}\ :\ 20,10 &
\{ 4, 5, 11, 14, 15\}\ :\ 19,10 &
\{ 4, 5, 11, 14, 18\}\ :\ 20,10\\
\{ 4, 5, 11, 15, 18, 22\}\ :\ 19,10 &
\{ 4, 5, 14, 15, 18, 19\}\ :\ 20,10 &
\{ 4, 6, 10, 15\}\ :\ 19,10\\
\{ 4, 6, 10, 19, 23\}\ :\ 20,10 &
\{ 4, 6, 11\}\ :\ 19,10 &
\{ 4, 6, 15, 19\}\ :\ 19,10\\
\{ 4, 7, 10\}\ :\ 19,10 &
\{ 4, 7, 11, 14\}\ :\ 18,10 &
\{ 4, 7, 11, 18, 22\}\ :\ 18,10\\
\{ 4, 7, 14, 18\}\ :\ 18,10 &
\{ 4, 10, 11, 14\}\ :\ 18,10 &
\{ 4, 10, 11, 15\}\ :\ 17,10\\
\{ 4, 10, 14, 15, 19\}\ :\ 17,10 &
\{ 4, 11, 14, 15, 18\}\ :\ 16,10 &
\{ 5, 6, 10, 15\}\ :\ 18,10\\
\{ 5, 6, 10, 19, 23\}\ :\ 19,10 &
\{ 5, 6, 11\}\ :\ 18,10 &
\{ 5, 6, 15, 19\}\ :\ 18,10\\
\{ 5, 7, 10\}\ :\ 18,10 &
\{ 5, 7, 11, 14\}\ :\ 17,10 &
\{ 5, 7, 11, 18, 22\}\ :\ 17,10\\
\{ 5, 7, 14, 18\}\ :\ 17,10 &
\{ 5, 10, 11, 14\}\ :\ 17,10 &
\{ 5, 10, 11, 15\}\ :\ 16,10\\
\{ 5, 10, 14, 15, 19\}\ :\ 16,10 &
\{ 5, 11, 14, 15, 18\}\ :\ 15,10 &
\{ 6, 7\}\ :\ 16,10\\
\{ 6, 10, 11\}\ :\ 16,10 &
\{ 6, 10, 15, 19\}\ :\ 15,10 &
\{ 6, 11, 15\}\ :\ 14,10\\
\{ 7, 10, 11\}\ :\ 15,10 &
\{ 7, 10, 14\}\ :\ 14,10 &
\{ 7, 11, 14, 18\}\ :\ 13,10\\
\{ 10, 11, 14, 15\}\ :\ 12,10 &
\{ 2, 4, 5, 6, 19\}\ :\ 19,11 &
\{ 2, 4, 5, 6, 23\}\ :\ na,11\\
\{ 2, 4, 5, 15\}\ :\ 20,11 &
\{ 2, 4, 5, 19, 23\}\ :\ 21,11 &
\{ 2, 4, 6, 15\}\ :\ 19,11\\
\{ 2, 4, 6, 19, 23\}\ :\ 20,11 &
\{ 2, 4, 11\}\ :\ 19,11 &
\{ 2, 4, 15, 19\}\ :\ 19,11\\
\{ 2, 5, 6, 15\}\ :\ 18,11 &
\{ 2, 5, 6, 19, 23\}\ :\ 19,11 &
\{ 2, 5, 11\}\ :\ 18,11\\
\{ 2, 5, 15, 19\}\ :\ 18,11 &
\{ 2, 6, 11\}\ :\ 17,11 &
\{ 2, 6, 15, 19\}\ :\ 17,11\\
\{ 2, 7\}\ :\ 16,11 &
\{ 2, 11, 15\}\ :\ 15,11 &
\{ 3, 4, 5, 7, 18\}\ :\ 18,11\\
\{ 3, 4, 5, 7, 22\}\ :\ na,11 &
\{ 3, 4, 5, 14\}\ :\ 19,11 &
\{ 3, 4, 5, 18, 22\}\ :\ 20,11\\
\{ 3, 4, 7, 14\}\ :\ 18,11 &
\{ 3, 4, 7, 18, 22\}\ :\ 19,11 &
\{ 3, 4, 10\}\ :\ 18,11\\
\{ 3, 4, 14, 18\}\ :\ 18,11 &
\{ 3, 5, 7, 14\}\ :\ 17,11 &
\{ 3, 5, 7, 18, 22\}\ :\ 18,11\\
\{ 3, 5, 10\}\ :\ 17,11 &
\{ 3, 5, 14, 18\}\ :\ 17,11 &
\{ 3, 6\}\ :\ 16,11\\
\{ 3, 7, 10\}\ :\ 16,11 &
\{ 3, 7, 14, 18\}\ :\ 15,11 &
\{ 3, 10, 14\}\ :\ 14,11\\
\{ 4, 5, 6, 10, 15\}\ :\ 17,11 &
\{ 4, 5, 6, 10, 19, 23\}\ :\ 19,11 &
\{ 4, 5, 6, 11\}\ :\ 18,11\\
\{ 4, 5, 6, 15, 19\}\ :\ 19,11 &
\{ 4, 5, 7, 10\}\ :\ 18,11 &
\{ 4, 5, 7, 11, 14\}\ :\ 16,11\\
\{ 4, 5, 7, 11, 18, 22\}\ :\ 17,11 &
\{ 4, 5, 7, 14, 18\}\ :\ 18,11 &
\{ 4, 5, 10, 11, 14\}\ :\ 17,11\\
\{ 4, 5, 10, 11, 15\}\ :\ 16,11 &
\{ 4, 5, 10, 14, 15, 19\}\ :\ 17,11 &
\{ 4, 5, 11, 14, 15, 18\}\ :\ 16,11\\
\{ 4, 6, 7\}\ :\ 17,11 &
\{ 4, 6, 10, 11\}\ :\ 16,11 &
\{ 4, 6, 10, 15, 19\}\ :\ 16,11\\
\{ 4, 6, 11, 15\}\ :\ 15,11 &
\{ 4, 7, 10, 11\}\ :\ 15,11 &
\{ 4, 7, 10, 14\}\ :\ 15,11\\
\{ 4, 7, 11, 14, 18\}\ :\ 14,11 &
\{ 4, 10, 11, 14, 15\}\ :\ 13,11 &
\{ 5, 6, 7\}\ :\ 16,11\\
\{ 5, 6, 10, 11\}\ :\ 15,11 &
\{ 5, 6, 10, 15, 19\}\ :\ 15,11 &
\{ 5, 6, 11, 15\}\ :\ 14,11\\
\{ 5, 7, 10, 11\}\ :\ 14,11 &
\{ 5, 7, 10, 14\}\ :\ 14,11 &
\{ 5, 7, 11, 14, 18\}\ :\ 13,11\\
\{ 5, 10, 11, 14, 15\}\ :\ 12,11 &
\{ 6, 7, 10\}\ :\ 13,11 &
\{ 6, 7, 11\}\ :\ 12,11\\
\{ 6, 10, 11, 15\}\ :\ 11,11 &
\{ 7, 10, 11, 14\}\ :\ 10,11 &
\{ 2, 3\}\ :\ 16,12\\
\{ 2, 4, 5, 6, 15\}\ :\ 16,12 &
\{ 2, 4, 5, 6, 19, 23\}\ :\ 18,12 &
\{ 2, 4, 5, 11\}\ :\ 17,12\\
\{ 2, 4, 5, 15, 19\}\ :\ 18,12 &
\{ 2, 4, 6, 11\}\ :\ 16,12 &
\{ 2, 4, 6, 15, 19\}\ :\ 17,12\\
\{ 2, 4, 7\}\ :\ 16,12 &
\{ 2, 4, 11, 15\}\ :\ 15,12 &
\{ 2, 5, 6, 11\}\ :\ 15,12\\
\{ 2, 5, 6, 15, 19\}\ :\ 16,12 &
\{ 2, 5, 7\}\ :\ 15,12 &
\{ 2, 5, 11, 15\}\ :\ 14,12\\
\{ 2, 6, 7\}\ :\ 14,12 &
\{ 2, 6, 11, 15\}\ :\ 13,12 &
\{ 2, 7, 11\}\ :\ 12,12\\
\{ 3, 4, 5, 7, 14\}\ :\ 15,12 &
\{ 3, 4, 5, 7, 18, 22\}\ :\ 17,12 &
\{ 3, 4, 5, 10\}\ :\ 16,12\\
\{ 3, 4, 5, 14, 18\}\ :\ 17,12 &
\{ 3, 4, 6\}\ :\ 16,12 &
\{ 3, 4, 7, 10\}\ :\ 15,12\\
\{ 3, 4, 7, 14, 18\}\ :\ 15,12 &
\{ 3, 4, 10, 14\}\ :\ 14,12 &
\{ 3, 5, 6\}\ :\ 15,12\\
\{ 3, 5, 7, 10\}\ :\ 14,12 &
\{ 3, 5, 7, 14, 18\}\ :\ 14,12 &
\{ 3, 5, 10, 14\}\ :\ 13,12\\
\{ 3, 6, 7\}\ :\ 13,12 &
\{ 3, 6, 10\}\ :\ 12,12 &
\{ 3, 7, 10, 14\}\ :\ 11,12\\
\{ 4, 5, 6, 7\}\ :\ 16,12 &
\{ 4, 5, 6, 10, 11\}\ :\ 14,12 &
\{ 4, 5, 6, 10, 15, 19\}\ :\ 15,12\\
\{ 4, 5, 6, 11, 15\}\ :\ 14,12 &
\{ 4, 5, 7, 10, 11\}\ :\ 13,12 &
\{ 4, 5, 7, 10, 14\}\ :\ 14,12\\
\{ 4, 5, 7, 11, 14, 18\}\ :\ 13,12 &
\{ 4, 5, 10, 11, 14, 15\}\ :\ 12,12 &
\{ 4, 6, 7, 10\}\ :\ 13,12\\
\{ 4, 6, 7, 11\}\ :\ 12,12 &
\{ 4, 6, 10, 11, 15\}\ :\ 11,12 &
\{ 4, 7, 10, 11, 14\}\ :\ 10,12\\
\{ 5, 6, 7, 10\}\ :\ 12,12 &
\{ 5, 6, 7, 11\}\ :\ 11,12 &
\{ 5, 6, 10, 11, 15\}\ :\ 10,12\\
\{ 5, 7, 10, 11, 14\}\ :\ 9,12 &
\{ 6, 7, 10, 11\}\ :\ 8,12 &
\{ 2, 3, 4\}\ :\ 15,13\\
\{ 2, 3, 5\}\ :\ 14,13 &
\{ 2, 3, 6\}\ :\ 13,13 &
\{ 2, 3, 7\}\ :\ 12,13\\
\{ 2, 4, 5, 6, 11\}\ :\ 13,13 &
\{ 2, 4, 5, 6, 15, 19\}\ :\ 15,13 &
\{ 2, 4, 5, 7\}\ :\ 14,13\\
\{ 2, 4, 5, 11, 15\}\ :\ 13,13 &
\{ 2, 4, 6, 7\}\ :\ 13,13 &
\{ 2, 4, 6, 11, 15\}\ :\ 12,13\\
\{ 2, 4, 7, 11\}\ :\ 11,13 &
\{ 2, 5, 6, 7\}\ :\ 12,13 &
\{ 2, 5, 6, 11, 15\}\ :\ 11,13\\
\{ 2, 5, 7, 11\}\ :\ 10,13 &
\{ 2, 6, 7, 11\}\ :\ 9,13 &
\{ 3, 4, 5, 6\}\ :\ 14,13\\
\{ 3, 4, 5, 7, 10\}\ :\ 12,13 &
\{ 3, 4, 5, 7, 14, 18\}\ :\ 13,13 &
\{ 3, 4, 5, 10, 14\}\ :\ 12,13\\
\{ 3, 4, 6, 7\}\ :\ 12,13 &
\{ 3, 4, 6, 10\}\ :\ 11,13 &
\{ 3, 4, 7, 10, 14\}\ :\ 10,13\\
\{ 3, 5, 6, 7\}\ :\ 11,13 &
\{ 3, 5, 6, 10\}\ :\ 10,13 &
\{ 3, 5, 7, 10, 14\}\ :\ 9,13\\
\{ 3, 6, 7, 10\}\ :\ 8,13 &
\{ 4, 5, 6, 7, 10\}\ :\ 11,13 &
\{ 4, 5, 6, 7, 11\}\ :\ 10,13\\
\{ 4, 5, 6, 10, 11, 15\}\ :\ 9,13 &
\{ 4, 5, 7, 10, 11, 14\}\ :\ 8,13 &
\{ 4, 6, 7, 10, 11\}\ :\ 7,13\\
\{ 5, 6, 7, 10, 11\}\ :\ 6,13 &
\{ 2, 3, 4, 5\}\ :\ 12,14 &
\{ 2, 3, 4, 6\}\ :\ 11,14\\
\{ 2, 3, 4, 7\}\ :\ 10,14 &
\{ 2, 3, 5, 6\}\ :\ 10,14 &
\{ 2, 3, 5, 7\}\ :\ 9,14\\
\{ 2, 3, 6, 7\}\ :\ 8,14 &
\{ 2, 4, 5, 6, 7\}\ :\ 10,14 &
\{ 2, 4, 5, 6, 11, 15\}\ :\ 9,14\\
\{ 2, 4, 5, 7, 11\}\ :\ 8,14 &
\{ 2, 4, 6, 7, 11\}\ :\ 7,14 &
\{ 2, 5, 6, 7, 11\}\ :\ 6,14\\
\{ 3, 4, 5, 6, 7\}\ :\ 9,14 &
\{ 3, 4, 5, 6, 10\}\ :\ 8,14 &
\{ 3, 4, 5, 7, 10, 14\}\ :\ 7,14\\
\{ 3, 4, 6, 7, 10\}\ :\ 6,14 &
\{ 3, 5, 6, 7, 10\}\ :\ 5,14 &
\{ 4, 5, 6, 7, 10, 11\}\ :\ 4,14\\
\{ 2, 3, 4, 5, 6\}\ :\ 7,15 &
\{ 2, 3, 4, 5, 7\}\ :\ 6,15 &
\{ 2, 3, 4, 6, 7\}\ :\ 5,15\\
\{ 2, 3, 5, 6, 7\}\ :\ 4,15 &
\{ 2, 4, 5, 6, 7, 11\}\ :\ 3,15 &
\{ 3, 4, 5, 6, 7, 10\}\ :\ 2,15\\
\{ 2, 3, 4, 5, 6, 7\}\ :\ 0,16 & &\\
\hline
\caption*{\small Dimensions and deviations for  $\rk=2,
\Ga=\lan 4,6,13\ran$}
\end{longtable}
\hspace {0.1cm}
}

\subsubsection{\sf The superpolynomial at
\texorpdfstring{$a\!=\!0$}{a=0}}
The formulas above used together
with the list of types (can be $A=\!$ affine or
$X$) give the complete portion $a=0$ of
the corresponding DAHA superpolynomial.
Recall that conjecturally the latter coincides with
the stable reduced Khovanov-Rozansky polynomial
for $C\!ab(13,2)T(3,2)$ colored by $\om_2$.
Though the theory of colored (and reduced)
$KhR$\~polynomials is not sufficiently developed
(and no formulas are known for iterated knots).
One has:
$$\rr=\{3,2\},\ss=\{2,1\},\ \,
\h_{\{3,2\},\{2,1\}}(\om_2;q,t,a\!=\!0)=$$
%\smallskip
\renewcommand{\baselinestretch}{0.5}
{\small
\(
1+q t+q^2 t+q^3 t+q t^2+2 q^2 t^2+2 q^3 t^2+2 q^4 t^2+q^2 t^3
+3 q^3 t^3+4 q^4 t^3+3 q^5 t^3+q^2 t^4+2 q^3 t^4+5 q^4 t^4
+6 q^5 t^4+4 q^6 t^4+q^3 t^5+3 q^4 t^5+7 q^5 t^5+8 q^6 t^5
+4 q^7 t^5+q^3 t^6+2 q^4 t^6+5 q^5 t^6+10 q^6 t^6+9 q^7 t^6
+2 q^8 t^6+q^4 t^7+3 q^5 t^7+7 q^6 t^7+12 q^7 t^7+7 q^8 t^7
+q^9 t^7+q^4 t^8+2 q^5 t^8+5 q^6 t^8+10 q^7 t^8+14 q^8 t^8
+4 q^9 t^8+q^5 t^9+3 q^6 t^9+7 q^7 t^9+12 q^8 t^9+11 q^9 t^9
+q^{10} t^9+q^5 t^{10}+2 q^6 t^{10}+5 q^7 t^{10}+10 q^8 t^{10}
+13 q^9 t^{10}+6 q^{10} t^{10}+q^6 t^{11}+3 q^7 t^{11}
+7 q^8 t^{11}+12 q^9 t^{11}+10 q^{10} t^{11}+q^{11} t^{11}
+q^6 t^{12}+2 q^7 t^{12}+5 q^8 t^{12}+10 q^9 t^{12}
+13 q^{10} t^{12}+3 q^{11} t^{12}+q^7 t^{13}+3 q^8 t^{13}
+7 q^9 t^{13}+11 q^{10} t^{13}+9 q^{11} t^{13}+q^7 t^{14}
+2 q^8 t^{14}+5 q^9 t^{14}+10 q^{10} t^{14}+11 q^{11} t^{14}
+2 q^{12} t^{14}+q^8 t^{15}+3 q^9 t^{15}+7 q^{10} t^{15}
+10 q^{11} t^{15}+4 q^{12} t^{15}+q^8 t^{16}+2 q^9 t^{16}
+5 q^{10} t^{16}+9 q^{11} t^{16}+8 q^{12} t^{16}+q^9 t^{17}
+3 q^{10} t^{17}+7 q^{11} t^{17}+8 q^{12} t^{17}+2 q^{13} t^{17}
+q^9 t^{18}+2 q^{10} t^{18}+5 q^{11} t^{18}+8 q^{12} t^{18}
+3 q^{13} t^{18}+q^{10} t^{19}+3 q^{11} t^{19}+6 q^{12} t^{19}
+5 q^{13} t^{19}+q^{10} t^{20}+2 q^{11} t^{20}+5 q^{12} t^{20}
+6 q^{13} t^{20}+q^{14} t^{20}+q^{11} t^{21}+3 q^{12} t^{21}
+5 q^{13} t^{21}+q^{14} t^{21}+q^{11} t^{22}+2 q^{12} t^{22}
+4 q^{13} t^{22}+3 q^{14} t^{22}+q^{12} t^{23}+3 q^{13} t^{23}
+3 q^{14} t^{23}+q^{12} t^{24}+2 q^{13} t^{24}+3 q^{14} t^{24}
+q^{13} t^{25}+2 q^{14} t^{25}+q^{15} t^{25}+q^{13} t^{26}
+2 q^{14} t^{26}+q^{15} t^{26}+q^{14} t^{27}+q^{15} t^{27}
+q^{14} t^{28}+q^{15} t^{28}+q^{15} t^{29}+q^{15} t^{30}
+q^{16} t^{32}.
\)
}
\renewcommand{\baselinestretch}{1.2}

The last monomial here corresponds to the last entry
in the table for dimensions and deviations,
namely $\{ 2, 3, 4, 5, 6, 7\}\ :\ dim=0,\ \kap=16.$
\smallskip

{\sf Some other approaches.}
Let us conclude the paper with some comments of informal nature
concerning other methods related to superpolynomials.
In \cite{CDGZ}, the uncolored {\em Alexander polynomials\,}
for algebraic knots %(for $t\!=\!q, a\!=\!-1$ in the DAHA parameters) 
are expressed in terms of the valuation semigroup $\Ga$:
$\h_{\r,\rk=1}(q,t\!=\!q,a\!=\!-1)/(1\!-\!q)=\sum_{\ga\in \Ga}q^\ga$ 
in the DAHA parameters. A major step was the 
Oblomkov-Schende conjecture \cite{ObS}, extended to
{\em colored\,} 
HOMFLY-PT polynomials and proved in \cite{Ma}.
\vskip 0.2cm

Generally,
the connection of the DAHA-superpolynomials
to the HOMFLY-PT polynomials ($t=q$) was established
in \cite{CJ} for torus knots, at the end of
\cite{ChD1} for iterated torus knots using Rosso-Jones
formula (it was done for $A_1$, but this approach can
be extended to $A_n$), in
\cite{MoS} for torus iterated
knots using {\em the Skein}, and in \cite{ChD2} for torus 
iterated links. In 
\cite{MoS}, the Skein was related to the 
{\em Elliptic Hall algebra\,}
for $t=q$ (see \cite{BSch,SV}), which is isomorphic to the spherical 
DAHA.
The refined theory, which is for $t\neq q$, is generally 
that of the 
Khovanov-Rozansky polynomials, but there is an important
physics direction too \cite{DGR, AS, GS, DMMSS}; for
instance, the super-duality in \cite{GS} is due to  
the action of $\C^\ast\times \C^\ast$ in $M_5$ theory.   

The Skein played the key role in \cite{BeS}.
This paper conjecturally relates spherical DAHA with
arbitrary (not only torus iterated) knots; formulas
provided there are different from ours. For instance, the 
relation of the DAHA-Jones polynomials for $C^\vee C_1$ from 
in \cite{CJJ} to those in \cite{BeS} is unclear to us, 
including the case of torus knots $T(2n+1,2)$. Also, our 
DAHA-Jones polynomials for iterated torus knots from
\cite{ChD1} do not match those from \cite{Sam}. At 
the end of \cite{CJJ}, the $C^\vee C_1$ formulas there 
for proper parameters are identified with the 
superpolynomials for $T(2n+1,2)$ colored by arbitrary rows.
We note that there is no formal definition of the Skein with 
$q\neq t$ (the case of an elliptic curve times a segment) by now;  
see \cite{MoS,BeS} and the end of \cite{ChD2}.   

Let us mention
here paper \cite{Mel}, where the relation of the 
Khovanov-Rozansky stable polynomials to DAHA superpolynomials
(conjectured in \cite{CJ}), and to Gorsky's formulas \cite{Gor1,
Gor2} too,  was 
verified for uncolored torus knots using {\em Soergel modules};
this identification requires significant amount of calculations
of combinatorial nature.
Soergel (bi-) modules play an important role in 
\cite{GNR}. Generally, arbitrary (uncolored)
knots can be considered in this approach.

Also, some refined invariants of
knots on Riemann surfaces of genus $2$
where constructed in \cite{ArS} using refined $6j$ symbols. 
Furthermore, there is an interpretation of 
uncolored superpolynomials
via the Hilbert scheme of $\C^2$ in \cite{GoN}, and several papers
and conjectures involve rational DAHA (see e.g.
\cite{GORS}); this is for uncolored torus knots.

The approaches we mentioned are almost always restricted to 
uncolored knots (our paper is devoted to those colored
by columns). The new geometric construction we present here
is a ``refined continuation" of \cite{CDGZ,ObS,Ma}. It is based
on our definition of ``compactified Jacobians in any ranks".
The latter are of obvious independent interest; not much is known
about the varieties of (torsion free) sheaves over singular
curves in ranks greater than $1$.

%\medskip
{\sf Acknowledgements.} We thank
David Kazhdan very much for various discussions and the referee for 
useful suggestions.
The first author thanks RIMS for the invitation and
hospitality.
%\vfill
%\eject

\comment{
Our special thanks go to Vivek Shende for clarifying discussions
on the topics related to \cite{ORS,GSh} and Mikhail
Khovanov for many conversations on his and Rozansky's theory.
We thank Semen Artamonov for our using his unique software for
calculating colored HOMFLYPT polynomials, Peter Samuelson for
sending us his work before it was posted and him and Yuri Berest
for discussions. The first author thanks Andras Szenes for the
talks on the algebraic-geometric aspects of our construction and
the
University of Geneva for the invitation and hospitality. We
thank the referees for their attention to our paper,
thorough reports and important suggestions.

The paper was mainly written at RIMS. The first author thanks
Hiraku Nakajima and RIMS for the invitation and hospitality and
the Simons Foundation;
the second author is grateful for the invitation to the
school \& conference "Geometric Representation Theory"
at RIMS (July 21- August 1, 2014) and generous RIMS' support.
I.D. acknowledges partial support from.
}

\vskip -4cm
%\medskip
\bibliographystyle{unsrt}

%\vfill\eject

%%%%%!!!!! APPENDICES-COMMENT

\appendix
\setcounter{equation}{0}
\section{\sc Dimensions for
\texorpdfstring{{\mathversion{bold}$\Ga=\lan 4,6,15\ran$}}
{v=9}}
\label{app:noaff}
Let us now provide {\em potential dimensions\,} $\partial-d$
for $\ell=0,1$ in the case of $\r=\C[[z^4,z^6+z^9]]$ and
$\Ga=\lan 4,6,15\ran$. Recall that they coincide with
actual dimensions of cells apart from the types $M,N,L$; the table
of non-affine cells must be used here. The
formatting of the table is as follows.

For $0$\~flags ($\ell=0$), there is no $+$ and after $D_0^\dag$
we put  $\partial -d$ and then the deviation $|D_0|$.
The pairs $\{D_0^\dag,D_1^\dag\}$ extending this $D_0^\dag$
are listed next and marked by $+$; no
deviations are needed since $|D_1|=|D_0|+1$. The
ordering is from left to right and then downward.
The symbol ``$na$" means that the flag is non-admissible.
In this table, $\de=9$, $\partial=4\de=36$, the
maximal deviation is $2\de=18$.
\smallskip

{\tiny
\hskip -0.7cm
\begin{longtable}{|l|l|l|}
\caption*
{\small
Potential dimensions for $\r=\C[[z^4,z^6+z^{9}]]$, $\ell=0,1$:}\\
\hline
\{\ \}\ :\ 36,0 \ $\rightsquigarrow$&
+\{ 34\}:\ 37 \ $\rightsquigarrow$&
+\{ 35\}:\ 36 \ $\rightsquigarrow$\\
\{ 34\}\ :\ 35,1 &
+\{ 4\}:\ na &
+\{ 22\}:\ 37\\
+\{ 26\}:\ 36 &
+\{ 34, 35\}:\ 35 &
\{ 35\}\ :\ 34,1\\
+\{ 5\}:\ na &
+\{ 23\}:\ 36 &
+\{ 27\}:\ 35\\
+\{ 34, 35\}:\ 34 &
\{ 4\}\ :\ na,2 &
\{ 5\}\ :\ na,2\\
\{ 22\}\ :\ 35,2 &
+\{ 4, 22\}:\ 37 &
+\{ 22, 26\}:\ 36\\
+\{ 22, 35\}:\ 35 &
\{ 23\}\ :\ 34,2 &
+\{ 5, 23\}:\ 36\\
\end{longtable}
\hspace {0.1cm}
}

Let us provide the full DAHA superpolynomial
for $\r=\C[[z^4,z^6+z^9]]$ and $\rk=2$. It
presumably coincides with the so-called reduced
stable Khovanov-Rozansky polynomial for
$C\!ab(15,2)T(3,2)$ colored by $\om_2$ upon
some change of variables $q,t,a$ and normalization.
Though we must mention that the theory of
reduced and colored $KhR$\~polynomials is not
sufficiently established (even for simple
torus knots).
The formula is:

$$\rr=\{3,2\},\ss=\{2,3\},
\h_{\{3,2\},\{2,3\}}(\om_2;q,t,a)=$$
%\smallskip
\renewcommand{\baselinestretch}{0.5}
{\small
\(
1+q^2 t+q^3 t+q^4 t+q^5 t+q^6 t+q^7 t+q^4 t^2+q^5 t^2+2 q^6 t^2+
2 q^7 t^2+3 q^8 t^2+3 q^9 t^2+3 q^{10} t^2+q^{11} t^2+q^{12} t^2
+q^6 t^3+q^7 t^3+2 q^8 t^3+3 q^9 t^3+4 q^{10} t^3+5 q^{11} t^3
+6 q^{12} t^3+5 q^{13} t^3+3 q^{14} t^3+2 q^{15} t^3+q^8 t^4
+q^9 t^4+2 q^{10} t^4+3 q^{11} t^4+5 q^{12} t^4+6 q^{13} t^4
+8 q^{14} t^4+8 q^{15} t^4+8 q^{16} t^4+4 q^{17} t^4+2 q^{18} t^4
+q^{10} t^5+q^{11} t^5+2 q^{12} t^5+3 q^{13} t^5+5 q^{14} t^5
+7 q^{15} t^5+9 q^{16} t^5+10 q^{17} t^5+11 q^{18} t^5
+9 q^{19} t^5+4 q^{20} t^5+2 q^{21} t^5+q^{12} t^6+q^{13} t^6
+2 q^{14} t^6+3 q^{15} t^6+5 q^{16} t^6+7 q^{17} t^6
+10 q^{18} t^6+11 q^{19} t^6+13 q^{20} t^6+12 q^{21} t^6
+9 q^{22} t^6+3 q^{23} t^6+q^{24} t^6+q^{14} t^7+q^{15} t^7
+2 q^{16} t^7+3 q^{17} t^7+5 q^{18} t^7+7 q^{19} t^7
+10 q^{20} t^7+12 q^{21} t^7+14 q^{22} t^7+14 q^{23} t^7
+10 q^{24} t^7+6 q^{25} t^7+q^{26} t^7+q^{16} t^8+q^{17} t^8
+2 q^{18} t^8+3 q^{19} t^8+5 q^{20} t^8+7 q^{21} t^8
+10 q^{22} t^8+12 q^{23} t^8+15 q^{24} t^8+14 q^{25} t^8
+11 q^{26} t^8+4 q^{27} t^8+q^{28} t^8+q^{18} t^9+q^{19} t^9
+2 q^{20} t^9+3 q^{21} t^9+5 q^{22} t^9+7 q^{23} t^9
+10 q^{24} t^9+12 q^{25} t^9+15 q^{26} t^9+15 q^{27} t^9
+7 q^{28} t^9+2 q^{29} t^9+q^{20} t^{10}+q^{21} t^{10}
+2 q^{22} t^{10}+3 q^{23} t^{10}+5 q^{24} t^{10}+7 q^{25} t^{10}
+10 q^{26} t^{10}+12 q^{27} t^{10}+15 q^{28} t^{10}
+10 q^{29} t^{10}+4 q^{30} t^{10}+q^{22} t^{11}+q^{23} t^{11}
+2 q^{24} t^{11}+3 q^{25} t^{11}+5 q^{26} t^{11}+7 q^{27} t^{11}
+10 q^{28} t^{11}+12 q^{29} t^{11}+10 q^{30} t^{11}
+5 q^{31} t^{11}+q^{24} t^{12}+q^{25} t^{12}+2 q^{26} t^{12}
+3 q^{27} t^{12}+5 q^{28} t^{12}+7 q^{29} t^{12}+10 q^{30} t^{12}
+8 q^{31} t^{12}+4 q^{32} t^{12}+q^{26} t^{13}+q^{27} t^{13}
+2 q^{28} t^{13}+3 q^{29} t^{13}+5 q^{30} t^{13}+7 q^{31} t^{13}
+6 q^{32} t^{13}+3 q^{33} t^{13}+q^{28} t^{14}+q^{29} t^{14}
+2 q^{30} t^{14}+3 q^{31} t^{14}+5 q^{32} t^{14}+4 q^{33} t^{14}
+2 q^{34} t^{14}+q^{30} t^{15}+q^{31} t^{15}+2 q^{32} t^{15}
+3 q^{33} t^{15}+2 q^{34} t^{15}+q^{35} t^{15}+q^{32} t^{16}
+q^{33} t^{16}+2 q^{34} t^{16}+q^{35} t^{16}+q^{34} t^{17}
+q^{35} t^{17}+q^{36} t^{18}
+a^6 \bigl(q^{27}+q^{29} t+q^{30} t
+q^{31} t^2+q^{32} t^2+q^{33} t^2+q^{33} t^3+q^{34} t^3
+q^{35} t^3+q^{36} t^3+q^{35} t^4+q^{36} t^4+q^{37} t^4
+q^{37} t^5+q^{38} t^5+q^{39} t^6\bigr)
+a^5 \bigl(q^{20}+q^{21}
+q^{22}+q^{23}+q^{24}+q^{25}+q^{22} t+2 q^{23} t+3 q^{24} t
+3 q^{25} t+3 q^{26} t+3 q^{27} t+q^{28} t+q^{24} t^2
+2 q^{25} t^2+4 q^{26} t^2+5 q^{27} t^2+5 q^{28} t^2+5 q^{29} t^2
+3 q^{30} t^2+q^{31} t^2+q^{26} t^3+2 q^{27} t^3+4 q^{28} t^3
+6 q^{29} t^3+7 q^{30} t^3+7 q^{31} t^3+5 q^{32} t^3+3 q^{33} t^3
+q^{34} t^3+q^{28} t^4+2 q^{29} t^4+4 q^{30} t^4+6 q^{31} t^4
+8 q^{32} t^4+8 q^{33} t^4+5 q^{34} t^4+3 q^{35} t^4+q^{36} t^4
+q^{30} t^5+2 q^{31} t^5+4 q^{32} t^5+6 q^{33} t^5+7 q^{34} t^5
+7 q^{35} t^5+4 q^{36} t^5+q^{37} t^5+q^{32} t^6+2 q^{33} t^6
+4 q^{34} t^6+6 q^{35} t^6+6 q^{36} t^6+4 q^{37} t^6+q^{38} t^6
+q^{34} t^7+2 q^{35} t^7+4 q^{36} t^7+4 q^{37} t^7+2 q^{38} t^7
+q^{39} t^7+q^{36} t^8+2 q^{37} t^8+2 q^{38} t^8+q^{39} t^8
+q^{38} t^9+q^{39} t^9\bigr)
+a^4 \bigl(q^{14}+q^{15}+2 q^{16}
+2 q^{17}+3 q^{18}+2 q^{19}+2 q^{20}+q^{21}+q^{22}+q^{16} t
+2 q^{17} t+4 q^{18} t+6 q^{19} t+8 q^{20} t+9 q^{21} t
+8 q^{22} t+6 q^{23} t+4 q^{24} t+2 q^{25} t+q^{18} t^2
+2 q^{19} t^2+5 q^{20} t^2+8 q^{21} t^2+13 q^{22} t^2
+15 q^{23} t^2+17 q^{24} t^2+13 q^{25} t^2+10 q^{26} t^2
+5 q^{27} t^2+2 q^{28} t^2+q^{20} t^3+2 q^{21} t^3+5 q^{22} t^3
+9 q^{23} t^3+15 q^{24} t^3+20 q^{25} t^3+23 q^{26} t^3
+22 q^{27} t^3+17 q^{28} t^3+11 q^{29} t^3+5 q^{30} t^3
+2 q^{31} t^3+q^{22} t^4+2 q^{23} t^4+5 q^{24} t^4+9 q^{25} t^4
+16 q^{26} t^4+22 q^{27} t^4+28 q^{28} t^4+27 q^{29} t^4
+24 q^{30} t^4+15 q^{31} t^4+8 q^{32} t^4+3 q^{33} t^4+q^{34} t^4
+q^{24} t^5+2 q^{25} t^5+5 q^{26} t^5+9 q^{27} t^5+16 q^{28} t^5
+23 q^{29} t^5+29 q^{30} t^5+30 q^{31} t^5+24 q^{32} t^5
+15 q^{33} t^5+6 q^{34} t^5+2 q^{35} t^5+q^{26} t^6+2 q^{27} t^6
+5 q^{28} t^6+9 q^{29} t^6+16 q^{30} t^6+23 q^{31} t^6
+30 q^{32} t^6+26 q^{33} t^6+19 q^{34} t^6+9 q^{35} t^6
+3 q^{36} t^6+q^{28} t^7+2 q^{29} t^7+5 q^{30} t^7+9 q^{31} t^7
+16 q^{32} t^7+23 q^{33} t^7+24 q^{34} t^7+18 q^{35} t^7
+9 q^{36} t^7+3 q^{37} t^7+q^{30} t^8+2 q^{31} t^8+5 q^{32} t^8
+9 q^{33} t^8+16 q^{34} t^8+17 q^{35} t^8+14 q^{36} t^8
+6 q^{37} t^8+2 q^{38} t^8+q^{32} t^9+2 q^{33} t^9+5 q^{34} t^9
+9 q^{35} t^9+11 q^{36} t^9+8 q^{37} t^9+3 q^{38} t^9+q^{39} t^9
+q^{34} t^{10}+2 q^{35} t^{10}+5 q^{36} t^{10}+5 q^{37} t^{10}
+4 q^{38} t^{10}+q^{39} t^{10}+q^{36} t^{11}+2 q^{37} t^{11}
+2 q^{38} t^{11}+q^{39} t^{11}+q^{38} t^{12}\bigr)
+a^3 \bigl(q^9+q^{10}+2 q^{11}+3 q^{12}+3 q^{13}+3 q^{14}
+3 q^{15}+2 q^{16}+q^{17}+q^{18}+q^{11} t+2 q^{12} t+4 q^{13} t
+7 q^{14} t+10 q^{15} t+12 q^{16} t+13 q^{17} t+12 q^{18} t
+9 q^{19} t+6 q^{20} t+3 q^{21} t+q^{22} t+q^{13} t^2
+2 q^{14} t^2+5 q^{15} t^2+9 q^{16} t^2+15 q^{17} t^2
+21 q^{18} t^2+26 q^{19} t^2+27 q^{20} t^2+24 q^{21} t^2
+18 q^{22} t^2+10 q^{23} t^2+5 q^{24} t^2+q^{25} t^2+q^{15} t^3
+2 q^{16} t^3+5 q^{17} t^3+10 q^{18} t^3+17 q^{19} t^3
+26 q^{20} t^3+36 q^{21} t^3+41 q^{22} t^3+40 q^{23} t^3
+34 q^{24} t^3+22 q^{25} t^3+12 q^{26} t^3+5 q^{27} t^3
+q^{28} t^3+q^{17} t^4+2 q^{18} t^4+5 q^{19} t^4+10 q^{20} t^4
+18 q^{21} t^4+28 q^{22} t^4+41 q^{23} t^4+51 q^{24} t^4
+54 q^{25} t^4+49 q^{26} t^4+36 q^{27} t^4+22 q^{28} t^4
+10 q^{29} t^4+4 q^{30} t^4+q^{31} t^4+q^{19} t^5+2 q^{20} t^5
+5 q^{21} t^5+10 q^{22} t^5+18 q^{23} t^5+29 q^{24} t^5
+43 q^{25} t^5+56 q^{26} t^5+63 q^{27} t^5+60 q^{28} t^5
+45 q^{29} t^5+27 q^{30} t^5+12 q^{31} t^5+4 q^{32} t^5
+q^{33} t^5+q^{21} t^6+2 q^{22} t^6+5 q^{23} t^6+10 q^{24} t^6
+18 q^{25} t^6+29 q^{26} t^6+44 q^{27} t^6+58 q^{28} t^6
+66 q^{29} t^6+62 q^{30} t^6+43 q^{31} t^6+23 q^{32} t^6
+9 q^{33} t^6+2 q^{34} t^6+q^{23} t^7+2 q^{24} t^7+5 q^{25} t^7
+10 q^{26} t^7+18 q^{27} t^7+29 q^{28} t^7+44 q^{29} t^7
+58 q^{30} t^7+63 q^{31} t^7+52 q^{32} t^7+30 q^{33} t^7
+13 q^{34} t^7+3 q^{35} t^7+q^{25} t^8+2 q^{26} t^8+5 q^{27} t^8
+10 q^{28} t^8+18 q^{29} t^8+29 q^{30} t^8+44 q^{31} t^8
+54 q^{32} t^8+48 q^{33} t^8+31 q^{34} t^8+13 q^{35} t^8
+3 q^{36} t^8+q^{27} t^9+2 q^{28} t^9+5 q^{29} t^9+10 q^{30} t^9
+18 q^{31} t^9+29 q^{32} t^9+40 q^{33} t^9+38 q^{34} t^9
+24 q^{35} t^9+11 q^{36} t^9+2 q^{37} t^9+q^{29} t^{10}
+2 q^{30} t^{10}+5 q^{31} t^{10}+10 q^{32} t^{10}
+18 q^{33} t^{10}+26 q^{34} t^{10}+25 q^{35} t^{10}
+16 q^{36} t^{10}+6 q^{37} t^{10}+q^{38} t^{10}+q^{31} t^{11}
+2 q^{32} t^{11}+5 q^{33} t^{11}+10 q^{34} t^{11}
+15 q^{35} t^{11}+14 q^{36} t^{11}+8 q^{37} t^{11}
+3 q^{38} t^{11}+q^{33} t^{12}+2 q^{34} t^{12}+5 q^{35} t^{12}
+8 q^{36} t^{12}+6 q^{37} t^{12}+3 q^{38} t^{12}+q^{39} t^{12}
+q^{35} t^{13}+2 q^{36} t^{13}+3 q^{37} t^{13}+2 q^{38} t^{13}
+q^{37} t^{14}+q^{38} t^{14}\bigr)
+a^2 \bigl(q^5+q^6+2 q^7+2 q^8+3 q^9+2 q^{10}+2 q^{11}+q^{12}
+q^{13}+q^7 t+2 q^8 t+4 q^9 t+6 q^{10} t+9 q^{11} t+11 q^{12} t
+11 q^{13} t+10 q^{14} t+8 q^{15} t+5 q^{16} t+2 q^{17} t
+q^{18} t+q^9 t^2+2 q^{10} t^2+5 q^{11} t^2+8 q^{12} t^2
+14 q^{13} t^2+19 q^{14} t^2+25 q^{15} t^2+25 q^{16} t^2
+25 q^{17} t^2+18 q^{18} t^2+12 q^{19} t^2+5 q^{20} t^2
+2 q^{21} t^2+q^{11} t^3+2 q^{12} t^3+5 q^{13} t^3+9 q^{14} t^3
+16 q^{15} t^3+24 q^{16} t^3+34 q^{17} t^3+41 q^{18} t^3
+43 q^{19} t^3+38 q^{20} t^3+27 q^{21} t^3+16 q^{22} t^3
+6 q^{23} t^3+2 q^{24} t^3+q^{13} t^4+2 q^{14} t^4+5 q^{15} t^4
+9 q^{16} t^4+17 q^{17} t^4+26 q^{18} t^4+39 q^{19} t^4
+50 q^{20} t^4+60 q^{21} t^4+56 q^{22} t^4+47 q^{23} t^4
+30 q^{24} t^4+16 q^{25} t^4+6 q^{26} t^4+2 q^{27} t^4
+q^{15} t^5+2 q^{16} t^5+5 q^{17} t^5+9 q^{18} t^5+17 q^{19} t^5
+27 q^{20} t^5+41 q^{21} t^5+55 q^{22} t^5+69 q^{23} t^5
+73 q^{24} t^5+63 q^{25} t^5+46 q^{26} t^5+25 q^{27} t^5
+11 q^{28} t^5+3 q^{29} t^5+q^{30} t^5+q^{17} t^6+2 q^{18} t^6
+5 q^{19} t^6+9 q^{20} t^6+17 q^{21} t^6+27 q^{22} t^6
+42 q^{23} t^6+57 q^{24} t^6+74 q^{25} t^6+80 q^{26} t^6
+75 q^{27} t^6+52 q^{28} t^6+29 q^{29} t^6+11 q^{30} t^6
+3 q^{31} t^6+q^{19} t^7+2 q^{20} t^7+5 q^{21} t^7+9 q^{22} t^7
+17 q^{23} t^7+27 q^{24} t^7+42 q^{25} t^7+58 q^{26} t^7
+75 q^{27} t^7+83 q^{28} t^7+72 q^{29} t^7+48 q^{30} t^7
+21 q^{31} t^7+7 q^{32} t^7+q^{33} t^7+q^{21} t^8+2 q^{22} t^8
+5 q^{23} t^8+9 q^{24} t^8+17 q^{25} t^8+27 q^{26} t^8
+42 q^{27} t^8+58 q^{28} t^8+76 q^{29} t^8+75 q^{30} t^8
+58 q^{31} t^8+29 q^{32} t^8+10 q^{33} t^8+q^{34} t^8+q^{23} t^9
+2 q^{24} t^9+5 q^{25} t^9+9 q^{26} t^9+17 q^{27} t^9
+27 q^{28} t^9+42 q^{29} t^9+58 q^{30} t^9+67 q^{31} t^9
+57 q^{32} t^9+31 q^{33} t^9+11 q^{34} t^9+q^{35} t^9
+q^{25} t^{10}+2 q^{26} t^{10}+5 q^{27} t^{10}+9 q^{28} t^{10}
+17 q^{29} t^{10}+27 q^{30} t^{10}+42 q^{31} t^{10}
+50 q^{32} t^{10}+47 q^{33} t^{10}+26 q^{34} t^{10}
+9 q^{35} t^{10}+q^{36} t^{10}+q^{27} t^{11}+2 q^{28} t^{11}
+5 q^{29} t^{11}+9 q^{30} t^{11}+17 q^{31} t^{11}
+27 q^{32} t^{11}+35 q^{33} t^{11}+32 q^{34} t^{11}
+18 q^{35} t^{11}+6 q^{36} t^{11}+q^{29} t^{12}+2 q^{30} t^{12}
+5 q^{31} t^{12}+9 q^{32} t^{12}+17 q^{33} t^{12}
+21 q^{34} t^{12}+20 q^{35} t^{12}+10 q^{36} t^{12}
+3 q^{37} t^{12}+q^{31} t^{13}+2 q^{32} t^{13}+5 q^{33} t^{13}
+9 q^{34} t^{13}+12 q^{35} t^{13}+10 q^{36} t^{13}
+4 q^{37} t^{13}+q^{38} t^{13}+q^{33} t^{14}+2 q^{34} t^{14}
+5 q^{35} t^{14}+5 q^{36} t^{14}+4 q^{37} t^{14}+q^{38} t^{14}
+q^{35} t^{15}+2 q^{36} t^{15}+2 q^{37} t^{15}+q^{38} t^{15}
+q^{37} t^{16}\bigr)
+a \bigl(q^2+q^3+q^4+q^5+q^6+q^7+q^4 t+2 q^5 t+3 q^6 t
+4 q^7 t+5 q^8 t+6 q^9 t+5 q^{10} t+3 q^{11} t+2 q^{12} t
+q^{13} t+q^6 t^2+2 q^7 t^2+4 q^8 t^2+6 q^9 t^2+9 q^{10} t^2
+12 q^{11} t^2+14 q^{12} t^2+13 q^{13} t^2+10 q^{14} t^2
+7 q^{15} t^2+3 q^{16} t^2+q^{17} t^2+q^8 t^3+2 q^9 t^3
+4 q^{10} t^3+7 q^{11} t^3+11 q^{12} t^3+16 q^{13} t^3
+21 q^{14} t^3+24 q^{15} t^3+23 q^{16} t^3+18 q^{17} t^3
+11 q^{18} t^3+5 q^{19} t^3+q^{20} t^3+q^{10} t^4+2 q^{11} t^4
+4 q^{12} t^4+7 q^{13} t^4+12 q^{14} t^4+18 q^{15} t^4
+25 q^{16} t^4+31 q^{17} t^4+35 q^{18} t^4+32 q^{19} t^4
+22 q^{20} t^4+13 q^{21} t^4+5 q^{22} t^4+q^{23} t^4+q^{12} t^5
+2 q^{13} t^5+4 q^{14} t^5+7 q^{15} t^5+12 q^{16} t^5
+19 q^{17} t^5+27 q^{18} t^5+35 q^{19} t^5+42 q^{20} t^5
+44 q^{21} t^5+36 q^{22} t^5+23 q^{23} t^5+11 q^{24} t^5
+4 q^{25} t^5+q^{26} t^5+q^{14} t^6+2 q^{15} t^6+4 q^{16} t^6
+7 q^{17} t^6+12 q^{18} t^6+19 q^{19} t^6+28 q^{20} t^6
+37 q^{21} t^6+46 q^{22} t^6+51 q^{23} t^6+46 q^{24} t^6
+33 q^{25} t^6+17 q^{26} t^6+6 q^{27} t^6+q^{28} t^6+q^{16} t^7
+2 q^{17} t^7+4 q^{18} t^7+7 q^{19} t^7+12 q^{20} t^7
+19 q^{21} t^7+28 q^{22} t^7+38 q^{23} t^7+48 q^{24} t^7
+54 q^{25} t^7+50 q^{26} t^7+36 q^{27} t^7+17 q^{28} t^7
+5 q^{29} t^7+q^{30} t^7+q^{18} t^8+2 q^{19} t^8+4 q^{20} t^8
+7 q^{21} t^8+12 q^{22} t^8+19 q^{23} t^8+28 q^{24} t^8
+38 q^{25} t^8+49 q^{26} t^8+55 q^{27} t^8+47 q^{28} t^8
+29 q^{29} t^8+11 q^{30} t^8+2 q^{31} t^8+q^{20} t^9
+2 q^{21} t^9+4 q^{22} t^9+7 q^{23} t^9+12 q^{24} t^9
+19 q^{25} t^9+28 q^{26} t^9+38 q^{27} t^9+49 q^{28} t^9
+51 q^{29} t^9+36 q^{30} t^9+16 q^{31} t^9+3 q^{32} t^9
+q^{22} t^{10}+2 q^{23} t^{10}+4 q^{24} t^{10}+7 q^{25} t^{10}
+12 q^{26} t^{10}+19 q^{27} t^{10}+28 q^{28} t^{10}
+38 q^{29} t^{10}+44 q^{30} t^{10}+37 q^{31} t^{10}
+18 q^{32} t^{10}+4 q^{33} t^{10}+q^{24} t^{11}+2 q^{25} t^{11}
+4 q^{26} t^{11}+7 q^{27} t^{11}+12 q^{28} t^{11}
+19 q^{29} t^{11}+28 q^{30} t^{11}+34 q^{31} t^{11}
+30 q^{32} t^{11}+16 q^{33} t^{11}+3 q^{34} t^{11}
+q^{26} t^{12}+2 q^{27} t^{12}+4 q^{28} t^{12}+7 q^{29} t^{12}
+12 q^{30} t^{12}+19 q^{31} t^{12}+24 q^{32} t^{12}
+22 q^{33} t^{12}+11 q^{34} t^{12}+2 q^{35} t^{12}
+q^{28} t^{13}+2 q^{29} t^{13}+4 q^{30} t^{13}+7 q^{31} t^{13}
+12 q^{32} t^{13}+16 q^{33} t^{13}+14 q^{34} t^{13}
+7 q^{35} t^{13}+q^{36} t^{13}+q^{30} t^{14}+2 q^{31} t^{14}
+4 q^{32} t^{14}+7 q^{33} t^{14}+9 q^{34} t^{14}
+8 q^{35} t^{14}+3 q^{36} t^{14}+q^{32} t^{15}+2 q^{33} t^{15}
+4 q^{34} t^{15}+5 q^{35} t^{15}+3 q^{36} t^{15}+q^{37} t^{15}
+q^{34} t^{16}+2 q^{35} t^{16}+2 q^{36} t^{16}+q^{37} t^{16}
+q^{36} t^{17}+q^{37} t^{17}\bigr).
\)
}
\renewcommand{\baselinestretch}{1.2}
\smallskip
\end{document}